\documentclass[twoside,12pt]{amsart}
\usepackage{amsthm,amsmath,amssymb,amscd,enumerate, pb-diagram,ifthen}
\newtheorem{theorem}[subsubsection]{Theorem}
\newtheorem{corollary}[subsubsection]{Corollary}
\newtheorem{lemma}[subsubsection]{Lemma}

\newtheorem{assump}[subsubsection]{Assumption}
\newtheorem{property}[subsubsection]{Property}

\newtheorem{proposition}[subsubsection]{Proposition}
\theoremstyle{definition}

\newtheorem{definition}[subsubsection]{Definition}

\theoremstyle{remark}

\theoremstyle{remark}

\newtheorem{remark}[subsubsection]{Remark}

\newtheorem{example}[subsubsection]{Example}
\numberwithin{equation}{subsubsection}

\newtheorem{thm}{Theorem}
\newtheorem{cor}{Corollary}
\newtheorem{prop}{Proposition}
\newtheorem{rem}{Remark}

\newcommand{\bpsi}{\bar{\psi}}
\newcommand{\Mbar}{\overline{\M}}

\newcommand{\com}{\mathbb{C}}
\newcommand{\pibar}{\bar{\pi}}
\newcommand{\fbar}{\bar{f}}

\newcommand{\X}{\mathcal{X}}
\newcommand{\Y}{\mathcal{Y}}
\newcommand{\Z}{\mathcal{Z}}

\newcommand{\K}{\mathcal{K}}
\newcommand{\M}{\mathcal{M}}
\newcommand{\C}{\mathcal{C}}
\newcommand{\A}{\mathcal{A}}
\newcommand{\B}{\mathcal{B}}

\newcommand{\D}{\mathcal{D}}
\newcommand{\F}{\mathcal{F}}
\newcommand{\U}{\mathcal{U}}
\newcommand{\sP}{\mathcal{P}}

\newcommand{\bq}{\mathbf{q}}
\newcommand{\bp}{\mathbf{p}}
\newcommand{\fb}{\mathbf{f}}

\newcommand{\bc}{\mathbf{c}}
\newcommand{\bt}{\mathbf{t}}
\newcommand{\sI}{\mathcal{I}}
\newcommand{\sL}{\mathcal{L}}
\newcommand{\sH}{\mathcal{H}}
\newcommand{\sO}{\mathcal{O}}
\newcommand{\tch}{\widetilde{ch}}
\newcommand{\ttd}{\widetilde{Td}}
\newcommand{\orbcup}{\cdot_{orb}}

\newcommand\ev{\operatorname{ev}}
\newcommand\h{\hbar}
\newcommand\ch{\operatorname{ch}}

\def\<{\left\langle}
\def\>{\right\rangle}

\title{Orbifold Quantum Riemann-Roch, Lefschetz and Serre}
\author{Hsian-Hua Tseng}
\address{Department of Mathematics, 
University of Wisconsin-Madison, 
Van Vleck Hall, 480 Lincoln Drive, 
Madison, WI 53706-1388, USA}
\email{tseng@math.wisc.edu}

\date{\today}
\begin{document}
\begin{abstract} 
Given a vector bundle $F$ on a smooth Deligne-Mumford stack $\X$ and an invertible multiplicative characteristic class $\bc$, we define orbifold Gromov-Witten invariants of $\X$ twisted by $F$ and $\bc$. We prove a ``quantum Riemann-Roch theorem'' (Theorem \ref{qrr}) which expresses the generating function of the twisted invariants in terms of the generating function of the untwisted invariants. A quantum Lefschetz hyperplane theorem is derived from this by specializing to genus zero. As an application, we determine the relationship between genus-0 orbifold Gromov-Witten invariants of $\X$ and that of a complete intersection, under additional assumptions. This provides a way to verify mirror symmetry predictions for some complete intersection orbifolds.
\end{abstract}

\maketitle
\tableofcontents

\section{Introduction}\label{intro}

Our main goal is to extend the Quantum Riemann-Roch theorem of Coates-Givental \cite{CG} in Gromov-Witten theory to the case of algebraic orbifold target spaces (i.e. smooth Deligne-Mumford stacks). As applications, we prove Quantum Serre duality for Deligne-Mumford stacks and a general form of Quantum Lefschetz Hyperplane Theorem for Deligne-Mumford stacks. Results leading towards mirror symmetry of complete intersection orbifolds are also deduced as consequences. 

\subsection{Background: Gromov-Witten Theory of Stacks}
We work over the field of complex numbers $\com$. A {\em Deligne-Mumford stack} $\X$ is a category fibered in groupoids which satisfies several rather complicated conditions. For the precise definition and detailed discussions about properties of Deligne-Mumford stacks, we refer to \cite{LMB} and \cite{Vi}. It is known \cite{KeMo} that a (separated) Deligne-Mumford stack $\X$ has a coarse moduli space $X$ which is in general an algebraic space. For any closed point $x\in X$ there is an \'etale neighborhood $U_x\to X$ of $x$ such that the pullback $U_x\times_X \X$ is a stack of the form $[V_x/\Gamma_x]$ with $V_x$ affine and $\Gamma_x$ a finite group. Thus one may view a Deligne-Mumford stack as a geometric object locally a quotient of an affine scheme by a finite group, just like one would view a scheme as a geometric object locally an affine scheme. This viewpoint is in analogy with the notion of {\em orbifolds} in differential geometry: A complex orbifold is a topological space $X$ together with a choice of  an open neighborhood $U_x \ni x$ for each $x\in X$, an open subset $V_x\subset \com^D$, and a finite group $\Gamma_x$ acting on $V_x$ such that $U_x$ is homeomorphic to a quotient $V_x/\Gamma_x$ of $V_x$ by a finite group $\Gamma_x$ and the collection $\{U_x,V_x,\Gamma_x\}_{x\in X}$ satisfies some compatibility conditions concerning $\Gamma_x$-actions on overlaps. 
 
In this paper we work with Deligne-Mumford stacks, but in view of the analogy mentioned above, the term ``orbifold'' will also be used. By abuse of language, we will treat the terms ``orbifold'' and ``smooth Deligne-Mumford stack'' as synonymous\footnote{We do {\em not} assume that a Deligne-Mumford stack has trivial generic stabilizers, unless otherwise mentioned.}.

Let $\X$ be a smooth Deligne-Mumford stack with projective coarse moduli space $X$.  The {\em inertia stack} $I\X:=\X\times_{\Delta, \X\times\X, \Delta}\X$ associated to $\X$ plays an important role in the theory of stacks. Locally at $x\in \X$, the inertia stack $I\X$ consists of connected components labeled by conjugacy classes of elements $g\in \Gamma_x$. Each connected component is described locally as a quotient $[V_x^g/C_{\Gamma_x}(g)]$, where $V_x^g\subset V_x$ denotes the locus fixed by $g$ and $C_{\Gamma_x}(g)\subset \Gamma_x$ denotes the centralizer of $g$. Objects in the category underlying $I\X$ are pairs $(x,g)$ with $x$ an object in $\X$ and $g\in Aut_\X(x)$. There is a canonical projection $q$ from $I\X$ to $\X$. Also, $I\X$ contains $\X$ as the component corresponding to choosing $g$ to be the identity element in $\Gamma_x$. See Section \ref{orbifold} for more details. 

The construction of Gromov-Witten invariants as intersection numbers on the moduli spaces of stables maps was generalized to symplectic orbifolds by Chen-Ruan \cite{CR2} and to Deligne-Mumford stacks by Abramovich-Graber-Vistoli \cite{AGV}, \cite{AGV2}. A summary of the basics of Gromov-Witten theory for stacks will be given in Section \ref{defi}. The ideas central to their constructions are:  (1) the domain curves $\C$ of a stable map $\C\to\X$ to a stack can be {\em orbicurves}, i.e. they can have nontrivial stack structures at marked points and nodes; (2) the stable maps $\C\to \X$ are required to respect the stack structures of $\C$ and $\X$, i.e. they should be representable morphisms.  

In this paper, we consider a variant of Gromov-Witten theory for stacks. Suppose that $\X$ satisfies Assumption \ref{global_quotient_assumption} below. Given a complex vector bundle $F$ on $\X$ and an invertible multiplicative characteristic class 
$\bc (\cdot )$
of complex vector bundles, we define {\em twisted orbifold Gromov-Witten invariants} using these data. These twisted invariants can be encoded in a generating function, called {\em $(\bc,F)$-twisted total descendant potential} of $\X$, which is defined as follows:
$$\D_{(\bc,F)}(\bt): = \exp{\left(\sum_{g=0}^{\infty}\hbar^{g-1} \sum_{n, d}\frac{Q^d}{n!}\int_{[\Mbar_{g,n}(\X,d)]^{w}} \bc(F_{g,n,d}) \ \wedge_{i=1}^n \sum_{k=0}^{\infty} \ev_i^*(t_k) \bpsi_i^k \right)}.$$
Let us explain the notations in this definition.       
Integration in this formula is performed over the weighted {\em virtual fundamental class} $[\Mbar_{g,n}(\X,d)]^w$ in the moduli space $\Mbar_{g,n}(\X,d)$ of degree-$d$ stable maps to $\X$ from genus-$g$ orbicurves with sections to all $n$ marked gerbes. The cohomology classes $t_k \in H^*(I\X,{\mathbb C})$ for $k=0,1,2,...$, are pulled back to the moduli space by the evaluation maps\footnote{Due to presence of stack structures on the domain curves, the evaluation of stable maps at marked points takes values in $I\X$ (rather then $\X$). Note that $\X$ is a component of $I\X$.} $\ev_i: \Mbar_{g,n}(\X,d) \to I\X$, $i=1,...,n$. The classes $\bpsi_i$ are the first Chern classes of the {\em universal cotangent line}\footnote{These are the cotangent line of the underlying coarse curve, {\em not} the orbicurve. See Section \ref{vir_class_and_descendant} for details.} bundles over the moduli spaces $\Mbar_{g,n}(\X,d)$. The ``twisting factor'' $\bc (F_{g,n,d})$ is the characteristic class $\bc$ applied to the virtual bundle $F_{g,n,d}\in K^0(\Mbar_{g,n}(\X,d))$, which is constructed as follows: 
Consider the universal family of orbifold stable maps,
$$\begin{CD}
\C_{g,n}(\X,d) @> \ev>>  \X \\
@V{f}VV  \\
\Mbar_{g,n}(\X,d).
\end{CD}$$
By definition, $f$ is a family of nodal orbicurves which are the source curves of the orbifold stable maps, and the restrictions of $\ev$ to the fibers give rise the stable maps which $\Mbar_{g,n}(\X,d)$ parametrizes. We put\footnote{It follows from the results of \cite{AGOT} that the map $f$ is a local complete intersection morphism. Therefore the K-theoretic push-forward $Rf_*$ of a bundle has a locally free resolution and thus defines an element in the Grothendieck group $K^0$. See Appendix \ref{vbdle} for more discussion on this.}
$$F_{g,n,d}:= Rf_* (\ev^* F) \in K^0(\Mbar_{g,n}(\X,d)).$$
$Q^d$ is an element in the Novikov ring $\Lambda_{nov}$ (see Section \ref{untwis}), and $\hbar$ is a formal variable. 
Finally, $\D_{(\bc,F)}(\bt)$ depends on $(t_0,t_1,t_2,...)$ and we package them as $\bt(z)=\sum_{k\geq 0}t_kz^k$. (Although we denote it by $\D_{(\bc, F)}(\bt)$, the descendant potential does {\em not} depend on $z$.)

The ``untwisted'' {\em total descendant potential} $\D_\X$ of $\X$, which encodes usual orbifold Gromov-Witten invariants, is defined by the defining equation of $\D_{(\bc,F)}$ with the twisting factors $\bc(F_{g,n,d})$ replaced by $1$. Details of the definition of twisted orbifold Gromov-Witten invariants will be given in Section \ref{twinv}. 

\subsection{Main Result: Orbifold Quantum Riemann-Roch}
The main result of this paper, orbifold Quantum Riemann-Roch theorem, expresses the twisted orbifold Gromov-Witten invariants in terms of the usual invariants. To state the result we need the following quantization formalism introduced into Gromov-Witten theory in \cite{Gi3}. Here we give a brief summary, see Section \ref{giv} for a detailed treatment. Let $H:=H^*(I\X, \com)$ be the cohomology (super-)space of the inertia stack. The space $H$ is equipped with the symmetric inner product (called the orbifold Poincar\'e pairing) $$(a,b)_{orb} := \int_{I\X} a\wedge I^*b,\quad a, b\in H,$$
where $I$ is an involution on $I\X$ induced by the inversion $g \mapsto g^{-1}, \forall g\in Aut_\X(x), x\in\X$. Fix an additive basis $\{\phi_\alpha\}$ of $H$ and let $\{\phi^\alpha\}$ be the dual basis with respect to $(\,,\,)_{orb}$. Introduce the space $\sH := H\otimes \Lambda_s\{z,z^{-1}\}$ of convergent Laurent series in $z$ (see Section \ref{givuntw}). Following \cite{Gi3} and \cite{CG}, we equip $\sH$ with the $\Lambda_s$-valued even symplectic form  
$$\Omega (f,g):= \text{Res}_{z=0}(f(-z),g(z))_{orb}\ dz,\quad f,g\in \sH.$$
The Lagrangian polarization $\sH =\sH_{+}\oplus \sH_{-}$, with $\sH_{+}=H\otimes \Lambda_s\{z\}$ and $\sH_{-}=z^{-1}(H\otimes \Lambda_s\{z^{-1}\})$, identifies $(\sH,\Omega)$ with the cotangent bundle $T^*\sH_{+}$, see Section \ref{givuntw}. Let $p_a^\mu, q_b^\nu$ be Darboux coordinates of $(\sH, \Omega)$ with respect to this polarization, as introduced in Section \ref{givuntw}. Put $p_k:=\sum_\mu p_k^\mu\phi^\mu, q_k=\sum_\nu q_k^\nu\phi_\nu$ and $(\bp:=\sum_{k\geq 0} p_k (-z)^{-k-1}, \bq:=\sum_{k\geq 0}q_kz^k)$.
 
Let $ch_k(\cdot)$ denote the degree $2k$ component of the Chern character. We may view $\bc(\cdot)=\exp(\sum_ks_kch_k(\cdot))$ as a family of characteristic classes depending on variables $s=(s_0, s_1,...)$. As $s$ varies, the twisted descendent potentials $\D_{(\bc,F)}$ define a family $\D_{s}$ of elements in the Fock space\footnote{See Section \ref{givuntw} for definition.} of  formal functions on $\sH_{+}$ using the following convention: For $\bt(z)=t_0+t_1z+t_2z^2+... \in H\otimes\Lambda_s\{z\}$, we identify $\bt$ with the Darboux coordinates $\bq \in H\otimes\Lambda_s\{z\}$ via
$$\bq(z) = \sqrt{\bc (F^{(0)})}\ (\bt(z)-{\bf 1}z),$$
where $F^{(0)}$ is the vector bundle on $I\X$ whose fiber at $(x,g)\in I\X$ is the subspace of $F|_x$ on which $g$ acts with eigenvalue $1$, and ${\bf 1}\in H^*(I\X, \com)$ is the unit cohomology class of the principal component $\X\subset I\X$ (see Section \ref{orbifold}).
Then put
$$\D_{s}(\bq):= \D_{(\bc,F)}(\bt).$$
In other words, $\D_{(\bc,F)}$ is now viewed as a function in $q_0=\sqrt{\bc(F^{(0)})}t_0, q_1=\sqrt{\bc(F^{(0)})}(t_0-{\bf 1})$ and $q_k= \sqrt{\bc(F^{(0)})}t_k, k\geq 2$.
Note that $\D_s|_{s_0=s_1=...=0}=\D_\X$.

We also need to define certain vector bundles on the inertia stack. The inertia stack is a disjoint union $$I\X=\coprod_{i\in \sI} \X_i,$$ where $\sI$ is a index set. For any $(x,g)\in \X_i$ let $r_i$ denote the order of the element $g\in Aut_\X(x)$. To a vector bundle $F$ on $\X$, define $F_i^{(l)}$ over $\X_i$ to be the vector bundle whose fiber $F_i^{(l)}|_{(x,g)}$ at $(x,g)\in \X_i$ is the subspace of $F|_x$ on which $g$ acts with eigenvalue $\exp{(2\pi \sqrt{-1} l/r_i)}$. See Section \ref{vbundle} for more details. Also observe that $H^*(I\X,\com)=\oplus_{i\in \sI}H^*(\X_i,\com)$.

The following is the main result of this paper.

\begin{thm}[Orbifold Quantum Riemann-Roch, see Theorem \ref{qrr}]\label{qrr_orb_intro} 
{\em $$\D_{s}  \approx  \hat{\Delta} \D_{\X}.$$
Here $\Delta: \sH\to \sH$ is the operator given by {\em ordinary} multiplication by 
\[ \Delta = \sqrt{\bc (F^{(0)})} \prod_{i\in \sI} \exp \left(\sum_{0\leq l\leq r_i-1}\sum_{k\geq 0} s_k \sum_{m\geq 0} \frac{B_m(l/r_i)}{m!} \ch_{k+1-m}(F_i^{(l)})z^{m-1}\right),\]
}
and $\hat{\Delta}$ is the differential operator obtained by quantizing $\Delta$.
\end{thm}
\medskip

We now explain the ingredients in the Theorem. 

\begin{enumerate}
\item
The symbol $\approx $ stands for ``equal up to a scalar factor depending on $s$'' which will be explicitly described in Section \ref{oqrr}; see Theorem \ref{qrr}.

\item
Here $B_m(x)$ are the {\em Bernoulli polynomials} defined by $$\frac{te^{tx}}{e^t-1}=\sum_{m\geq 0} \frac{B_m(x)t^m}{m!}.$$ For example, $B_0(x)=1, B_1(x)=x-1/2, B_2(x)=x^2-x+1/6$. 

\item
The operators on $\sH$ defined as multiplication by $\ch_{k+1-m}(F_i^{(l)})z^{m-1}$ over the component $\X_i$ of $I\X$ turns out to be anti-symmetric with respect to the form $\Omega $ and thus define infinitesimal linear symplectic transformations on $\sH$, see Corollary \ref{infsymp}. The quantized operator $\hat{\Delta}$ on the Fock space is defined as follows: The operator $$\log \Delta:=\frac{1}{2}\sum_{k\geq 0}s_k\ch_k(F^{(0)})+\sum_{i\in \sI}\sum_{0\leq l\leq r_i-1}\sum_{k\geq 0} s_k \sum_{m\geq 0} \frac{B_m(l/r_i)}{m!} \ch_{k+1-m}(F_i^{(l)})z^{m-1}$$ is infinitesimally symplectic. We define $\hat{\Delta}:=\exp (\widehat{\log\Delta})$, where $\widehat{\log\Delta}$ is the differential operator defined by quantizing the quadratic Hamiltonians of $\log\Delta$ following the standard rule in Darboux coordinates: 
$$(q_{\alpha}q_{\beta})\hat{\ }:=\h^{-1} q_{\alpha}q_{\beta},\ \  
   (p_{\alpha}p_{\beta})\hat{\ }:=\h \partial_{q_{\alpha}}\partial_{q_{\beta}}, \ \ 
   (q_{\alpha}p_{\beta})\hat{\ }:= q_{\alpha} \partial_{q_{\beta}}.$$ 
See Section \ref{quantization} for more details on the quantization procedure.
\end{enumerate}

\begin{rem}
\hfill
\begin{enumerate}
\renewcommand{\labelenumi}{(\roman{enumi})}
\item
When the target space $\X$ is a {\em manifold}, $\Delta$ is simplified to
$$\exp \left(\sum_{k\geq 0} s_k \sum_{m\geq 0} \frac{B_{2m}(0)}{(2m)!} \ch_{k+1-2m}(F) z^{2m-1}\right).$$
Thus our main Theorem recovers the Quantum Riemann-Roch theorem of Coates-Givental \cite{CG}. Their proof is based on the Grothendieck-Riemann-Roch (GRR) theorem applied to a family of nodal curves and thus goes back to Mumford \cite{M} and Faber-Pandharipande \cite{FaP}. Our proof of Theorem 1 relies on an appropriate generalization, in the spirit of Kawasaki \cite{Ka}, of the GRR formula valid for morphisms between Deligne-Mumford stacks. This version of the GRR formula, explained in Appendix \ref{grr}, is known to hold in algebraic context (it is a result of B. Toen). It is tempting to extend our results to almost K\"ahler orbifolds, but we are unable to do so at this moment. The case of almost K\"ahler manifolds is treated in Appendix B of \cite{C}.

\item
The {\em Bernoulli numbers} $B_{2m}(0)$ arise naturally in the formula of Coates-Givental due to the use of the GRR formula. Peculiarly, the values $B_m(l/r)$ of the Bernoulli polynomials featuring in our main result {\em do not} seem to arise in the generalization of the GRR formula to the case of orbifolds. It would be interesting to have a conceptual understanding of the presence of Bernoulli polynomials in our result.
\end{enumerate}  
\end{rem}
   
\medskip

Theorem \ref{qrr_orb_intro} has some immediate consequences in genus zero. The genus zero $(\bc,F)$-twisted descendant potential is defined as
$$\F_{(\bc,F)}^0:=\sum_{n, d}\frac{Q^d}{n!}\int_{[\Mbar_{0,n}(\X,d)]^{w}} \bc(F_{0,n,d})\wedge_{i=1}^n \sum_{k=0}^{\infty} \ev_i^*(t_k) \bpsi_i^k.$$
It is viewed as a element in the Fock space in the way described above. The genus zero orbifold Gromov-Witten potential $\F_\X^0$ is defined by the above equation with the twisting factor $\bc(F_{0,n,d})$ replaced by $1$. The graphs of the differentials of $\F_{(\bc,F)}^0$ and $\F_\X^0$ are two (formal germs of) Lagrangian submanifolds $\sL_s=\sL_{(\bc,F)}$ and $\sL_\X$ of the symplectic vector space $\sH$. Theorem \ref{qrr_orb_intro} yields a relationship between these two Lagrangian submanifold germs.

\begin{cor}[=Corollary \ref{qrrg0}]
$$\sL_s=\Delta \sL_\X.$$
\end{cor}

The genus $0$ orbifold Gromov-Witten potential $\F_\X^0$ is known to satisfy three sets of partial differential equations:  the string equation (SE), the dilaton equation (DE),  and topological recursion relations (TRR), see Section \ref{univeq}. According to Givental (see \cite{Gi5}, Theorem 1), this is equivalent to the following property of the Lagrangian submanifold germ $\sL_\X$:\\
\\
{\bf Property ($\star$):} $\sL_\X$ is the germ of a Lagrangian cone with the vertex at the origin and such that its tangent spaces $L$ are tangent to $\sL_\X$ exactly along $zL$. See (\ref{overruled}) for its precise meaning.
\\


The property ($\star$) is formulated in terms of the symplectic structure $\Omega$ and the operator of multiplication by $z$. It does not depend on the choice of polarization. Therefore it is invariant under the action of the {\em twisted loop group}, which consists of $End(H^*(I\X))$-valued formal Laurent series $M$ in $z^{-1}$ satisfying\footnote{Here $*$ denotes the adjoint with respect to $(\cdot,\cdot)_{orb}$.} $M^*(-z)M(z)=1$. One checks that $\Delta$ defines an element in the twisted loop group. This yields the following corollary:
\begin{cor}
The Lagrangian submanifold $\sL_s$ satisfies property ($\star$). In other words, twisted orbifold Gromov-Witten invariants in genus zero satisfy the axioms (TRR), (SE), and (DE) of genus zero theory. 
\end{cor}

\subsection{Applications of Quantum Riemann-Roch}
\subsubsection{Quantum Serre duality} 
Consider the orbifold Gromov-Witten theory twisted by the dual vector bundle $F^\vee$ and the dual class $$\bc^\vee (\cdot):=\exp{\left(\sum_{k\geq 0} (-1)^{k+1} s_k ch_k(\cdot)\right)}.$$ Theorem \ref{qrr_orb_intro} implies the following ``Quantum Serre duality''.

\begin{cor}[=Theorem \ref{qserre}]\label{qserreintro}
Let $\bt^\vee(z)=\bc(F)\bt(z)+({\bf 1}-\bc(F))z$. Then we have 
$$\D_{(\bc^\vee,F^\vee)}(\bt^\vee)\approx\D_{(\bc, F)}(\bt).$$\end{cor}

\medskip
See Theorem \ref{qserre} for the precise $s$-dependent scalar factor.

We may equip the bundle $F$ with an $\com^*$-action given by scaling the fibers. We are also interested in the special case of twisting by the equivariant Euler class $e(\cdot)$ with respect to this $\com^*$-action. Let the dual bundle $F^\vee$ be equipped with the dual $\com^*$-action and let $e^{-1}(\cdot)$ be the inverse $\com^*$-equivariant Euler class. 

Let $M: H^*(I\X)\to H^*(I\X)$ be the operator defined as follows: On the cohomology $H^*(\X_i)$ of a component $\X_i\subset I\X$, $M$ is defined to be multiplication by the number $(-1)^{-age(F_i)+\frac{1}{2}rank F_i^{mov}}$. See Section \ref{qserre_for_euler} for more details, including definitions of $age(F_i)$, $(q^*F)^{inv}$, and $F_i^{mov}$. 

Put $\bt^*(z)=z+(-1)^{\frac{1}{2}rank (q^*F)^{inv}}M e(F)(\bt (z)-{\bf 1}z)$ and define a change $$\diamond: Q^d\mapsto Q^d(-1)^{\<c_1(F),d\>}, \quad Q^d\in \Lambda_{nov}, $$ in the Novikov ring. The following Proposition is deduced from Corollary \ref{qserreintro}.

\begin{prop}[=Theorem \ref{eulerserre}]
$$\D_{(e^{-1},F^\vee)}(\bt^*,Q)\approx\D_{(e,F)}(\bt,\diamond Q).$$
\end{prop}
See Theorem \ref{eulerserre} for the precise constant factor.

\subsubsection{Quantum Lefschetz}
Again we consider the $\com^*$-action on the bundle $F$ given by scaling the fibers. Let $\lambda$ denote the equivariant parameter. We now consider the genus zero theory of the special case of twisting by $\com^*$-equivariant Euler class $e$ of this action. We assume that $F$ is pulled back from the coarse moduli space $X$. In this situation, the operator $\Delta$ is closely related to asymptotics of the Gamma function:
$$\Delta \sim \frac{1}{\sqrt{e(F)}}\prod_{i=1}^{N}\frac{1}{\sqrt{2\pi z}}\int_0^\infty e^{\frac{-x+(\lambda+q^*\rho_i) \ln x}{z}}dx,$$ 
where $\rho_1,...,\rho_N$ are Chern roots of $F$.

The intersection of $\sL_\X$ with the affine subspace $-z+z\sH_-$ defines a function $J_\X(t,-z)$ called the $J$-function: For $t\in H^*(I\X)$, define
$$J_\X(t,-z):=-z+t+d_\bq \F_\X^0|_{\bq=t-z},$$
see Definition \ref{jfunction} for more explanation.

\begin{thm}[=Theorem \ref{iandj}]
Let $F$ be a vector bundle which is a direct sum of line bundles pulled back from the coarse moduli space $X$. Let $\rho_1,...,\rho_N$ be the Chern roots of $F$. Let a formal function $I(t,z)$ of $t\in H$ be given as in Definition \ref{hypermod}. Then the family $$t\mapsto I(t,-z), \quad t\in H$$ lies on the cone $\sL_{(e,F)}$. In view of the property ($\star$), the cone $\sL_{(e,F)}$ is determined by this family.
\end{thm}

This is an abstract form of the Quantum Lefschetz Hyperplane Theorem for Deligne-Mumford stacks, which generalizes previous results for varieties (see \cite{Gi2}, \cite{B}, \cite{L}, \cite{Ga}, \cite{CG}, \cite{LLY}). 

\begin{remark}
The formal function $I(t,z)$ may also be written as follows:
$$I(t,z)=\prod_{i=1}^N \frac{\int_0^\infty e^{x/z}J_\X(t+(\lambda+\rho_i)\ln x,z) dx}{\int_0^\infty e^{\frac{x-(\lambda+q^*\rho_i)\ln x}{z}}dx},$$ 
where the integrals are interpreted as their stationary phase asymptotics as $z\to 0$. To see this, first rewrite $J_\X$ in the above equation using the string and divisor equations, then use integration by parts.
\end{remark}

\subsection{Towards a mirror theorem for orbifolds}
Many examples of Calabi-Yau varieties in the mathematics and physics literatures are constructed as complete intersections in toric varieties, and many of them have quotient singularities. In dimension at most three, one may avoid dealing with singular Calabi-Yaus by taking crepant resolutions. In higher dimension, this is not possible in general since crepant resolutions may not exist. Therefore it is desirable to work directly with varieties with quotient singularities. The structure of quotient singularities on varieties is naturally described via Deligne-Mumford stacks.

A motivation to introduce twisted Gromov-Witten invariants is to compute Gromov-Witten invariants of complete intersections and verify predictions from mirror symmetry of Calabi-Yau manifolds (for example quintic threefolds in $\mathbb{P}^4$). This approach first appeared in the work of Kontsevich \cite{Ko}. Ever since the formulation of Quantum Lefschetz hyperplane principle (see e.g. \cite{Gi1-5}, \cite{Ki}, \cite{L}), the verification of mirror symmetry predictions for complete intersections has been divided into two independent parts:
\begin{enumerate}
\item
Compute Gromov-Witten invariants for the ambient spaces;
\item
Understand relationships between Gromov-Witten invariants of the complete intersections and those of the ambient spaces.
\end{enumerate}

One motivation of the present paper is to prove mirror symmetry predictions for orbifolds using this approach. The Quantum Lefschetz hyperplane theorem proved in this paper establishes part (2) for orbifold target spaces under additional assumptions. A more useful version of quantum Lefschetz theorem for orbifolds is proven in \cite{ccit}. So far, works on part (1) have been most successful in the case of toric varieties. The toric mirror construction (see for instance \cite{CK}) applied to a toric orbifold $\X$ yields conjectural mirror pairs of Calabi-Yau orbifolds as complete intersections in toric orbifolds. Under additional convexity assumptions, some twisted orbifold Gromov-Witten invariants are related to orbifold Gromov-Witten invariants of the complete intersections. Thus our Quantum Lefschetz Hyperplane Theorem gives relations between genus-0 orbifold Gromov-Witten invariants of those Calabi-Yau orbifolds and the invariants of the ambient toric orbifolds, see Corollary \ref{comp}. Once the orbifold Gromov-Witten invariants of toric orbifolds are computed (i.e. part (1) is settled), our result yields information about genus-0 orbifold Gromov-Witten invariants of the Calabi-Yau complete intersection orbifolds. This will eventually lead to verifications of mirror symmetry prediction for toric complete intersection orbifolds. Using the results of \cite{ccit}, the case of complete intersections in weighted projective spaces is treated in \cite{CCLT}. We hope to return to other cases in the future.

\subsection{Plan of the paper}
The rest of the paper is organized as follows. Sections \ref{defi} and \ref{giv} contain most of the preparatory materials. In Section \ref{defi} we present some definitions and properties used throughout this paper. Section \ref{orbifold} and \ref{vbundle} contain discussions on important notions of stacks needed in this paper. Properties of orbifold cohomology are reviewed in Section \ref{oprod}. Section \ref{modulistablemap} is devoted to the moduli spaces of orbifold stable maps, on which orbifold Gromov-Witten theory is based. In Section \ref{orbgwtheory} we review the orbifold Gromov-Witten theory constructed in \cite{CR2} and \cite{AGV, AGV2}. We introduce the twisted orbifold Gromov-Witten invariants in Section \ref{twinv}. In Section \ref{giv} we explain how Givental's symplectic vector space formalism \cite {Gi3}, \cite{Gi5} can be applied to twisted and untwisted orbifold Gromov-Witten theory. In Section \ref{oqrr} we state the orbifold Quantum Riemann-Roch theorem (Theorem \ref{qrr}). This is used to derive Quantum Lefschetz Hyperplane Principle in Section \ref{qlh} and \ref{ci}. Orbifold Quantum Serre duality is proved in Section \ref{qse}. Section \ref{pfqrr} contains a proof of Theorem \ref{qrr}. We discuss a Grothendieck-Riemann-Roch formula for Deligne-Mumford stacks in Appendix \ref{grr}. Appendix \ref{vbdle} concerns properties of the virtual bundle $F_{g,n,d}$. Some calculation concerning the quantized operators are given in Appendices \ref{quantized_operator} and \ref{cocycle_calculation}. In Appendix \ref{pf_of_TRR} we present a proof of the topological recursion relation for genus $0$ orbifold Gromov-Witten theory.

\section*{Acknowledgments} The author is deeply grateful to A. Givental for his guidance, constant help and encouragement. Many thanks to D. Abramovich, T. Coates, T. Graber, and H. Iritani for numerous helpful discussions and suggestions on the subject and their interests in this work. Thanks to A. Kresch, Y.-P. Lee, M. Olsson, and B. Toen for many helpful discussions. The author is grateful for the referees' numerous helpful comments and suggestions, which greatly improved the paper. The first version of this paper forms the main part of the author's Ph.D. thesis. During the subsequent revision of this paper, the author was supported in part by postdoctoral fellowships from the Mathematical Science Research Institute (Berkeley, California) and Pacific Institute of Mathematical Sciences (Vancouver, Canada), and a visiting fellowship from Institut Mittag-Leffler (Djursholm, Sweden).

\section{Orbifolds and their Gromov-Witten theory}\label{defi}
In this section, we present some definitions, notations and properties which we use throughout.

\subsection{Orbifolds}\label{orbifold}
Throughout this paper, let $\X$ be a proper smooth Deligne-Mumford stack over the complex numbers $\com$ with projective coarse moduli space $X$. In this section, we discuss some general properties of $\X$ and fix notations throughout.


A friendly introduction to basic notions of stacks can be found in \cite{fan}. For comprehensive introductions to rigorous foundation of stacks the reader may consult \cite{E} and the Appendix of \cite{Vi}. A very detailed treatment of the theory of algebraic stacks can be found in \cite{LMB} (see also the forthcoming book \cite{BEFFGK}). The geometry of a stack of the form $[M/G]$ with $M$ a scheme and $G$ an algebraic group is essentially equivalent to the equivariant geometry of $M$ with respect to the $G$-action. Since almost all stacks we treat in this paper are of this form, keeping this interpretation in mind may help the readers unfamiliar with stacks understand this paper. 

Recall that a morphism $f:\X\to \Y$ of stacks is called {\em representable} if for every morphism $g:S\to \Y$ from a scheme $S$, the fiber product $S\times_{g,\Y,f}\X$ is a scheme. In particular, any morphism from a scheme to a stack is representable.

To a Deligne-Mumford stack $\X$ we can associate a {\em coarse moduli space} $X$ which is in general an algebraic space \cite{KeMo}. For a morphism $\X\to\Y$ of stacks, there is an induced morphism $X\to Y$ between their coarse moduli spaces. 
 
We now introduce the inertia stack associated to a stack $\X$, which plays a central role in Gromov-Witten theory for stacks.
\begin{definition}\label{inertia}
Let $\X$ be a Deligne-Mumford stack. The inertia stack $I\X$ associated to $\X$ is defined to be the fiber product $$I\X:=\X\times_{\Delta,\X\times \X,\Delta}\X, $$ where $\Delta: \X\to \X\times \X$ is the diagonal morphism. The objects in the category underlying $I\X$ can be described as follows:
\begin{equation*}
\begin{split}
Ob(I\X)&=\{(x,g)|x\in Ob(\X), g\in Aut_\X(x)\}\\
&=\{(x,H,g)|x\in Ob(\X), H\subset Aut_\X(x), g \mbox{ a generator of } H\}.
\end{split}
\end{equation*}
\end{definition}

\begin{remark}\label{inerem}
\hfill
\begin{enumerate}
\renewcommand{\labelenumi}{(\roman{enumi})}
\item
For a stack $\X$ over $\com$, $I\X$ is isomorphic to the stack of representable morphisms from a constant cyclotomic gerbe to $\X$,
\begin{equation}\label{inertia=cyclotomic_inertia}
I\X\simeq\coprod_{r\in \mathbb{N}} HomRep(B\mu_r,\X).
\end{equation}
At the level of objects, this means
$$Ob(I\X)=\{(x,H,\chi)|x\in Ob(\X), H\subset Aut(x), \chi : H\to\mu_r \mbox{ an isomorphism for some } r\}.$$
Since we work over $\com$ we will from now identify $\mu_r$ as the subgroup of $\com^*$ of $r$-th roots of $1$, and fix a generator $\mathfrak{u}_{r}:=\exp(2\pi\sqrt{-1}/r)$ of $\mu_r$. In doing so, the identification (\ref{inertia=cyclotomic_inertia}) can be described as follows. An object $(x, g)$ of $I\X$ over a scheme $S$ is identified with a representable morphism $S\times B\mu_r \to \X$ such that the image is $x$ and the induced group homomorphism $\mu_r\to Aut_\X(x)$ takes $\mathfrak{u}_r$ to $g$.

This description of $I\X$ will also be used. For more details, see \cite{AGV}, Section 4.4 and \cite{AGV2}, Section 3.2.
\item
There is a natural projection $q:I\X\to \X$. On objects we have $q((x,g))=x$.
\end{enumerate}
\end{remark}
An important observation is that the inertia stack $I\X$ is in general not connected (unless $\X$ is a connected algebraic space). We write 
$$I\X=\coprod_{i\in \sI}\X_i$$ for the decomposition of $I\X$ into a disjoint union of connected components. Here $\mathcal{I}$ is an index set.
Among all components there is a distinguished one (indexed by $0\in \sI$) $$\X_0:=\{(x,id)|x\in Ob(\X), id\in Aut(x) \text{ is the identity element}\},$$ which is isomorphic to $\X$. 

There is a natural involution $I:I\X\to I\X$ defined by interchanging the factors of $\X\times_{\X\times\X}\X$. On objects we have $I((x,g))=(x,g^{-1})$. The restriction of $I$ to $\X_i$ is denoted by $I_i$. The map $I_i$ is an isomorphism between $\X_i$ and another component which we denote by $\X_{i^I}$. It is clear that $\X_{{i^I}^I}=\X_i$. Also, the restriction of $I$ to the distinguished component $\X_0$ is the identity map $\X_0\to\X_0$. 

There is a locally constant function $ord :I\X\to \mathbb{Z}$ defined by sending $(x,g)$ to the order of $g$ in $Aut_\X(x)$. Let $r_i$ denote its value on the connected component $\X_i$. Note that $r_{i^I}=r_i$.
If we view $I\X$ as in Remark \ref{inerem} (i), it is easy to see that the value of $ord$ at $[B\mu_r\to\X]$ is $r$.


\begin{example}
Let $\X$ be of the form $[M/G]$ with $M$ a smooth variety and $G$ a finite group. We can take the index set $\sI$ to be the set $\{(g)|g\in G\}$ of conjugacy classes of $G$. In this case the centralizer $C_G(g)$ acts on the locus $M^g$ of $g$-fixed points. For the conjugacy class $(g)$ we have the component $$\X_{(g)}=[M^g/C_G(g)],$$ and the distinguished component is $[M^{id}/C_G(id)]=[M/G]$. The morphism $I_{(g)}$ is an isomorphism between $\X_{(g)}$ and $\X_{(g^{-1})}$. In our notation, $(g)^I=(g^{-1})$. Also, the value of the function $ord$ on the component $[M^g/C_G(g)]$ is the order of the element $g$ in $G$.
\end{example}

\subsection{Vector bundles on orbifolds}\label{vbundle}
Let $F$ be a vector bundle on $\X$. When we view $\X$ as a geometric object locally a quotient of an affine scheme by a finite group, we may view $F$ as an object on $\X$ locally an equivariant vector bundle on an affine scheme. In this section we discuss some properties of the pullback bundle $q^*F$, which is a vector bundle on the inertia stack $I\X$.

Denote by $(q^*F)_i$ the restriction to $\X_i$ of the pullback of $F$, i.e. $(q^*F)_i:=q^*F|_{\X_i}$. At a point $(x,g)\in\X_i$, the fiber of $(q^*F)_i$ admits an action of $g$, and is accordingly decomposed into a direct sum of eigenspaces of the $g$-action. This gives a global decomposition (see \cite{T}), $$(q^*F)_i=\bigoplus_{0\leq l< r_i}F_i^{(l)},$$ where $F_i^{(l)}$ is the eigen-subbundle with eigenvalue $\zeta_{r_i}^l$ and $\zeta_{r_i}=\exp(2\pi \sqrt{-1}/r_i)$ is a primitive $r_i$-th root of unity. We make the convention that $0 \leq l < r_i$. Define $(q^*F)_i^{inv}:=F_i^{(0)}$. Denote by $q^*F^{inv}$ the bundle over $I\X$ whose restriction to $\X_i$ is $(q^*F)_i^{inv}$. 

The following result addresses compatibility of the decomposition of $(q^*F)_i$ with pulling back by the involution $I_i: \X_i\to \X_{i^I}$. 
\begin{lemma}\label{bundleinv}
\hfill
\begin{enumerate}
\item \label{a}
$I_i^* (F_{i^I}^{(r_i-l)})=F_i^{(l)}$ for $0<l<r_i$.
\item \label{b}
$I_i^* (F_{i^I}^{(0)})=F_i^{(0)}$.
\end{enumerate}
\end{lemma}
\begin{proof}
We verify (\ref{a}). Let $S$ be a scheme and $(x,g^{-1})$ an $S$-valued point of $\X_{i^I}$. Denote by $\tilde{x}: S\to \X_{i^I}$ the morphism corresponding to $(x,g^{-1})$. Then the $S$-valued point $(x,g)$ of $\X_i$ corresponds to the morphism $I_{i^I}\circ \tilde{x}:S\to \X_i$.

Since $F_{i^I}^{(r_i-l)}|_{(x,g^{-1})}:=\tilde{x}^*(F_{i^I}^{(r_i-l)})$ is the subbundle of  $F_{i^I}|_{(x,g^{-1})}:=\tilde{x}^*F$ on which $g^{-1}$ acts with eigenvalue $\zeta_{r_i}^{r_i-l}$, $g$ acts on $F_{i^I}^{(r_i-l)}|_{(x,g^{-1})}$ with eigenvalue $\zeta_{r_i}^l$. Also, 
$$(I_i^*(F_{i^I}^{(r_i-l)}))|_{(x,g)}:=(I_{i^I}\circ\tilde{x})^*I_i^* (F_{i^I}^{(r_i-l)})=\tilde{x}^*(F_{i^I}^{(r_i-l)})=:F_{i^I}^{(r_i-l)}|_{(x,g^{-1})}.$$ Therefore $(I_i^* (F_{i^I}^{(r_i-l)}))|_{(x,g)}$ is the subbuundle of $\tilde{x}^*F$ on which $g$ acts with eigenvalue $\zeta_{r_i}^l$, which is $F_i^{(l)}|_{(x,g)}:=(I_{i^I}\circ \tilde{x})^*(F_i^{(l)})$.
Hence $I_i^* (F_{i^I}^{(r_i-l)})\subset F_i^{(l)}$. The same argument proves that $I_{i^I}^*(F_i^{(l)})\subset F_{i^I}^{(r_i-l)}$. Since $I_{i^I}\circ I_i$ is the identity map, we find $F_i^{(l)}=I_i^*I_{i^I}^*(F_i^{(l)})\subset I_i^*(F_{i^I}^{(r_i-l)})$. Thus $I_i^* (F_{i^I}^{(r_i-l)})= F_i^{(l)}$.

A similar argument proves (\ref{b}).
\end{proof}

We can describe the vector bundles $F_i^{(l)}$ using the identification (\ref{inertia=cyclotomic_inertia}). Each component $\X_i$ of $I\X$ can be viewed as the moduli stack of representable morphisms from constant $\mu_{r_i}$-gerbes to $\X$. Hence there is a universal family over $\X_i$:
\begin{equation*}
\begin{CD}
\X_i\times B\mu_{r_i} @> \rho >>  \X \\
@V{ }VV  \\
\X_i.
\end{CD}
\end{equation*}
Let $\gamma: \X_i\to \X_i\times B\mu_{r_i}$ be the morphism such that the map $\X_i\to \X_i$ to the first factor is the identity and the map $\X_i\to B\mu_{r_i}$ to the second factor\footnote{In other words, the map $\X_i\to B\mu_{r_i}$ to the second factor is the composition $\X_i\to \text{Spec}\, \com\to B\mu_{r_i}$ where $\X_i\to \text{Spec}\, \com$ is the constant map and $\text{Spec}\, \com\to B\mu_{r_i}$ is the atlas of $B\mu_{r_i}$.} corresponds to the trivial $\mu_{r_i}$-bundle over $\X_i$. The pull-back $\rho^*F$ admits an action of $\mathfrak{u}_{r_i}$. Let $(\rho^*F)^{(l)}$ be the eigen sub-bundle of $\rho^*F$ on which $\mathfrak{u}_{r_i}$ acts with eigenvalue $\zeta_{r_i}^l$. Then\footnote{Note that the bundle $(\rho^*F)^{(l)}$ over $\X_i\times B\mu_{r_i}$ can be viewed as a bundle over $\X_i$ with a $\mu_{r_i}$-action. In this point of view pulling back by $\gamma$ simply forgets the $\mu_{r_i}$-action.} we have 
\begin{equation}\label{eigen-bdle-description}
\gamma^*((\rho^*F)^{(l)})=F_i^{(l)}.
\end{equation}

\subsection{Orbifold Cohomology and Orbifold Cup Product}\label{oprod}
In this section we collect some facts about orbifold cohomology which we will use.
\begin{definition}
Following Chen-Ruan \cite{CR1}, the cohomology $H^*(I\X, \com)$ of the inertia stack is called the orbifold cohomology. 
\end{definition}
\begin{remark}
In general, the cohomology with rational coefficients of a stack can be defined as the (singular) cohomology of a geometric realization of the simplicial scheme associated to this stack. For our purpose we define the cohomology of a Deligne-Mumford stack as the (singular) cohomology of its coarse moduli space. In our setting these two definitions  are equivalent. See \cite{AGV2}, Section 2.2 for a detailed discussion. 
\end{remark}
\subsubsection*{Grading on Orbifold Cohomology}
According to \cite{CR1} (see also \cite{AGV}, \cite{AGV2}), the orbifold cohomology $$H^*(I\X, \com)=\oplus_{i\in \mathcal{I}}H^*(\X_i,\com)$$ is equipped with a grading different from the usual one. This grading is explained below.

\begin{definition}
For each component $\X_i$ of $I\X$, the age $age(\X_i)$ is defined as follows: Let $(x,g)\in \X_i$. The tangent space $T_x\X$ is decomposed into a direct sum $\oplus_{0\leq l<r_i}V^{(l)}$ of eigenspaces according to the $g$-action, where $V^{(l)}$ is the eigenspace with eigenvalue $\zeta_{r_i}^l$, $0\leq l<r_i$, and $\zeta_{r_i}=\exp(2\pi\sqrt{-1}\frac{1}{r_i})$. The age is defined to be $$age(\X_i):=\frac{1}{r_i}\sum_{0\leq l<r_i}l\cdot dim_\com V^{(l)}.$$ It is easy to see that this definition is independent of choices of $(x,g)\in \X_i$.
\end{definition}
The following Lemma follows directly from the definition.
\begin{lemma}[\cite{CR1}, Lemma 3.2.1]\label{agerel}
$$age(\X_i)+age(\X_{i^I})=dim_\com \X-dim_\com \X_i.$$
\end{lemma}
\begin{definition}
The orbifold degree of a class $a\in H^*(\X_i,\com)$ is defined to be
$$orbdeg(a):=deg(a)+2age(\X_i).$$
The orbifold degree gives a grading on $H^*(I\X,\com)$ different from the usual one.
\end{definition}

\subsubsection*{Orbifold Poincar\'e pairing}
Following \cite{CR1}, Section 3.3, the orbifold Poincar\'e pairing 
$$(\,\,\,,\,\,\,)_{orb}: H^*(I\X,\com)\times H^*(I\X,\com)\to \com$$ is defined as follows:
For $a\in H^*(\X_i,\com), b\in H^*(\X_{i^I},\com)$, define 
$$(a,b)_{orb}:=\int_{\X_i}a\wedge I_i^*b,$$
where $\int_{\X_i}$ stands for the pushforward map $H^*(\X_i,\com)\to H^*(\text{Spec }\com,\com)\simeq \com$. For other choices of classes $a,b$ supported on components of $I\X$ the pairing $(a,b)_{orb}$ is defined to be $0$. Obviously this definition extends linearly to a pairing on $H^*(I\X,\com)$.
 
The orbifold Poincar\'e pairing pairs cohomology classes from a component $\X_i$ with classes from the isomorphic component $\X_{i^I}$. The fact that it is a non-degenerate pairing follows from the fact that the usual Poincar\'e pairing on $H^*(\X_i,\com)$ is non-degenerate.

\begin{definition}\label{basis}
In what follows we often fix a homogeneous additive basis $\{\phi_\alpha\}$ of $H^*(I\X,\com)$ such that each $\phi_\alpha$ is supported on one component $\X_i$ of $I\X$. We denote by $\{\phi^\alpha\}$ the dual basis under orbifold Poincar\'e pairing.
\end{definition}

\subsubsection*{Orbifold Cup Product}
On $H^*(I\X,\com)$ there is a product structure, defined in \cite{CR1} and \cite{AGV}, called the orbifold cup product (or Chen-Ruan cup product), which is different from the ordinary cup product on $H^*(I\X,\com)$.
\begin{definition}\label{orb_cup_prod}
For $a,b\in H^*(I\X,\com)$, their orbifold cup product $a\orbcup b$ is defined as follows: For $c\in H^*(I\X,\com)$, $$(a\orbcup b,c)_{orb}:=\<a,b,c\>_{0,3,0},$$ where the right side is defined in Section \ref{untwis}.
\end{definition}
Together with the grading by orbifold degrees, $(H^*(I\X,\com), \orbcup)$ is a graded $\com$-algebra with unit ${\bf 1}\in H^0(\X)$.

In the following special case, we can compare the orbifold cup product with the ordinary cup product of $H^*(I\X,\com)$. 
\begin{lemma}\label{ocup}
For $a\in H^*(\X,\com)$ and $b\in H^*(\X_i,\com)$, the orbifold cup product $a\orbcup b$ is equal to the ordinary product $q^*a\cdot b$ in $H^*(I\X,\com)$.
\end{lemma}
\begin{proof}
Using the identification $$\Mbar_{0,3}(\X,0;i,i^I,0)\simeq \X_i\times B\mu_{r_i}\times B\mu_{r_i}$$ described in Remark \ref{spcase} below, we find that $a\orbcup b\in H^*(\X_i,\com)$. For $c\in H^*(\X_{i^I},\com)$, by Definition \ref{orb_cup_prod} we have 
$$(a\orbcup b,c)_{orb}=\int_{\X_i} q^*a\cdot b\cdot I_i^*c.$$
On the other hand, by definition of the orbifold Poincar\'e pairing we have $$(q^*a\cdot b,c)_{orb}=\int_{\X_i}(q^*a\cdot b)\cdot I_i^*c.$$ We conclude by the non-degeneracy of the pairing $(\,,\,)_{orb}$.
\end{proof}

\subsection{Moduli of Orbifold Stable Maps}\label{modulistablemap}
In this section we discuss some properties of the moduli stacks of orbifold stable maps. We also set up notations used throughout the paper.

Let $\Mbar_{g,n}(\X,d)$ be the moduli stack of $n$-pointed genus $g$ orbifold stable maps to $\X$ of degree $d$ {\em with sections to all gerbes} (see \cite{AGV}, Section 4.5). The stack $\Mbar_{g,n}(\X,d)$ parametrizes the following objects: 
 $$\begin{CD}
(\C,\{\Sigma_i\}) @> \mathfrak{f} >>  \X \\
@V{ }VV  \\
T,
\end{CD}$$ 
 where
\begin{enumerate}
\item $\C/T$ is a prestable genus $g$ balanced nodal orbicurve\footnote{In \cite{AV} and \cite{AGV}, this is called a balanced twisted curve.},
\item for $i=1,...,n$, the substack $\Sigma_i\subset \C$ is an \'etale cyclotomic gerbe over $T$ with a section (hence a trivialization), and
\item $\mathfrak{f}$ is a representable morphism whose induced map between coarse moduli spaces is a $n$-pointed genus $g$ stable map of degree $\mathfrak{f}_*[\C]=d\in \text{Eff}(\X)$. (The object $\text{Eff}(\X)$ is defined in Definition \ref{gw-generating-function} below. See \cite{AGV2}, Section 2.2 for the definition of $\mathfrak{f}_*$.) 
\end{enumerate}
A precise definition of balanced nodal orbicurves can be found in \cite{AV} and \cite{AGV}. The key idea is that an orbicurve is not a curve but a ``stacky'' version of curve: Nontrivial stack structures occur only at marked points or nodes. \'Etale locally near a marked point, an orbicurve over $\text{Spec}\,\com$ is isomorphic to the quotient $[\text{Spec}\,\com[z]/\mu_r]$ for some $r$, where the cyclic group $\mu_r$ acts on $\text{Spec}\,\com[z]$ via $z\mapsto \xi z, \xi\in \mu_r$. \'Etale locally near a node, an orbicurve over $\text{Spec}\,\com$ is isomorphic to the quotient $[\text{Spec}(\com[x,y]/(xy))/\mu_r]$ for some $r$, where $\mu_r$ acts on $\text{Spec}(\com[x,y]/(xy))$ via $x\mapsto \xi x, y\mapsto \xi^{-1} y, \xi\in \mu_r$.

An \'etale cyclotomic gerbe over $T$ with a section is identified (via the trivialization given by the section) with $T\times B\mu_r$, where $B\mu_r\simeq[\text{Spec}\,\com/\mu_r]$ is the classifying stack associated to the finite group $\mu_r$.

\begin{remark}
In \cite{AV} and \cite{AGV, AGV2}, orbifold stable maps are called twisted stable maps. Since the word ``twisted'' is used in a different context in this paper, we use the term ``orbifold stable maps'' instead.
\end{remark}

For each $i=1,...,n$ there is an {\em evaluation map} $ev_i:\Mbar_{g,n}(\X,d)\to I\X$ defined as follows: 
$ev_i$ sends an object $\mathfrak{f}:(\C,\{\Sigma_i\})\to \X$ to $\mathfrak{f}|_{\Sigma_i}:\Sigma_i\to\X$. Since $\mathfrak{f}|_{\Sigma_i}$ is a map from a constant cyclotomic gerbe $T\times B\mu_r$ to $\X$, $\mathfrak{f}|_{\Sigma_i}$ is an object in $I\X$ by the description of $I\X$ in Remark \ref{inerem} (i). We obtain a morphism $ev_i:\Mbar_{g,n}(\X,d)\to I\X$.

The stack $\Mbar_{g,n}(\X,d)$ can be decomposed according to images of the evaluation maps. Define $$\Mbar_{g,n}(\X,d;i_1,...,i_n):=ev_1^{-1}(\X_{i_1})\cap...\cap ev_n^{-1}(\X_{i_n})=\cap_{j=1}^n ev_j^{-1}(\X_{i_j}).$$
We have $$\Mbar_{g,n}(\X,d)=\coprod_{i_1,...,i_n\in\mathcal{I}}\Mbar_{g,n}(\X,d;i_1,...,i_n).$$
This decomposition according to images of the evaluation maps will be important for us: later on in our computations we will need explicit control on stack structures at the marked points.

The universal family over the moduli stack $\Mbar_{g,n}(\X,d)$ also admits a modular description. Let $$\Mbar_{g,n+1}(\X,d)':=ev_{n+1}^{-1}(\X_0)\subset \Mbar_{g,n+1}(\X,d)$$ denote the open-and-closed substack consisting of orbifold stable maps with trivial stack structure on the $(n+1)$-st marked point. According to \cite{AV}, Corollary 9.1.3, there is a morphism 
$$f:\Mbar_{g,n+1}(\X,d)'\to \Mbar_{g,n}(\X,d)$$ which forgets the $(n+1)$-st marked point. Moreover, $f$ exhibits $\Mbar_{g,n+1}(\X,d)'$ as the universal family over $\Mbar_{g,n}(\X,d)$, and $ev_{n+1}:\Mbar_{g,n+1}(\X,d)'\to \X_0\simeq \X$ is the universal orbifold stable map. Similarly, the universal family over the substack $\Mbar_{g,n}(\X,d;i_1,...,i_n)$ is $$\Mbar_{g,n+1}(\X,d;i_1,...,i_n,0)\to\Mbar_{g,n}(\X,d;i_1,...,i_n).$$

\begin{remark}\label{K_gn_Xd}
There is another moduli stack $\K_{g,n}(\X,d)$ studied in \cite{AV}, which parametrizes orbifold stable maps {\em without} trivializing the gerbes. Over $\K_{g,n}(\X,d)$ there are $n$ universal gerbes $\mathfrak{S}_j, 1\leq j\leq n$, corresponding to the marked points, and the fiber product over $\K_{g,n}(\X,d)$ of these gerbes is $\Mbar_{g,n}(\X,d)$, see \cite{AGV}, Section 4.5. We will not use this stack $\K_{g,n}(\X,d)$ to construct orbifold Gromov-Witten theory.
\end{remark}

\begin{remark}\label{spcase}
We discuss briefly the special case $(g,n,d)=(0,3,0)$. According to \cite{AGV}, Section 6.2, the evaluation maps of $\K_{0,3}(\X,0)$ can be taken to have target $I\X$. Moreover, by \cite{Ca}, Lemma 7.7, $\K_{0,3}(\X,0)'$ is isomorphic to $I\X$. See also the proof of \cite{AGV2}, Proposition 8.2.1. Under this isomorphism, $ev_1$ is identified with the identity map $I\X\to I\X$, $ev_2$ is identified with $I:I\X\to I\X$, and $ev_3$ is identified with $q:I\X\to \X$. The space $\Mbar_{0,3}(\X,0;i_1,i_2,0)$ is empty if $i_2\neq i_1^I$. We have an isomorphism $$\Mbar_{0,3}(\X,0;i,i^I,0)\simeq \X_i\times B\mu_{r_i}\times B\mu_{r_i}.$$
\end{remark}

\subsubsection{Marked Points and Nodes}
The marked points define divisors in the universal family $\Mbar_{g,n+1}(\X,d)'$. Let $$\D_j\subset\Mbar_{g,n+1}(\X,d)'$$ be the  $j$-th universal gerbe over $\Mbar_{g,n}(\X,d)$. By definition, $\D_j$ is the pullback to $\Mbar_{g,n}(\X,d)$ of the gerbe $\mathfrak{S}_j$ over $\K_{g,n}(\X,d)$. Since $\Mbar_{g,n}(\X,d)$ is the fiber product of all the $\mathfrak{S}_j$'s, the pullback gerbe $\D_j$ admits a canonical section and is thus trivialized by this section. So for each $1\leq j\leq n$ there is a section $\Mbar_{g,n}(\X,d)\to \Mbar_{g,n+1}(\X,d)'$ corresponding to the $j$-th marked point. The image of this section is $\D_j$.

The identification of $\Mbar_{g,n+1}(\X,d)'$ as the universal family over $\Mbar_{g,n}(\X,d)$ implies that $\D_j$ can be described as a moduli space parametrizing maps $\mathfrak{f}: (\C, \{\Sigma_i\})\to \X$ with the following property: the domain has a distinguished balanced node $\Sigma\subset \C$ separating two parts $\C_0$ and $\C_1$. The marked points $\Sigma_j$ and $\Sigma_{n+1}$ lie on $\C_1$ and the other marked points lie on $\C_0$. $\mathfrak{f}|_{\C_0}: (\C_0, \{\Sigma_i\}_{i\neq j, n+1}, \Sigma)\to \X$ is an $n$-pointed orbifold stable map of genus $g$ and degree $d$, and $\mathfrak{f}|_{\C_1}: (\C_1, \Sigma,  \Sigma_j,\Sigma_{n+1})\to \X$ is a $3$-pointed orbifold stable map of genus $0$ and degree $0$.

Put $$\D_{j,(i_1,...,i_n)}:=\D_j\cap\Mbar_{g,n+1}(\X,d;i_1,...,i_n,0).$$
$\D_{j,(i_1,...,i_n)}$ is the $j$-th universal gerbe over $\Mbar_{g,n}(\X,d; i_1,...,i_n)$, which according to the discussion above is canonically trivialized. We have that $\D_{j,(i_1,...,i_n)}$ is isomorphic to $\Mbar_{g,n}(\X,d;i_1,...,i_n)\times B\mu_{r_{i_j}}$. Under this isomorphism, $f|_{\D_{j,(i_1,...,i_n)}}$ coincides with the projection to the first factor.

Let $\Z\subset \Mbar_{g,n+1}(\X,d)'$ be the locus of nodes in the universal family. $\Z$ is a disjoint union $$\Z=\Z^{irr}\coprod \Z^{red},$$ where $\Z^{irr}$ is the locus of non-separating nodes and $\Z^{red}$ is the locus of separating nodes. $\Z$ is of virtual codimension two in $\Mbar_{g,n+1}(\X,d)'$, and is a cyclotomic gerbe over $f(\Z)$.

There is a locally constant function $ord:\Z\to \mathbb{Z}$ defined by assigning to a node the order of its automorphism group: If a node is locally the quotient $[U/\mu_r]$ where $U$ is the curve $xy=t$ and the cyclic group $\mu_r$ of order $r$ acts via $(x, y) \mapsto (\zeta x, \zeta^{-1}y)$, $\forall \zeta \in \mu_r$, then $ord$ sends this node to the integer $r$.

Let $\Z_r:=ord^{-1}(r)\subset\Z$. Define
$$\Z_{(i_1,...,i_n)}:=\Z\cap \Mbar_{g,n+1}(\X,d;i_1,...,i_n,0),\,\,\,\Z_{r,(i_1,...,i_n)}:=\Z_r\cap \Mbar_{g,n+1}(\X,d;i_1,...,i_n,0).$$ The substacks $\Z_{(i_1,...,i_n)}^{irr}, \Z_{(i_1,...,i_n)}^{red}, \Z_{r,(i_1,...,i_n)}^{irr}, \Z_{r,(i_1,...,i_n)}^{red}\subset \Mbar_{g,n+1}(\X,d;i_1,...,i_n,0)$ are similarly defined.

\subsubsection{Stable maps to the coarse moduli space}
Let $\Mbar_{g,n}(X,d)$ be the moduli stack of $n$-pointed genus $g$ stable maps of degree $d$ to the {\em coarse moduli space} $X$. The universal family over $\Mbar_{g,n}(X,d)$ is $\Mbar_{g,n+1}(X,d)\to\Mbar_{g,n}(X,d)$, see for example \cite{BeMa}, Corollary 4.6. There is a morphism $$\pi_n: \Mbar_{g,n}(\X,d)\to\Mbar_{g,n}(X,d),$$ which sends an orbifold stable map to its induced stable map between coarse moduli spaces, see \cite{AV}, Theorem 1.4.1.

\subsection{Orbifold Gromov-Witten Theory}\label{orbgwtheory}
In this section we describe the Gromov-Witten theory for Deligne-Mumford stacks following \cite{AGV}, which is based on the stacks $\Mbar_{g,n}(\X,d)$. We refer the reader to \cite{Ca} and \cite{AGV2} for a construction of orbifold Gromov-Witten theory based on the stacks $\K_{g,n}(\X,d)$ (see Remark \ref{K_gn_Xd}). Both constructions yield the same orbifold Gromov-Witten invariants. The intersection theory for algebraic stacks needed here can be found in \cite{Vi} and \cite{Kr} (which has already been used to construct Gromov-Witten theory for varieties).

\subsubsection{Virtual fundamental classes and descendants}\label{vir_class_and_descendant}
The stack $\Mbar_{g,n}(\X,d)$ admits a perfect obstruction theory relative to the Artin stack of prestable pointed orbicurves (\cite{AGV}, Section 4.6). This obstruction theory is given by the object $(\mathbf{R}^\bullet f_*ev_{n+1}^*T\X)^\vee$ in the derived category $D(Coh(\Mbar_{g,n}(\X,d)))$. Results in \cite{BeFa} and \cite{Beh} apply to yield a virtual fundamental class $$[\Mbar_{g,n}(\X,d)]^{vir}\in H_*(\Mbar_{g,n}(\X,d),\mathbb{Q}).$$ The virtual fundamental class $[\Mbar_{g,n}(\X,d;i_1,...,i_n)]^{vir}$ may be obtained by restriction. As observed in \cite{AGV}, we need to use a {\em weighted} virtual fundamental class $[\Mbar_{g,n}(\X,d)]^w$ defined as follows: the restriction of $[\Mbar_{g,n}(\X,d)]^w$ to $\Mbar_{g,n}(\X,d;i_1,...,i_n)$, which we denote by $[\Mbar_{g,n}(\X,d;i_1,...,i_n)]^{w}$, is defined by $$[\Mbar_{g,n}(\X,d;i_1,...,i_n)]^{w}:=(\prod_{j=1}^nr_{i_j})[\Mbar_{g,n}(\X,d;i_1,...,i_n)]^{vir}.$$ We refer to \cite{AGV}, Section 4.6 for more details.

We now define the descendant classes. For each $i=1,..., n$, the universal family $\Mbar_{g,n+1}(X,d)\to\Mbar_{g,n}(X,d)$ has a section $$\sigma_i:\Mbar_{g,n}(X,d)\to\Mbar_{g,n+1}(X,d),$$ which corresponds to the $i$-th marked point (note that here we consider the moduli stack of stable maps to the {\em coarse moduli space} $X$). Recall that the $i$-th tautological line bundle is defined to be the pullback of the relative dualizing sheaf $\omega_{\Mbar_{g,n+1}(X,d)/\Mbar_{g,n}(X,d)}$ by $\sigma_i$,
$$L_i:=\sigma_i^*\omega_{\Mbar_{g,n+1}(X,d)/\Mbar_{g,n}(X,d)},$$ see for example \cite{Ma}. Let $\psi_i=c_1(L_i)$ and $$\bpsi_i:=\pi_n^*\psi_i\in H^2(\Mbar_{g,n}(\X,d),\mathbb{Q}).$$ These $\bpsi_i$ are the {\em descendant classes} of $\Mbar_{g,n}(\X,d)$.
Note that our choice of descendant classes differs from those of \cite{CR2} by constants.

\subsubsection{Untwisted Theory}\label{untwis}
We are now ready to define the invariants, following \cite{CR2} and \cite{AGV}. Let $a_j\in H^{p_j}(\X_{i_j},\com), j=1,...,n$ be cohomology classes and $k_1,...,k_n$ nonnegative integers. We define the descendant orbifold Gromov-Witten invariants to be
$$\<a_1\bpsi^{k_1},...,a_n\bpsi^{k_n}\>_{g,n,d}:=\int_{[\Mbar_{g,n}(\X,d; i_1,...,i_n)]^w}(ev_i^*a_1)\bpsi_1^{k_1}...(ev_n^*a_n)\bpsi_n^{k_n}.$$
Here $\int_{[\Mbar_{g,n}(\X,d; i_1,...,i_n)]^w}$ stands for capping with the virtual fundamental class $[\Mbar_{g,n}(\X,d; i_1,...,i_n)]^w$ followed by pushing forward to $\text{Spec }\com$. The symbol $\<\ldots\>_{g,n,d}$ is by definition multi-linear in its entries.

The invariant $\<a_1\bpsi^{k_1},...,a_n\bpsi^{k_n}\>_{g,n,d}$ is zero unless
\begin{equation}\label{vir_dimension}
\frac{1}{2}(orbdeg(a_1)+...+orbdeg(a_n))+k_1+...+k_n=(1-g)(dim_\com \X-3)+n+\int_d c_1(T_\X),
\end{equation}
where $orbdeg(a_j)=p_j+2age(\X_{i_j})$ is the orbifold degree defined in Section \ref{oprod}.
This follows from the formula for virtual dimension of $\Mbar_{g,n}(\X,d; i_1,...,i_n)$, which follows from \cite{AGV2}, Theorem 7.2.1. 
\begin{remark} \label{novikov_and_super}


The cohomology $H^*(I\X, \com)$ is viewed as a super vector space. For simplicity we systematically ignore the signs that may come out. It is straightforward to include the signs in our results (c.f. \cite{CG}).
\end{remark}

We can form generating functions to encode these invariants. 
\begin{definition}\label{gw-generating-function}
Let $\mathbf{t}=\mathbf{t}(z)=t_0+t_1z+t_2z^2+...\in H^*(I\X)[z]$. Define
$$\<\mathbf{t},...,\mathbf{t}\>_{g,n,d}=\<\mathbf{t}(\bpsi),...,\mathbf{t}(\bpsi)\>_{g,n,d}:=\sum_{k_1,...,k_n\geq 0}\<t_{k_1}\bpsi^{k_1},...,t_{k_n}\bpsi^{k_n}\>_{g,n,d}.$$
The {\em total descendant potential} is defined to be $$\D_\X(\bt):=\exp\left(\sum_{g\geq 0} \hbar^{g-1}\F_\X^g(\bt)\right),$$ where
$$\F_\X^g(\bt):=\sum_{n\geq 0,d\in \text{Eff}(\X)}\frac{Q^{d}}{n!}\<\mathbf{t},...,\mathbf{t}\>_{g,n,d}.$$ Here $\hbar$ is a formal variable, and $Q^d$ is an element of the {Novikov ring} $\Lambda_{nov}$ which is a completion of the group ring $\com[\text{Eff}(\X)]$ of the semi-group $\text{Eff}(\X)$ of effective curve classes (i.e. classes in $H_2(\X, \mathbb{Q})$ represented by images of representable maps from complete stacky curves to $\X$). The completion is done with respect to an additive valuation
$$v\left(\sum_{d\in \text{Eff}(\X)}c_dQ^d\right)=\text{min}_{c_d\neq 0}\int_d c_1(L)$$
defined by the ample polarization $L$ of the {\em coarse moduli space} $X$ which we choose once and for all.

$\F_\X^g(\bt)$ is called the genus-$g$ descendant potential. It is regarded as a $\Lambda_{nov}$-valued formal power series in the variables $t_k^\alpha$ where $$t_k=\sum_\alpha t_k^\alpha\phi_\alpha\in H^*(I\X,\com), \quad k\geq 0.$$ 
\end{definition}

\begin{remark}
Following the treatement in \cite{AGV2}, Section 2.2, the homology group $H_2(\X, \mathbb{Q})$ with rational coefficients is defined to be the homology group $H_2(X, \mathbb{Q})$ of the coarse moduli space $X$. For this reason degree of effective curve classes in $\X$ are identified with degrees of effective curve classes in $X$, and we will use the term interchangably.
\end{remark}

\begin{lemma}[c.f. \cite{C}, Lemma 1.3.1]\label{D_X_well_defined}
$\D_\X$ is well-defined as a formal power series in the variables $t_k^\alpha$ taking values in $\Lambda_{nov}[[\hbar,\hbar^{-1}]]$. 
\end{lemma}

\begin{proof}
We follows the argument of \cite{C}, Lemma 1.3.1, which treats the manifold case. 
First of all, the expression 
\begin{equation}\label{partition_function}
\sum_{g\geq 0}\hbar^{g-1}\F_\X^g(\bt)
\end{equation}
 is well-defined as a formal power series in $t_k^\alpha$ taking values in $\Lambda_{nov}[[\hbar,\hbar^{-1}]]$. Given a monomial $\hbar^{g-1}Q^d\prod_{1\leq i\leq n}(t_{k_i}^{\alpha_i})^{j_i}$, we define its degree to be the triple $(g-1, \sum_{1\leq i\leq n}j_i, d)$. The coefficient of a monomial of degree $(a,b,c)$ that occurs in (\ref{partition_function}) is a (non-zero) orbifold Gromov-Witten invariant coming from the moduli space $\Mbar_{a+1,b}(\X, c)$.  One observes that
\begin{enumerate}
\item
Since $\Mbar_{a+1,b}(\X, c)$ is finite dimensional, in each degree only finitely many monomials can occur in (\ref{partition_function}); 
\item
Since $\Mbar_{0,0}(\X, 0)$ and $\Mbar_{1,0}(\X,0)$ are empty, if a monomial of degree $(a,b,0)$ occurs in (\ref{partition_function}), then at least one of $a$ and $b$ is strictly positive. 
\end{enumerate}
Now, a monomial of degree $(a,b,c)$ occurs in $\D_\X=\exp(\sum_{g\geq 0}\hbar^{g-1}\F_\X^g(\bt))$ only if there are monomials of degrees $(a_1,b_1,c_1),..., (a_N, b_N, c_N)$ in (\ref{partition_function}) such that $$a_1+...+a_N=a, \quad b_1+...+b_N=b,\quad c_1+...+c_N=c.$$
By the observations above, there are only finitely many choices of $\{(a_i, b_i,c_i)\}$. The result follows.
\end{proof}

We remark that the orbifold Gromov-Witten theory considered here differs slightly from those in \cite{CR2} and \cite{AGV}: we work with algebraic stacks while \cite{CR2} works with symplectic orbifolds. But unlike \cite{AGV}, we work with cohomology instead of Chow ring. One reason for this is that Poincar\'e duality holds for cohomology, but not for Chow rings in general. A definition of cohomological orbifold Gromov-Witten invariants of Deligne-Mumford stacks using the moduli stack $\K_{g,n}(\X,d)$ can be found in \cite{AGV2}. This definition is equivalent to ours.

\subsubsection{Universal equations in orbifold Gromov-Witten theory}\label{univeq}
The Gromov-Witten invariants for varieties are known to satisfy four sets of universal equations\footnote{
These universal equations can be rewritten as differential equations of the generating functions.}: string equation (SE), divisor equation (DIV), dilaton equation (DE), and topological recursion relations (TRR). One may find the precise forms of these equations in for instance \cite{Ma}. The proof of these equations is based on comparisons of descendant classes on various moduli spaces related by forgetful maps. These four sets of equations hold in orbifold Gromov-Witten theory as well, and they take the same form as those in Gromov-Witten theory for varieties. More precisely, we have\\
{\bf String equation:}
$$
\<a_1\bpsi^{k_1},...,a_n\bpsi^{k_n},{\bf 1}\>_{g,n+1,d}=\sum_{j=1}^n\<a_1\bpsi^{k_1},...,a_j\bpsi^{k_j-1},...,a_n\bpsi^{k_n}\>_{g,n,d};$$
{\bf Divisor equation:} for $\gamma\in H^2(\X, \com)$,
\begin{equation*}
\begin{split}
&\<a_1\bpsi^{k_1},...,a_n\bpsi^{k_n}, \gamma\>_{g, n+1,d}=\left(\int_d\gamma\right)\<a_1\bpsi^{k_1},...,a_n\bpsi^{k_n}\>_{g,n,d}\\
&+\sum_{j=1}^n\<a_1\bpsi^{k_1},...,(\gamma\cdot_{orb} a_j)\bpsi^{k_j-1},...,a_n\bpsi^{k_n}\>_{g,n,d};
\end{split}
\end{equation*}
{\bf Dilaton equation:} 
$$
\<a_1\bpsi^{k_1},...,a_n\bpsi^{k_n}, {\bf 1}\bpsi\>_{g, n+1,d}=(2g-2+n)\<a_1\bpsi^{k_1},...,a_n\bpsi^{k_n}\>_{g,n,d};
$$
{\bf Topological recursion relations (in genus zero):} for $t\in H^*(I\X)$, define 
$$\<\<a_1\bpsi^{k_1},...,a_n\bpsi^{k_n}\>\>_0:=\sum_{k=0}^\infty\sum_{d\in \text{Eff}(\X)}\frac{Q^d}{k!}\<a_1\bpsi^{k_1},...,a_n\bpsi^{k_n}, t,..., t\>_{0,n+k,d}.$$ Then
$$
\<\<\phi_{\alpha_1}\bpsi^{k_1+1},\phi_{\alpha_2}\bpsi^{k_2}, \phi_{\alpha_3}\bpsi^{k_3}\>\>_0=\sum_\alpha\<\<\phi_{\alpha_1}\bpsi^{k_1},\phi_\alpha\>\>_0\<\<\phi^\alpha, \phi_{\alpha_2}\bpsi^{k_2},\phi_{\alpha_3}\bpsi^{k_3}\>\>_0,
$$
where $\{\phi_\alpha\}$ is an additive basis of $H^*(I\X)$ and $\{\phi^\alpha\}$ its dual basis under orbifold Poincar\'e pairing. In these equations the term $\bpsi^{-1}$ is defined to be $0$.

Proofs of (SE), (DIV) and (DE) can be found in \cite{AGV2}. The key observation is that, since our choice of descendant classes are pulled back from moduli space of stable maps to the {\em coarse moduli space}, the comparisons of various descendant classes remain unchanged. See \cite{AGV2} for more details. A proof of (TRR) is given in Appendix \ref{pf_of_TRR}.

\subsubsection{Twisted Theory}\label{twinv}
We now introduce twisted orbifold Gromov-Witten invariants. We will make the following 
\begin{assump}\label{global_quotient_assumption}
$\X$ is a quotient of a smooth quasi-projective scheme by a linear algebraic group.
\end{assump}

Given a vector bundle $F$ over $\X$ and an invertible multiplicative characteristic class $\bc(\cdot)=\exp(\sum_k s_k ch_k(\cdot))$. We define the ``twisting factor'' as follows.
\begin{definition}
For a vector bundle $F$ on $\X$, define $$F_{g,n,d}:=f_*ev_{n+1}^*F,$$ where $f_*$ is the K-theoretic pushforward. Assumption \ref{global_quotient_assumption} and the results of \cite{AGOT} imply that the map $f$ is a local complete intersection morphism. Therefore the K-theoretic push-forward $f_*$ of a bundle has a locally free resolution and thus defines an element in the Grothendieck group $K^0$. Hence $F_{g,n,d}$ is an element in $K^0(\Mbar_{g,n}(\X,d))$. Its restriction to $\Mbar_{g,n}(\X,d;i_1,...,i_n)$, which is an element in $K^0(\Mbar_{g,n}(\X,d;i_1,...,i_n))$, is denoted by $F_{g,n,d,(i_1,..,i_n)}$. 
The cohomology classes $\bc(F_{g,n,d})$ and $\bc(F_{g,n,d,(i_1,..,i_n)})$ are called the twisting factors. 
\end{definition}
More detailed discussions and properties of $F_{g,n,d}$ can be found in Appendix \ref{vbdle}.

We define the $(\bc,F)$-twisted descendant orbifold Gromov-Witten invariants to be 
$$\<a_1\bpsi^{k_1},...,a_n\bpsi^{k_n}\>_{g,n,d,(\bc,F)}:=\int_{[\Mbar_{g,n}(\X,d;i_1,...,i_n)]^w}(ev_i^*a_1)\bpsi_1^{k_1}...(ev_n^*a_n)\bpsi_n^{k_n}\bc(F_{g,n,d,(i_1,..,i_n)}),$$
where $a_1,...,a_n$ are as in Section \ref{untwis}. The symbol $\<\ldots\>_{g,n,d, (\bc, F)}$ is by definition multi-linear in its entries. Again, these invariants can be packaged into generating functions.

\begin{definition}
Define $$\<\mathbf{t},...,\mathbf{t}\>_{g,n,d,(\bc,F)}=\<\mathbf{t}(\bpsi),...,\mathbf{t}(\bpsi)\>_{g,n,d,(\bc,F)}:=\sum_{k_1,...,k_n\geq 0}\<t_{k_1}\bpsi^{k_1},...,t_{k_n}\bpsi^{k_n}\>_{g,n,d,(\bc, F)}.$$
The $(\bc,F)$-{\em twisted total descendant potential} is defined to be $$\D_{(\bc,F)}(\bt):=\exp\left(\sum_{g\geq 0} \hbar^{g-1}\F_{(\bc,F)}^g(\bt)\right),$$
where
$$\F_{(\bc,F)}^g(\bt):=\sum_{n\geq 0,d\in \text{Eff}(\X)}\frac{Q^{d}}{n!}\<\mathbf{t},...,\mathbf{t}\>_{g,n,d,(\bc,F)}.$$ $\F_{(\bc, F)}^g(\bt)$ is regarded as a formal power series in the variables $t_k^\alpha$ taking values in the ring $\Lambda_s$. The ring $\Lambda_s$ is defined to be the completion of $\com[\text{Eff}(\X)][s_0, s_1,...]$ with respect to the additive valuation $$v\left(\sum_{d\in \text{Eff}(\X)}c_dQ^d\right)=\text{min}_{c_d\neq 0}\int_d c_1(L),\quad v(s_k)=k+1 \quad (\text{Here }L \text{ is the chosen ample line bundle on }X).$$ 
\end{definition}
The total descendant potential $\D_{(\bc, F)}$ is well-defined as a formal power series in $t_k^\alpha$ taking values in $\Lambda_s[[\hbar, \hbar^{-1}]]$. The proof of Lemma \ref{D_X_well_defined} can be easily adjusted to treat this case.

\section{Givental's symplectic space formalism}\label{giv}
A. Givental introduces a symplectic vector space formalism to describe Gromov-Witten theory (see \cite{Gi3}, \cite{CG}, \cite{Gi5}). In this formalism many properties of Gromov-Witten invariants can be studied using linear symplectic transformations of a certain symplectic vector space, making them more geometric. In this section we explain how this formalism is applied to orbifold Gromov-Witten theory. We will present this in detail for both twisted and untwisted theories.

To take care of certain convergence issues, we will make use of the following definition.
\begin{definition}[c.f. \cite{ccit}]
Let $R$ be a topological ring with an additive valuation $v: R\setminus \{0\}\to \mathbb{R}$. Define the space of $R$-valued {\em convergent Laurent series} in $z$ to be 
$$R\{z,z^{-1}\}:=\left\{\sum_{n\in \mathbb{Z}}r_nz^n \text{ : } r_n\in R, v(r_n)\to\infty \text{ as } |n|\to \infty\right\}.$$
Note that $R\{z,z^{-1}\}$ is a ring if $R$ is complete. Also put 
\begin{equation*}
\begin{split}
R\{z\}:=\left\{\sum_{n\geq 0}r_nz^n \text{ : } r_n\in R, v(r_n)\to\infty \text{ as } n\to \infty\right\},\\
R\{z^{-1}\}:=\left\{\sum_{n\leq 0}r_nz^n \text{ : } r_n\in R, v(r_n)\to\infty \text{ as } -n\to \infty\right\}.
\end{split}
\end{equation*}
\end{definition}

\subsection{Givental's formalism for untwisted theory}\label{givuntw}

Consider the space $$\sH:=H^*(I\X,\com)\otimes \Lambda_s\{z,z^{-1}\}$$ of orbifold-cohomology-valued convergent Laurent series. There is a $\Lambda_s$-valued symplectic form on $\sH$ given by 
$$\Omega(f,g)=\text{Res}_{z=0}(f(-z),g(z))_{orb}dz, \quad \text{for } f,g\in\sH.$$ 

Consider the following polarization 
\begin{equation}\label{standard_polarization}
\begin{split}
&\sH=\sH_+\oplus\sH_-,\\
&\sH_+:=H^*(I\X,\com)\otimes \Lambda_s\{z\} \text{ and } \sH_-:=z^{-1}H^*(I\X,\com)\otimes\Lambda_s\{z^{-1}\}. 
\end{split}
\end{equation}
This identifies $\sH$ with $\sH_+\oplus \sH_+^\star$, where $\sH_+^\star$ is the dual $\Lambda_s$-module. (We may thus think of $\sH$ as the cotangent bundle $T^*\sH_+$.) Both $\sH_+$ and $\sH_-$ are Lagrangian subspaces with respect to $\Omega$.

Introduce a Darboux coordinate system $\{p_a^\mu,q_b^\nu\}$ on $(\sH,\Omega)$ with respect to the polarization (\ref{standard_polarization}). Namely, in these coordinates, a general point in $\sH$ takes the form 
$$\sum_{a\geq 0}\sum_\mu p_a^\mu\phi^\mu(-z)^{-a-1}+\sum_{b\geq 0}\sum_\nu q_b^\nu\phi_\nu z^b.$$
Put $p_a=\sum_\mu p_a^\mu\phi^\mu$ and $q_b=\sum_\nu q_b^\nu \phi_\nu$. Denote 
\begin{equation*}
\begin{split}
&\bp=\bp(z):=\sum_{k\geq 0} p_k(-z)^{-k-1}=p_0(-z)^{-1}+p_1(-z)^{-2}+...;\\
&\bq=\bq(z):=\sum_{k\geq 0}q_kz^k=q_0+q_1z+q_2z^2+....
\end{split}
\end{equation*}
For $\bt(z)\in\sH_+=H^*(I\X,\com) \otimes\Lambda_s\{z\}$ introduce a shift $\bq(z)=\bt(z)-{\bf 1}z$ called the dilaton shift. 

Define the Fock space $\mbox{\it Fock}$ to be the space of formal functions\footnote{This means formal power series in variables $t_k^\alpha$ where $t_k=\sum_\alpha t_k^\alpha\phi_\alpha$.} in $\bt(z)\in \sH_+$ taking values in $\Lambda_s[[\hbar,\hbar^{-1}]]$. In other words, $\mbox{\it Fock}$ is the space of formal functions on $\sH_+$ in the formal neighborhood of $\bq=-{\bf 1}z$. The descendant potential $\D_\X(\bt)$ is regarded as an element in {\it Fock} via the dilaton shift.

The generating function $\F_\X^0$ of genus-$0$ orbifold Gromov-Witten invariants, which is defined in a formal neighborhood of $-{\bf 1}z$, defines a formal germ of Lagrangian submanifold 
$$\sL_\X:=\{(\bp,\bq)|\bp=d_\bq\F_\X^0\}\subset \sH=T^*\sH_+,$$ which is just the graph of the differential of $\F_\X^0$. Equivalently $\sL_\X$ is defined by all equations of the form $p_a^\mu=\frac{\partial \F_\X^0}{\partial q_a^\mu}$.

By \cite{Gi5}, Theorem 1, string and dilaton equations and topological recursion relations imply that $\sL_\X$ satisfies the following properties. 
\begin{theorem}[c.f. \cite{ccit1}]
$\sL_\X$ is the formal germ of a Lagrangian cone with vertex at the origin such that each tangent space $T$ to the cone is tangent to the cone exactly along $zT$. In other words, if $N$ is a formal neighbourhood in $\sH$ of the unique geometric point on $\sL_\X$, then
\begin{equation}\label{overruled}
\begin{split}
    &\text{(1) } T \cap \sL_\X = zT \cap N ;\\
    &\text{(2) } \text{for each } \fb \in zT \cap N , \text{ the tangent space to } \sL_\X \text{ at } \fb \text{ is } T ;\\
    &\text{(3) } \text{if } T = T_\fb \sL_\X \text{ then } f \in zT\cap N .
\end{split}
\end{equation}
\end{theorem}
The statements in (\ref{overruled}) are valid in the context of formal geometry. So for instance $T\cap \sL_\X=zT\cap N$ means that any formal family of elements of $\sH$ which is both a family of elements of $T$ and of $\sL_\X$ is also a family of elements of both $zT$ and $N$, and vice versa. Also, these statements imply that the tangent spaces $T$ of $\sL_\X$ are closed under multiplication by $z$. Moreover, because $T/zT$ is isomorphic to $H^*(I\X,\com)$, it follows from (\ref{overruled}) that $\sL_\X$ is the union of the (finite-dimensional) family of germs of (infinite-dimensional) linear subspaces 
$$\{zT\cap N | T \text{ is a tangent space of } \sL_\X \}.$$


\begin{definition}\label{jfunction}
Following \cite{Gi2}, we define the $J$-function $J_\X(t,z)$ as follows,
$$J_\X(t,z)=z+t+\sum_{n\geq 1,d\in\text{Eff}(\X)}\frac{Q^d}{(n-1)!}\sum_{k\geq 0,\,\alpha}\<t,...,t,\phi_\alpha\bpsi^k\>_{0,n,d}\frac{\phi^\alpha}{z^{k+1}}.$$  
This is a formal power series in coordinates $t^\alpha$ of $t=\sum_\alpha t^\alpha\phi_\alpha\in H^*(I\X, \com)$ taking values in $\sH$. The point of $\sL_\X$ above $-z+t\in \sH_+$ is $J_\X(t,-z)$. 

For each $k\geq 0$, the coefficient of the $z^{-1-k}$ term in $J_\X(t,z)$ takes values in $H^*(I\X,\com)\otimes \Lambda_s$. According to the decomposition $H^*(I\X, \com)=\oplus_{i\in\sI}H^*(\X_i, \com)$, we write $$J_\X(t,z)=(J_i(t,z))\quad \mbox{where}\quad J_i(t,z) \mbox{ takes values in } H^*(\X_i,\com)\otimes \Lambda_s\{z,z^{-1}\}.$$ We further decompose $J_i$ according to degrees, $$J_i(t,z)=\sum_{d\in\text{Eff}(\X)} J_{i,d}(t,z)Q^d.$$
\end{definition}
This $J$-function plays an important role in the genus-0 theory. For example:
\begin{lemma}\label{J_transverse}
The union of the (finite-dimensional) family $$t\mapsto zT_{J_\X(t,-z)}\sL_\X\cap N, \quad t \text{ in a formal neighborhood of zero in } H^*(I\X, \com)\otimes\Lambda_s,$$ of germs of linear subspaces is $\sL_\X$. 
\end{lemma}
\begin{proof}
According to the discussion above, we just need to prove that every tangent space $T$ of $\sL_\X$ is of the form $T_{J_\X(\tau, -z)}\sL_\X$ for some $\tau\in H^*(I\X, \com)\otimes \Lambda_s$. This can be found in \cite{ccit1}, Proposition 2.16.
\end{proof}

\begin{remark}
In untwisted Gromov-Witten theory one usually use the Novikov ring $\Lambda_{nov}$ as the ground ring. Since we will need to compare untwisted theory with twisted theory, we choose to work with the larger ground ring $\Lambda_s$. This only requires minor notational changes applied to discussions in Section \ref{orbgwtheory}.
\end{remark}

\subsection{Givental's formalism for twisted theory}
The formalism for twisted theory requires a twisted version of the pairing on $H^*(I\X, \com)$ which we call the $(\bc,F)$-twisted orbifold Poincar\'e pairing $(\,\,\,,\,\,\,)_{(\bc,F)}$. It is defined by
$$(a,b)_{(\bc,F)}:=\int_{\X_i}a\wedge I_i^*b\wedge \bc((q^*F)^{inv}_i),\quad \text{for }a\in H^*(\X_i, \com), b\in H^*(\X_{i^I}, \com).$$
For other choices of $a,b$ the pairing $(a,b)_{(\bc,F)}$ is defined to be $0$.

Consider another symplectic vector space $(\sH_{(\bc,F)},\Omega_{(\bc,F)})$, where $\sH_{(\bc,F)}=\sH$ and the $\Lambda_s$-valued symplectic form $\Omega_{(\bc,F)}$ is given by
$$\Omega(f,g)_{(\bc,F)}=\text{Res}_{z=0}(f(-z),g(z))_{(\bc,F)}dz,$$ for $f,g\in\sH_{(\bc, F)}$. 

\begin{lemma}
The symplectic vector spaces $(\sH,\Omega)$ and $(\sH_{(\bc,F)},\Omega_{(\bc,F)})$ are identified via the map
\begin{equation}\label{ide}
(\sH_{(\bc, F)},\Omega_{(\bc,F)})\to(\sH,\Omega)
\end{equation}
defined by $a\mapsto a\sqrt{\bc((q^*F)^{inv})},$ where $a\sqrt{\bc((q^*F)^{inv})}$ is the {\em ordinary} cup product in $H^*(I\X,\com)$.
\end{lemma} 
\begin{proof}
For $a, b\in H^*(I\X)$, we have
\begin{equation*}
\begin{split}
&(a\sqrt{\bc((q^*F)^{inv})},b\sqrt{\bc((q^*F)^{inv})})_{orb}=\int_{I\X}a\sqrt{\bc((q^*F)^{inv})}\wedge I^*(b\sqrt{\bc((q^*F)^{inv})})\\
&=\int_{I\X}a\sqrt{\bc((q^*F)^{inv})}\wedge I^*b\wedge I^*\sqrt{\bc((q^*F)^{inv})}=\int_{I\X}a\sqrt{\bc((q^*F)^{inv})}\wedge I^*b\wedge \sqrt{\bc(I^*((q^*F)^{inv}))}\\
&=\int_{I\X}a\sqrt{\bc((q^*F)^{inv})}\wedge I^*b\wedge \sqrt{\bc((q^*F)^{inv})}=\int_{I\X}a\wedge I^*b\wedge\bc((q^*F)^{inv})=(a,b)_{(\bc,F)}.
\end{split}
\end{equation*}
Here the fact $I^*((q^*F)^{inv})=(q^*F)^{inv}$ is used, see Lemma \ref{bundleinv} (\ref{b}).
\end{proof}

We equip $\sH_{(\bc,F)}$ with the same polarization as that of $\sH$, namely $\sH_{(\bc,F)}=(\sH_{(\bc,F)})_+\oplus (\sH_{(\bc,F)})_-$ with $(\sH_{(\bc,F)})_\pm=\sH_\pm$. This polarization also identifies $\sH_{(\bc,F)}$ with $(\sH_{(\bc, F)})_+\oplus (\sH_{(\bc, F)})_+^\star$, where $(\sH_{(\bc, F)})_+^\star$ is the dual $\Lambda_s$-module. (We may thus think $\sH_{(\bc,F)}$ as the cotangent bundle $T^*(\sH_{(\bc,F)})_+$.) 

We define the twisted dilaton shift to be $\bq(z)=\sqrt{\bc((q^*F)^{inv})}(\bt(z)-{\bf 1}z)$, where the {\em ordinary} cup product in $H^*(I\X,\com)$ is used. Via the twisted dilaton shift the twisted total descendant potential $\D_{(\bc,F)}(\bt)$ is regarded as an element in the Fock space, the space of $\Lambda_s[[\hbar,\hbar^{-1}]]$-valued formal functions on $\sH_+$ in the formal neighborhood of $\bq=-\sqrt{\bc((q^*F)^{inv})}{\bf 1}z$.

Similar to the untwisted case, the twisted genus-$0$ potential $\F^0_{(\bc,F)}$, which is defined in a formal neighborhood of $-{\bf 1}z\in \sH_+$, defines a (formal germ of) Lagrangian submanifold
$$\sL_{(\bc,F)}:=\{(\bp,\bq)|\bp=d_\bq\F_{(\bc,F)}^0\}\subset \sH.$$
Here $\F_{(\bc,F)}^0$ is first regarded as an element in the Fock space of functions on $(\sH_{(\bc,F)})_+\subset (\sH_{(\bc,F)},\Omega_{(\bc,F)})$ via the untwisted dilaton shift. Define a (formal germ of) Lagrangian submanifold $\tilde{\sL}_{(\bc,F)}\subset (\sH_{(\bc,F)},\Omega_{(\bc,F)})$ by the graph of its differential. Second, the map (\ref{ide}) identifies this Lagrangian submanifold with the submanifold $\sL_{(\bc,F)}\subset (\sH,\Omega)$.

We remark that it is not a priori clear whether the Lagrangian submanifold $\sL_{(\bc,F)}$ satisfies (\ref{overruled}) or not. This will be a consequence of our main theorem, see Corollary \ref{qrrg0}.

\begin{definition}\label{twjfunction}
The twisted $J$-function $J_{(\bc,F)}(t,z)$ is defined as follows:
\begin{equation*}
\begin{split}
(J_{(\bc,F)}(t,z),a)_{(\bc,F)}&:=(z+t,a)_{(\bc,F)}+\sum_{n\geq 0,d\in\text{Eff}(\X)}\frac{Q^d}{n!}\<t,...,t,\frac{a}{z-\bpsi}\>_{0,n+1,d,(\bc,F)}\\
&=(z+t,a)_{(\bc,F)}+\sum_{n\geq 0,d\in\text{Eff}(\X)}\sum_{k\geq 0}\frac{Q^d}{n!}\<t,...,t,a\bpsi^k\>_{0,n+1,d,(\bc,F)}\frac{1}{z^{k+1}}.
\end{split}
\end{equation*}
Again, the twisted $J$-function is a formal power series in coordinates $t^\alpha$ of $t=\sum_\alpha t^\alpha\phi_\alpha\in H^*(I\X,\com)$ taking values in $\sH_{(\bc, F)}$.
\end{definition}


\subsection{Quantization of Quadratic Hamiltonians}\label{quantization}
Givental \cite{Gi3} observed that many interesting relations in Gromov-Witten theory be expressed in simple forms by applying the Weyl quantization, which is a standard way to produce (projective) Fock space representations of the Heisenberg Lie algebra, to his symplectic space formalism. In this section, we describe this quantization of quadratic Hamiltonian procedure. This quantization procedure allows us to write the quantum Riemann-Roch formula in a simple form.

Let $A:\sH\to\sH$ be a linear infinitesimally symplectic transformation, i.e. $\Omega (Af, g)+\Omega (f, Ag)=0$ for all $f, g\in \sH$. When $f\in \sH$ is written in Darboux coordinates, the quadratic Hamiltonian $$f\mapsto\frac{1}{2}\Omega(Af,f),$$ is a series of homogeneous degree two monomials in Darboux coordinates $p_a^\alpha,q_b^\alpha$. Define the quantization of quadratic monomials as 
$$\widehat{q_a^\mu q_b^\nu}=\frac{q_a^\mu q_b^\nu}{\hbar},\,\, \widehat{q_a^\mu p_b^\nu}=q_a^\mu\frac{\partial}{\partial q_b^\nu},\,\,\widehat{p_a^\mu p_b^\nu}=\hbar\frac{\partial}{\partial q_a^\mu}\frac{\partial}{\partial q_b^\nu}.$$ 
Extending linearly, this defines a quadratic differential operator $\widehat{A}$, called the quantization of $A$. The differential operators $\widehat{q_aq_b},\widehat{q_ap_b}, \widehat{p_ap_b} $ act on $\mbox{\it Fock}$. Since the quadratic Hamiltonian of $A$ may contain infinitely many monomials, the quantization $\widehat{A}$ do not act on ${\it Fock}$ in general. The quantization of a symplectic transformation of the form $\exp(A)$, with $A$ infinitesimally symplectic, is defined to be $\exp (\widehat{A})=\sum_{k\geq 0}\frac{\widehat{A}^k}{k!}$. In general, $\exp(\widehat{A})$ is not well-defined. However the operator that occurs in our quantum Riemann-Roch formula does act on the descendant potential. 

For infinitesimal symplectomorphisms $A$ and $B$, there is the following relation 
$$[\widehat{A},\widehat{B}]=\{ A,B \}^\wedge+ \C (h_A,h_B),$$ where $\{\cdot,\cdot\}$ is the Lie bracket, $[\cdot,\cdot]$ is the supercommutator, and $h_A$ (respectively $h_B$) is the quadratic Hamiltonian of $A$ (respectively $B$). A direct calculation shows that the cocycle $\C$ is given by
\begin{equation*}
\begin{split}
&\C(p_a^\mu p_b^\nu,q_a^\mu q_b^\nu)=-\C(q_a^\mu q_b^\nu,p_a^\mu p_b^\nu)=1+\delta^{\mu\nu}\delta_{ab},\\
&\C=0\mbox{ on any other pair of quadratic Darboux monomials.}
\end{split}
\end{equation*}
For simplicity, we write $\C(A,B)$ for $\C (h_A,h_B)$.

Some universal equations in orbifold Gromov-Witten theory can be expressed as differential equations satisfied by the total descendant potential $\D_\X$. These differential equations can often be written in very simple forms using the quantization formalism. For example, 
\begin{lemma}
The string equation can be written as 
\begin{equation}\label{string_eq}
\widehat{\left(\frac{1}{z}\right)}\D_\X=0.
\end{equation}
\end{lemma}

\begin{proof}
This is proved in the same way as that for varieties (which can be found in \cite{C}, Example 1.3.3.2). We explain the details for the readers' convenience. 

Put $\bt_i(z)=\sum_{j\geq 0}t_{ij}z^j\in H^*(I\X)[[z]]$.
The string equation in cases $(g,n,d)\neq (0,3,0), (1,1,0)$ can be written as
$$\<\bt_1(\bpsi),...,\bt_{n-1}(\bpsi), {\bf 1}\>_{g,n,d}=\sum_{i=1}^{n-1}\<\bt_1(\bpsi),...,\left[\frac{\bt_i(\bpsi)}{\bpsi}\right]_+ ,... ,\bt_{n-1}(\bpsi)\>_{g,n-1,d},$$
where $$\left[\frac{\bt_i(\bpsi)}{\bpsi}\right]_+=\sum_{j\geq 1}t_{ij}\bpsi^{j-1}.$$
Summing over $g,n,d$ yields 
\begin{equation*}
\begin{split}
&\quad \sum_{g,n,d}\frac{Q^d\hbar^{g-1}}{(n-1)!}\<\bt(\bpsi),...,\bt(\bpsi),{\bf 1}\>_{g,n,d}\\
&=\sum_{g,n,d}\frac{Q^d\hbar^{g-1}}{(n-1)!}\<\left[\frac{\bt(\bpsi)}{\bpsi}\right]_+,\bt(\bpsi),...,\bt(\bpsi)\>_{g,n,d}+\frac{1}{2\hbar}\<\bt(\bpsi),\bt(\bpsi),{\bf 1}\>_{0,3,0}+\<{\bf 1}\>_{0,1,0}.
\end{split}
\end{equation*}
This gives 
$$-\frac{1}{2\hbar}q_0^\alpha g_{\alpha\beta}q_0^\beta-\sum_{g,n,d}\frac{Q^d\hbar^{g-1}}{(n-1)!}\<\left[\frac{\bq(\bpsi)}{\bpsi}\right]_+,\bt(\bpsi),...,\bt(\bpsi)\>_{g,n,d}=0,$$ where $g_{\alpha\beta}=(\phi_\alpha,\phi_\beta)_{orb}$.

A direct calculation shows that this is (\ref{string_eq}).
\end{proof}

In the proof of Theorem \ref{qrr} we will encounter quantizations of operators of the form $A=Bz^m$ with $B\in \text{End}(H^*(I\X))$. An explicit expression of $\widehat{A}$ may be found by a straightforward computation. This is worked out in \cite{C}, Example 1.3.3.1, to which we refer the readers for details. See also Appendix \ref{quantized_operator}.

\subsection{Loop space interpretation}\label{loop_interpret}
In this Section we sketch an interpretation of Givental's formalism in terms of loop spaces. The interpretation is topological in nature, so we work with the topological stack underlying the Deligne-Mumford stack $\X$ (which we still denote by $\X$). We should point out that while this interpretation sheds some light on the conceptual meaning of this formalism, one can work with the formalism without knowing this loop space interpretation. 

Let $L\X=Map(S^1,\X)$ be the stack of loops in $\X$. Definition and properties of $L\X$ can be found in e.g \cite{BGNX}. Loop rotation yields an $S^1$-action on $L\X$.  The stack $L\X^{S^1}$ of $S^1$-fixed loops is identified with the inertia stack $I\X$. 
One may think of $\sH_{+}$ as the $S^1$-equivariant cohomology of $L\X$ expressed in terms of the cohomology of the space $L\X^{S^1}\simeq I\X$ and the first Chern class $z$ of the universal line bundle $L$ over $BS^1$. 

\section{Quantum Riemann-Roch}\label{oqrr}
As in Section \ref{twinv}, consider a characteristic class $\bc$ which is multiplicative and invertible. Since the logarithm of $\bc$ is additive, it is a linear combination of components of the Chern character.  Hence we may write $$\bc(\cdot)=\exp\left(\sum_{k\geq 0} s_k ch_k(\cdot)\right).$$ For convenience, set $s_{-1}=0$. We regard $s_k$ as parameters and consider the twisted descendant potentials $\D_s:=\D_{(\bc,F)}$ as a family of elements in the Fock space of functions on $\sH_+$ depending\footnote{The $s_k$-dependence of $\D_{(\bc,F)}$ occurs in two places: the twisting factor $\bc(F_{g,n,d})$ and the twisted dilaton shift.} on $s_k$. We have $\D_s=\D_\X$ when all $s_k=0$. In this Section we formulate our main result, Theorem \ref{qrr}, which expresses $\D_s$ in terms of $\D_\X$.

\subsection{Some Infinitesimal Symplectic Operators}
In this section we introduce certain operators acting on $\sH$ which will be used in the subsequent sections.

Recall that the Bernoulli polynomials $B_m(x)$ are defined by 
$$\frac{te^{tx}}{e^t-1}=\sum_{m\geq 0}\frac{B_m(x)t^m}{m!},$$ see for instance \cite{WW}, Section 7.20.
In particular we have $B_0(x)=1, B_1(x)=x-1/2$. The Bernoulli numbers $B_m$ are given by $B_m:=B_m(0)$.
The following Lemma is immediate from the definition.
\begin{lemma}\label{ber}
$B_m(1-x)=(-1)^mB_m(x)$.
\end{lemma}

\begin{definition}
For each integer $m\geq 0$, define an element $A_m\in H^*(I\X)=\oplus_{i\in \sI}H^*(\X_i)$ as follows:
The component of $A_m$ on $H^*(\X_i)$ is
$$A_m|_{\X_i}:=\sum_{0\leq l\leq r_i-1} ch(F_i^{(l)})B_m(l/r_i).$$
Let $(A_m)_k$ denote the degree $2k$ part of $A_m$,
$$(A_m)_k|_{\X_i}:=\sum_{0\leq l\leq r_i-1} ch_k(F_i^{(l)})B_m(l/r_i).$$
\end{definition}

Ordinary multiplication by $A_m$ defines an operator on $H^*(I\X)$. By abuse of notation, we denote this operator by $A_m$. The quantization of the operator $A_mz^{m-1}$ will appear in Theorem \ref{qrr}. The main goal of this section is to prove that $A_mz^{m-1}$ is infinitesimally symplectic, which is not a priori clear. It follows from the following result.
\begin{lemma}\label{self}
For $m\geq 1$, the operator $A_{2m+1}$ is anti-self-adjoint with respect to the usual or twisted orbifold Poincar\'e pairing. The operator $A_{2m}$ is self-adjoint with respect to the usual or twisted orbifold Poincar\'e pairing.
\end{lemma}
\begin{proof}
We prove the statements for the usual pairing. The proofs for the twisted pairing are identical.

For $a\in H^*(\X_i), b\in H^*(\X_{i^I})$ and $0<l<r_i$, we have, by Lemma \ref{bundleinv} (\ref{a}), the following:
$$(ch(F_{i}^{(l)})a,b)_{orb}=\int_{\X_{i}}ch(F_{i}^{(l)})a\wedge I^*b
=\int_{\X_{i}}a\wedge I^* ch(F_{i^I}^{(r_i-l)})I^*b=(a,ch(F_{i^I}^{(r_i-l)})b)_{orb}.$$
Multiplying by $B_{2m+1}(l/r_i)$ yields $$(B_{2m+1}(l/r_i)ch(F_{i}^{(l)})a,b)_{orb}=(a,B_{2m+1}(l/r_i)ch(F_{i^I}^{(r_i-l)})b)_{orb}\quad \text{for } 0<l<r_i.$$
By Lemma \ref{ber}, $B_{2m+1}(\frac{l}{r_i})=-B_{2m+1}(1-\frac{l}{r_i})=-B_{2m+1}(\frac{r_i-l}{r_i})$. Hence for $0 <l<r_i$ we have
\begin{equation}\label{asa}
(B_{2m+1}(l/r_i)ch(F_{i}^{(l)})a,b)_{orb}=-(a,B_{2m+1}(\frac{r_i-l}{r_i})ch(F_{i^I}^{(r_i-l)})b)_{orb}.
\end{equation}
By Lemma \ref{bundleinv} (\ref{b}), $(ch(F_{i}^{(0)})a,b)_{orb}=(a,ch(F_{i^I}^{(0)})b)_{orb}$. Since $B_{2m+1}(0)=0$ for $m\geq 1$, we have
\begin{equation}\label{asa0}
(B_{2m+1}(0)ch(F_{i}^{(0)})a,b)_{orb}=-(a,B_{2m+1}(0)ch(F_{i^I}^{(0)})b)_{orb}.
\end{equation}
Adding (\ref{asa}) for $l=1,...,r_i-1$ and (\ref{asa0}) yields
$$(A_{2m+1}|_{\X_{i}}a,b)_{orb}=-(a,A_{2m+1}|_{\X_{i^I}}b)_{orb},$$ which proves the statement about $A_{2m+1}$.

To prove the statement for $A_{2m}$, we start with
$$(B_{2m}(l/r_i)ch(F_{i}^{(l)})a,b)_{orb}=(a,B_{2m}(l/r_i)ch(F_{i^I}^{(r_i-l)})b)_{orb}\quad \text{for } 0<l<r_i,$$ and 
\begin{equation}\label{sa0}
(B_{2m}(0)ch(F_{i}^{(0)})a,b)_{orb}=(a,B_{2m}(0)ch(F_{i^I}^{(0)})b)_{orb}.
\end{equation}
By Lemma \ref{ber}, $B_{2m}(l/r_i)=B_{2m}(1-l/r_i)=B_{2m}(\frac{r_i-l}{r_i})$.
This implies that, for $0<l<r_i$,
\begin{equation}\label{sa}
(B_{2m}(l/r_i)ch(F_{i}^{(l)})a,b)_{orb}=(a,B_{2m}(\frac{r_i-l}{r_i})ch(F_{i^I}^{(r_i-l)})b)_{orb}.
\end{equation}
Adding (\ref{sa}) for $l=1,...,r_i-1$ and (\ref{sa0}) yields 
$$(A_{2m}|_{\X_{i}}a,b)_{orb}=(a,A_{2m}|_{\X_{i^I}}b)_{orb},$$ which proves the statement about $A_{2m}$.
\end{proof}

\begin{remark}
\hfill
\begin{enumerate}
\renewcommand{\labelenumi}{(\roman{enumi})}
\item Since $B_0(x)=1$, we have 
$$A_0|_{\X_{i}}=\sum_{0\leq l\leq r_i-1} ch(F_{i}^{(l)})=ch(q^*F)|_{\X_{i}}.$$ Thus multiplication by $A_0$ defines a self-adjoint operator with respect to both pairings.

\item The operator $A_1$ is not anti-self-adjoint. However note that $$A_1|_{\X_{i}}=\sum_{0\leq l\leq r_i-1}B_1(l/r_i) ch(F_{i}^{(l)})=B_1(0)ch(F_{i}^{(0)})+\sum_{l=1}^{r_i-1}B_1(l/r_i)ch(F_{i}^{(l)}).$$ We can use the arguments in the proof of Lemma \ref{self} to show that 
$$(\sum_{l=1}^{r_i-1}B_1(l/r_i)ch(F_{i}^{(l)})a,b)_{orb}=-(a,\sum_{l=1}^{r_i-1}B_1(1-l/r_i)ch(F_{i^I}^{(r_i-l)})b)_{orb}.$$ Using $B_1(0)=-1/2$ we rewrite this as
$$((A_1|_{\X_{i}}+\frac{1}{2}ch(F_{i}^{(0)}))a,b)_{orb}=-(a,(A_1|_{\X_{i^I}}+\frac{1}{2}ch(F_{i^I}^{(0)}))b)_{orb}.$$
In our notation, $ch(F_{i}^{(0)})=ch((q^*F)^{inv})|_{\X_{i}}$. Thus multiplication by $A_1+\frac{1}{2}ch((q^*F)^{inv})$ defines an anti-self-adjoint operator.

\item Lemma \ref{self} also holds if we replace $A_{2m+1}$ and $A_{2m}$ by $(A_{2m+1})_k$ and $(A_{2m})_k$ respectively.
\end{enumerate}
\end{remark}

\begin{corollary}\label{infsymp}
Multiplications by the following classes define infinitesimally symplectic transformations
on $(\sH,\Omega)$ and $(\sH_{(\bc,F)},\Omega_{(\bc,F)})$:
\begin{equation*}
\begin{split}
&A_{2m}z^{2m-1},\,\, A_{2m+1}z^{2m},\,\,\,m\geq 1;\,\, A_0/z,\,\, A_1+\frac{1}{2}ch((q^*F)^{inv});\\
&(A_{2m})_kz^{2m-1},\,\, (A_{2m+1})_kz^{2m},\,\,\,m\geq 1;\,\, (A_0)_k/z,\,\, (A_1)_k+\frac{1}{2}ch_k((q^*F)^{inv}).
\end{split}
\end{equation*}
\end{corollary}

\subsection{Orbifold Quantum Riemann-Roch formula}
Recall that in the definition of twisted orbifold Gromov-Witten invariants  in Section \ref{twinv}, we assume Assumption \ref{global_quotient_assumption} (i.e $\X$ is assumed to be a quotient of a smooth quasi-projective scheme by a linear algebraic group). This assumption is needed also for the application of Grothendieck-Riemann-Roch. To apply Grothendieck-Riemann-Roch formula for Deligne-Mumford stacks to the universal family of orbifold stable maps, we need the universal family to have certain properties. The required properties are proved in \cite{AGOT} for those $\X$ that satisfy Assumption \ref{global_quotient_assumption}. Many interesting stacks, for instance the toric Deligne-Mumford stacks \cite{BCS}, satisfy Assumption \ref{global_quotient_assumption}. 
Through the collective efforts of many works, including \cite{EHKV}, \cite{KV}, \cite{dJ}, it is now known that if $\X$ is a smooth, separated, generically tame Deligne-Mumford stack over $\com$ with quasi-projective coarse moduli space, then $\X$ satisfies Assumption \ref{global_quotient_assumption}. See \cite{kr_geom_dm}, Section 4 for a detailed account.

Now we state the orbifold quantum Riemann-Roch theorem. Its proof is deferred to Section \ref{pfqrr}.
\begin{theorem}[Orbifold Quantum Riemann-Roch]\label{qrr}
Let $\X$ be as in Assumption \ref{global_quotient_assumption}. 
Then we have
\begin{equation*}
\begin{split} 
&\quad \exp\left(-\frac{s_0}{2}\text{rank}\,F\<\bpsi\>_{1,1,0}+s_0\<c_1(F)\>_{1,1,0}\right)\D_s\\
&=\exp\left(\sum_{k\geq 0} s_k\left(\sum_{m>0} \frac{(A_m)_{k+1-m}z^{m-1}}{m!}+\frac{ch_k((q^* F)^{inv})}{2}\right)^\wedge\right)\exp\left(\sum_{k\geq 0} s_k\left(\frac{(A_0)_{k+1}}{z}\right)^\wedge\right)\D_\X.
\end{split}
\end{equation*}
\end{theorem}

This theorem expresses, in a rather nontrivial way, the twisted descendant potential $\D_s$ in terms of the untwisted potential $\D_\X$.

\begin{remark}
The right-hand side of Theorem \ref{qrr} is well-defined. The verification of this is a straightforward modification of \cite{C}, Proposition A.0.2 (and the fact that $\Lambda_s$ is equipped with a topology). We omit the details. 
\end{remark}

Passing to the quasi-classical limit $\hbar\to 0$, we find that applying the operator $\exp(\hat{A})$ to $\D_\X$ corresponds to transforming the Lagrangian submanifold $\sL_\X$ by the (unquantized) operator $\exp(A)$. Hence we have the following
\begin{corollary}\label{qrrg0}
The Lagrangian submanifolds $\sL_s:=\sL_{(\bc,F)}$ and $\sL_\X$ are related by
\begin{equation*}
\begin{split}
\sL_s&=\exp \left(\sum_{k\geq 0} s_k\left(\sum_{m+h=k+1; m,h\geq 0}\frac{(A_m)_hz^{m-1}}{m!}+\frac{ch_k((q^*F)^{inv})}{2}\right)\right)\sL_\X\\ 
&=\exp \left(\sum_{m,h\geq 0} s_{m+h-1}\frac{(A_m)_hz^{m-1}}{m!}+\sum_{k\geq 0} s_k\frac{ch_k((q^*F)^{inv})}{2}\right)\sL_\X.
\end{split}
\end{equation*}
In particular, $\sL_s$ is the germ of a Lagrangian cone swept out by a finite dimensional family of subspaces (i.e. (\ref{overruled}) holds for $\sL_s$).
\end{corollary}

When $\X$ is a variety, the inertia stack $I\X$ is just $\X$ itself and Theorem \ref{qrr} reduces to \cite{CG}, Theorem 1 of Coates-Givental.
An interesting feature of Theorem \ref{qrr} is the presence of values of Bernoulli polynomials (see the definition of elements $A_m$) in place of Bernoulli numbers which appear in the quantum Riemann-Roch theorem for varieties (\cite{CG}, Theorem 1). It would be interesting to find a conceptual way to explain why this is the case.

\begin{remark}[Loop space interpretation]
There is a heuristic interpretation of the operator $\Delta=\exp(\sum_{k\geq 0} s_k (\sum_{m\geq 0} \frac{(A_m)_{k+1-m}z^{m-1}}{m!}+\frac{ch_k((q^* F)^{inv})}{2}))$ in terms of loop space $L\X$ (Section \ref{loop_interpret}), which we sketch below. On each component $\X_i\subset I\X\simeq L\X^{S^1}$, the $S^1$-action on $\X_i$ is trivial. This action is related to the $S^1$-action on the coarse space $X_i$ via the $r_i$-fold cover $S^1\to S^1$. We have $\X_i\times_{S^1}ES^1\simeq \X_i\times BS^1$. Let $pr_1, pr_2$ be the projections to factors and let $L^{1/r_i}$ denote the pullback by $pr_2$ of the universal line bundle over $BS^1$. Define $\F$ to be the $S^1$-equivariant vector bundle over $I\X\times BS^1$ whose restriction to $\X_i\times BS^1$ is $\bigoplus_{0\leq l\leq r_i-1} pr_1^*F_i^{(l)} \otimes (L^{1/r_i})^{\otimes l}$.

Consider the infinite product
$$\sqrt{\bc (F^{(0)})} \prod_{m=1}^{\infty} 
\bc \left(\F \otimes L^{-m}\right).$$ 
We interpret this as follows: Let $s(x):=\sum_{k\geq 0}s_k \frac{x^k}{k!}$. Note that if $x=c_1(\mathfrak{L})$ is the first Chern class of a line bundle $\mathfrak{L}$, then $s(x)=\sum_ks_k\ch_k(\mathfrak{L})=\log \bc(\mathfrak{L})$. We write 
$$\sum_{m>0} s\left(x+\frac{l}{r}z-mz\right) 
= \frac{e^{\frac{l}{r}z\frac{\partial}{\partial x}}\ z\frac{\partial}{\partial x}}{e^{z\frac{\partial}{\partial x}}-1} \left(z\frac{\partial}{\partial x}\right)^{-1}s(x) 
=\sum_{m\geq 0}\frac{B_m(\frac{l}{r})}{m!}\left(z\frac{\partial}{\partial x}\right)^{m-1}s(x).$$
Using this (and splitting principle) we expand 
$$\log \left(\prod_{m>0} \bc (F_i^{(l)}\otimes L^{\frac{l}{r_i}-m})\right) = 
\sum_{k\geq 0} s_k \sum_{m\geq 0} \frac{B_m(l/r_i)}{m!} \ch_{k+1-m}(F_i^{(l)})z^{m-1}.$$  
Thus the infinite product above gives rise the operator $\Delta$ after some simplification.
\end{remark}

\subsection{Relations to Hurwitz-Hodge integrals}
Let $G\subset SL_n(\com)$ be a finite subgroup. Consider a $G$-action on $\com^n$ without trivial factors so that $0\in \com^n$ is an isolated $G$-fixed point. Hurwitz-Hodge integrals (c.f. \cite{bgp}) arise in the study of orbifold Gromov-Witten theory of $[\com^n/G]$. More precisely, the components of $\Mbar_{g,n}([\com^n/G],0)$ parametrizing maps with images $[0/G]$ from orbicurves with stacky marked points may be identified with $\Mbar_{g,n}(BG)$, and the restriction of the virtual fundamental class is given by the Euler class $e(R^1f_*ev_{n+1}^*\mathcal{V})$, where $\mathcal{V}$ is the vector bundle over $BG$ defined by the $G$-representation $\com^n$, $f:\Mbar_{g,n+1}(BG)'\to \Mbar_{g,n}(BG)$ is the universal orbicurve, and $ev_{n+1}: \Mbar_{g,n+1}(BG)'\to BG$ is the universal orbifold stable map (see Section \ref{modulistablemap}). The integrals over $\Mbar_{g,n}(BG)$ of cohomology classes involving $e(R^1f_*ev_{n+1}^*\mathcal{V})$ are called Hurwitz-Hodge integrals. 

One may consider an equivariant version of this: Let $\com^*$ acts on $\com^n$ by scaling. This $\com^*$-action commutes with the $G$-action, hence descends to a $\com^*$-action on the stack $[\com^n/G]$. A $\com^*$-equivariant Hurwitz-Hodge integral $$\int_{[\Mbar_{g,n}(BG)]}(....)e_{\com^*}(R^1f_*ev_{n+1}^*\mathcal{V})$$ coincides with 
$$e_{\com^*}(R^0f_*ev_{n+1}^*\mathcal{V})\int_{[\Mbar_{g,n}(BG)]}(...)e_{\com^*}^{-1}(\mathcal{V}_{g,n,0}),$$ where $(...)$ denotes cohomology and/or descendant insertions. From this it is easy to conclude that Hurwitz-Hodge integrals can be determined by twisted orbifold Gromov-Witten invariants of $BG$. Theorem \ref{qrr} implies that descendant twisted orbifold Gromov-Witten invariants of $BG$ are determined by the usual descendant orbifold Gromov-Witten invariants of $BG$. In \cite{JK}, explicit formulas expressing descendant orbifold Gromov-Witten invariants of $BG$ in terms of descendant integrals on moduli stacks of  stable curves have been proven. It is interesting to combine these results to obtain formulas for Hurwitz-Hodge integrals. In genus zero, under additional assumptions, a procedure of explicitly computing $(e_{\com^*}^{-1},\mathcal{V})$-twisted orbifold Gromov-Witten invariants using information about usual orbifold Gromov-Witten invariants of $BG$ has been established, see \cite{ccit}. A method of computing Hurwitz-Hodge integrals directly using the Grothedieck-Riemann-Roch calculation in this paper has been developed by J. Zhou \cite{zhou}, and is used by him to prove the crepant resolution conjecture in higher genus for type $A$ surface singularities \cite{zhou2}. 

\section{Quantum Lefschetz}

\subsection{Twisting by Euler class}\label{qlh}
Consider the group $\com^*$ which acts trivially on $\X$ and on the vector bundle $F$ by scaling the fibers. Let $\lambda$ be the equivariant parameter and $e(\cdot)$ the $\com^*$-equivariant Euler class. In this section we consider the special case of twisting by $F$ and $\bc=e$. From $\lambda+x=\exp(\sum_{k\geq 0}s_k\frac{x^k}{k!})$, we find\footnote{Here we work over the ground ring $\Lambda_s$ with the values of $s_k$ specified by (\ref{s-values}).} 
\begin{equation}\label{s-values}
s_k =\left\{
 \begin{array}{rr}
 \ln\lambda , \quad &\mbox{$k=0$}\\
 \frac{(-1)^{k-1}(k-1)!}{\lambda^k},\quad &\mbox{$k>0$}.
 \end{array}\right.  
\end{equation}
Let $\rho_i^{l,j}$ be the Chern roots of $F_i^{(l)}$, $j=1,...,rank F_i^{(l)}$. The following is the case $\bc=e$ of Corollary \ref{qrrg0}.
\begin{corollary}\label{qrreu}
The Lagrangian cone $\sL_e:=\sL_{(e,F)}\subset\sH$ of the twisted theory is obtained from $\sL_\X$ by (ordinary) multiplication by the product over Chern roots $\rho_i^{l,j}$ of 
\[  \gamma_{\rho_i^{l,j}}(z) =\left\{
 \begin{array}{rr}
 \exp (\frac{(\rho_i^{l,j}+\lambda)\ln(\rho_i^{l,j}+\lambda)-(\rho_i^{l,j}+\lambda)}{z}+\ln\lambda(\frac{l}{r_i}-\frac{1}{2})+(\frac{l}{r_i}-\frac{1}{2})\ln(1+\frac{\rho_i^{l,j}}{\lambda})\\
+\sum_{m\geq 2}\frac{(-1)^mB_m(l/r_i)}{m(m-1)}(\frac{z}{\lambda+\rho_i^{l,j}})^{m-1})
 ,\\ \mbox{if $l\neq0$};\\
 
\exp(\frac{(\rho_i^{0,j}+\lambda)\ln(\rho_i^{0,j}+\lambda)-(\rho_i^{0,j}+\lambda)}{z}
+\sum_{m\geq 2}\frac{(-1)^mB_m}{m(m-1)}(\frac{z}{\lambda+\rho_i^{0,j}})^{m-1})
,\\ \mbox{if $l=0$}.\\
 \end{array}\right.  \]
\end{corollary}
\begin{proof}
We substitute the definition of $s_k$ into the statement of Corollary \ref{qrrg0} and express components $ch_h(F_i^{(l)})$ of the Chern characters using the Chern roots $\rho_i^{l,j}$. Then by using the identity $$\sum_{h\geq 0}s_{m+h-1}\frac{\rho^h}{h!}=\frac{d^{m-1}}{d\rho^{m-1}}\ln (\lambda+\rho)=\frac{(-1)^m(m-2)!}{(\lambda+\rho)^{m-1}}, \quad \text{for }m\geq 1$$ we check directly that the  $z^{m-1}$ terms, with $m\geq 1$, coincide with what are given in Corollary \ref{qrrg0}. For the $z^{-1}$ term, a direct calculation gives 
$$\frac{1}{z}\sum_{k\geq 0}s_kch_{k+1}(F_i^{(l)})=\frac{1}{z}\sum_{j}\left[((\rho_i^{l,j}+\lambda)\ln(\rho_i^{l,j}+\lambda)-(\rho_i^{l,j}+\lambda))-(\lambda\ln\lambda-\lambda)\right].$$
Since the operator $\frac{1}{z}$ preserves the cone $\sL_\X$, we may discard the term $\frac{\lambda\ln\lambda-\lambda}{z}$. The result follows.
\end{proof}

Our next goal is to extract from Corollary \ref{qrreu} more explicit information about genus zero invariants. For the rest of this section and Section \ref{ci}, we make the following assumption.

\begin{assump}\label{assumption_for_QLefschetz}
\hfill
\begin{enumerate}
\item The generic stabilizer of the stack $\X$ is trivial. 

\item The bundle $F$ is a direct sum $\oplus_j F_j$ of line bundles and each $F_j$ is a line bundle pulled back via the natural map $\pi: \X\to X$ to the coarse moduli space $X$.
\end{enumerate}
\end{assump}

\begin{remark}
\hfill
\begin{enumerate}
\renewcommand{\labelenumi}{(\roman{enumi})}
\item In the situation of Assumption \ref{assumption_for_QLefschetz}, the intersection index $\<c_1(F_j),\pi^*\beta\>:=c_1(F_j)\cdot\pi^*\beta$ is an integer for all effective curve classes $\beta$ of $X$. Let $L=\pi^*M$ be a line bundle on $\X$ that is pulled back from a line bundle $M$ on the coarse moduli space $X$. Then for any such $\beta$, we have $c_1(L)\cdot\pi^*\beta=c_1(M)\cdot\beta\in\mathbb{Z}$.
\item For each $i$ and $j$, the line bundle $q^*(F_j)|_{\X_i}$ has $\mu_{r_i}$-eigenvalue $1$. In other words, $q^*(F_j)|_{\X_i}=q^*(F_j)|_{\X_i}^{(0)}$.

\end{enumerate}
\end{remark}

We are interested in a more precise relationship between the $J$-function $J_\X$ and the twisted $J$-function $J_{(e,F)}$. We generalize the approach of \cite{CG}.
\begin{definition}\label{hypermod}
Put $\rho_{ji}:=c_1(q^*(F_j)|_{\X_i})\in H^2(\X_i)$ and $\rho_j:=c_1(F_j)\in H^2(\X)$. 
Define $I_F(t,z):=(I_F(t,z)_i)$ where
$$I_F(t,z)_i:=\sum_{d\in\text{Eff}(\X)} J_{i,d}(t,z)Q^d\prod_{j}\frac{\prod_{k=-\infty}^{\<\rho_j,d\>}(\lambda+\rho_{ji}+kz)}{\prod_{k=-\infty}^0(\lambda+\rho_{ji}+kz)}.$$
Following \cite{CG}, we call this the {\em hypergeometric modification} of $J_\X$.
\end{definition}

\begin{remark}
Assumption \ref{assumption_for_QLefschetz} is used to ensure that the intersection indices $\<c_1(F_j),d\>$ are integers, which is needed in order for the hypergeometric modification to be well-defined. (Note that for $d\in \text{Eff}(\X)$ there exists $\beta\in \text{Eff}(X)$ such that $d=\pi^*\beta$. If $\<c_1(F_j), d\>=\<c_1(F_j),\pi^*\beta\>$ are not integers, the product $\prod_{k=-\infty}^{\<\rho_j,d\>}(\lambda+\rho_{ji}+kz)/\prod_{k=-\infty}^0(\lambda+\rho_{ji}+kz)$ doesn't make sense.) 
\end{remark}

\begin{theorem}\label{iandj}
The family $$t\mapsto I_F(t,-z),\,\,\, t\in H^*(I\X)$$ of vectors in $(\sH_{(e,F)},\Omega_{(e,F)})$ lies on the Lagrangian submanifold $\tilde{\sL}_{(e,F)}$.
\end{theorem}
\begin{remark}
Theorem \ref{iandj} uses Assumption \ref{assumption_for_QLefschetz} in a essential way. A more general result of this kind is given in \cite{ccit}.
\end{remark}

\begin{proof}
This is a generalization of \cite{CG}, Theorem 2 (see also \cite{C}, Theorem 1.7.3). In view of Assumption \ref{assumption_for_QLefschetz} and Lemma \ref{ocup}, we may rewrite the operators $\gamma_{\rho_j^{l,i}}(z)$ in terms of the Chen-Ruan orbifold cup product (note that our assumption forces $l=0$). More precisely, 
$$\gamma_{\rho_j^{0,i}}(z)=\exp\left(\frac{(\rho_{j}+\lambda)\ln(\rho_{j}+\lambda)-(\rho_{j}+\lambda)}{z}
+\sum_{m\geq 2}\frac{(-1)^mB_m}{m(m-1)}\left(\frac{z}{\lambda+\rho_{j}}\right)^{m-1}\right)\cdot_{orb}|_{H^*(\X_i)}.$$
It is then straightforward to check that the argument of \cite{CG} and \cite{C} applies verbatim (of course with Corollary \ref{qrreu} replacing its manifold version). Details are left to the readers.

\end{proof}


\begin{corollary}\label{ij}
The tangent space $L_t$ to $\tilde{\sL}_{(e,F)}$ at the point $I_F(t,-z)$ is equal to the tangent space of $\sL_{(e,F)}$ at a unique point $J_{(e,F)}(\tau(t),-z)$, where $\tau(t)\in H^*(I\X, \com)\otimes\Lambda_s$.
\end{corollary}
\begin{proof}
Note that $I_F(t,z)\equiv J_\X(t,z)\,\,mod\, Q$. An easy calculation shows that the family $$t\mapsto I_F(t,-z),\,\,\, t\in H^*(I\X,\com)\otimes\Lambda_s $$ is transverse to $zL_t$ for every $t$. As pointed out in Corollary \ref{qrrg0}, (\ref{overruled}) holds for $\sL_{(e,F)}$. Thus the proof of \cite{ccit1}, Proposition 2.16 may be applied to show that $L_t$ is equal to the tangent space of $\sL_{(e,F)}$ at a unique point $J_{(e,F)}(\tau(t),-z)$.
\end{proof}
\begin{remark}
\hfill
\begin{enumerate}
\renewcommand{\labelenumi}{(\roman{enumi})}
\item Intuitively this Corollary may be interpreted as saying that the intersection of $zL_t$ with $\{-z+z\sH_-\}\cap\tilde{\sL}_{(e,F)}$ is equal to
$$J_{(e,F)}(\tau(t),-z)\in -z+\tau(t)+\sH_-,$$ 
where $\tau(t)\in H^*(I\X,\com)\otimes\Lambda_s$ is defined by this intersection.

\item This Corollary should be viewed as a procedure of computing $J_{(e,F)}$ from $I_F$. This procedure is related to Birkhoff factorization in the theory of loop groups. More precisely, this procedure applied to the first derivatives of $I_F$ is indeed an example of Birkhoff factorization.
\item The map $t\mapsto \tau=\tau(t)$ may be viewed as the ``mirror map''. This Corollary gives a geometric description of this map.
\end{enumerate}
\end{remark}

\subsection{Complete intersections}\label{ci}
In this Section we apply Corollary \ref{ij} to vector bundles with some positivity property to deduce relationships between orbifold Gromov-Witten invariants of a complete intersection orbifold and orbifold Gromov-Witten invariants of its ambient orbifold.

\begin{definition}
A line bundle $F$ over $\X$ is called {\em convex} if $H^1(\C,f^*F)=0$ for all $1$-pointed genus-$0$ orbifold stable maps $f:(\C, \Sigma)\to \X$.
\end{definition}
\begin{example}
Let $L:=\pi^*M$ be a line bundle on $\X$ that is the pullback of a line bundle $M$ on the coarse moduli space $X$. For an orbifold stable map $f:\C\to \X$ with induced map $\fbar: C\to X$ between coarse moduli spaces, we have
$$H^1(\C,f^*L)=H^1(\C,f^*\pi^*M)=H^1(\C,\pibar^*\fbar^*M)=H^1(C,\fbar^*M).$$
Here $\pibar:\C\to C$ is the map to the coarse curve. From this we see that the bundle $L$ is convex if $M$ is convex in the usual sense.
\end{example}

The following Proposition follows from \cite{KiKrPa}.
\begin{proposition}\label{virclass}
Let $F=\oplus_j F_j$ be a direct sum of convex line bundles. Let $\Y$ be the zero locus of a regular section of $F$, and $j_{0,n,d}: \Mbar_{0,n}(\Y,d)\to \Mbar_{0,n}(\X,d)$ the induced map. Then $j_{0,n,d*}[\Mbar_{0,n}(\Y,d)]^w=\mathfrak{e}(F_{0,n,d})\cap [\Mbar_{0,n}(\X,d)]^w$,  where $\mathfrak{e}(\cdot)$ denotes the non-equivariant Euler class. 
\end{proposition}
In the situation of Proposition \ref{virclass} let $j:\Y\to\X$ be the inclusion. Let $I_{\X,\Y}(t,z)$ and $J_{\X,\Y}(\tau,z)$ be the nonequivariant limits $\lambda\to 0$ of $I_F(t,z)$ and $J_{(e,F)}(\tau,z)$ respectively. Let $F_{0,n+1,d}'$ be the kernel of the evaluation map $F_{0,n+1,d}\to ev_{n+1}^*q^*F$ at the $(n+1)$-st marked point. Note that the image of the evaluation map is contained in $ev_{n+1}^*((q^*F)^{inv})$. The non-equivariant limit $J_{\X,\Y}$ can be written as 
$$J_{\X,\Y}(t,z)=z+t+\sum_{n\geq 0,d\in\text{Eff}(\X)}\frac{Q^d}{n!}ev_{n+1*}(ev_1^*t\cup...\cup ev_n^*t\cup\frac{\mathfrak{e}(F_{0,n+1,d}')}{z-\bpsi_{n+1}}).$$
Together with Proposition \ref{virclass} this implies that
\begin{equation}\label{compare}
\mathfrak{e}((q^*F)^{inv})J_{\X,\Y}(u,z)=j_*J_\Y(j^*u,z)
\end{equation} where on the right-hand side the Novikov rings should be changed according to $\text{Eff}(\Y)\to \text{Eff}(\X)$.

By taking the nonequivariant limit, we obtain
\begin{corollary}\label{lef}
Let $\X, \Y$ and $F$ be as in Proposition \ref{virclass}. Then $I_{\X,\Y}(t,-z)$ and $J_{\X,\Y}(\tau,-z)$ determine the same cone. Moreover, $J_{\X,\Y}(\tau,-z)$ is determined from $I_{\X,\Y}(t,-z)$ by the procedure described in Corollary \ref{ij}, followed by the mirror map $t\mapsto \tau$.
\end{corollary}
This is a generalization of ``Quantum Lefschetz Hyperplane Principle'' (see \cite{Gi2}, \cite{Ki}, \cite{B}, \cite{L}, \cite{Ga}, \cite{CG}) to Deligne-Mumford stacks.

We now restrict to the small parameter space $H^{\leq 2}(\X)$. We continue to assume that $F=\oplus_j F_j$ is a direct sum of convex line bundles.
\begin{proposition}
Let $\{\gamma_k\}$ be a basis for $H^{\leq 2}(\X)$. If $c_1(F)\leq c_1(T_\X)$, then for $t\in H^{\leq 2}(\X)$ we have an expansion $$I_{\X,\Y}(t,z)=zF(t)+\sum_k G^k(t)\gamma_k+O(z^{-1}),$$ where $F(t)$ and $G^k(t)$ are certain scalar-valued functions with $F(t)$ invertible. 
\end{proposition}
\begin{proof}
We have 
$$I_F(t,z)_i=z+t+\sum_{d>0} J_{i,d}(t,z)Q^d\prod_{j}\prod_{k=1}^{\<\rho_j,d\>}(\lambda+\rho_{ji}+kz) +O(z^{-1}).$$ 
Recall that
$$J_{i,d}(t,z)=\sum_{n\geq 0}\frac{Q^d}{n!}\sum_{k\geq 0,\alpha}\<t,...,t,\phi_\alpha\bpsi^k\>_{0,n+1,d}\frac{\phi^\alpha}{z^{k+1}},$$ where $\{\phi^\alpha\}$ is an additive basis of $H^*(\X_i)$.
We need to identify the highest power of $z$ in $J_{i,d}(t,z)$. For this one should take $t\in H^2(\X)$ and $orbdeg(\phi_\alpha)$ to be as large as possible. In view of Lemma \ref{agerel}, the largest possible orbifold degree is $2dim_\com\X$. Therefore, by (\ref{vir_dimension}), the largest power of $z$ in $J_{i,d}(t,z)$ is $1-\<c_1(T\X),\pi^*d\>$. 

The highest power of $z$ occurring in $$J_{i,d}(t,z)Q^d\prod_{j}\prod_{k=1}^{\<\rho_j,d\>}(\lambda+\rho_{ji}+kz)$$
is equal to $$1+\<c_1(F),\pi^*d\>-\<c_1(T_\X),\pi^*d\>.$$ By our assumption, this is at most $1$. If this is equal to $1$, then the class $\phi^\alpha$ has orbifold degree $0$. In order to have $z^0$ term, we must have $orbdeg(\phi_\alpha)\geq 2dim_\com\X-2$, which implies that $orbdeg(\phi^\alpha)\leq 2$. Also, we see that $F(t)\equiv 1 \,\,(mod \, Q)$. The Proposition follows.
\end{proof}

Since $J_{\X,\Y}(\tau,z)$ is characterized by the asymptotic $J_{\X,\Y}(\tau,z)=z+\tau+O(z^{-1}),$
by comparing the asymptotics of $I_{\X,\Y}$ and $J_{\X,\Y}$, we obtain
\begin{corollary}\label{comp}
If $c_1(F)\leq c_1(T_\X)$, then the restriction of $J_{\X,\Y}(\tau,z)$ to small parameter space $H^{\leq 2}(\X)$ is given by
$$J_{\X,\Y}(\tau,z)=\frac{I_{\X,\Y}(t,z)}{F(t)},\quad \text{where} \quad \tau=\sum_k \frac{G^k(t)}{F(t)}\gamma_k.$$
\end{corollary}
This may be regarded as a {\em mirror formula} for complete intersection orbifolds. Once the $J$-function of $\X$ is known, part of the $J$-function of $\Y$ that involves classes pulled back from $\X$ can be computed by Corollary \ref{comp} and (\ref{compare}).

\section{Quantum Serre duality}\label{qse}
The so-called ``Quantum Serre duality'' (\cite{Gi2-5}, \cite{CG}) is formulated as a relation between $(\bc,F)$-twisted invariants and invariants twisted by the ``dual data'' $(\bc^\vee, F^\vee)$ defined below. In this Section we prove such a relation for Deligne-Mumford stacks.

\subsection{General case}
We again consider the general case of twisting by a vector bundle $F$ over $\X$ and multiplicative invertible characteristic class $\bc(\cdot)=\exp{(\sum_k s_kch_k(\cdot))}$. Here we do not require Assumption \ref{assumption_for_QLefschetz}. Consider the dual case of twisting by the dual bundle $F^\vee$ and the class $$\bc^\vee (\cdot):=\exp \left(\sum_{k\geq 0} (-1)^{k+1} s_k ch_k(\cdot)\right).$$
Note that $\bc^\vee(F^\vee)=\frac{1}{\bc(F)}$. An application of Theorem \ref{qrr} yields the following relation between the potentials $\D_{(\bc, F)}$ and $\D_{(\bc^\vee,F^\vee)}$.

\begin{theorem}[Quantum Serre duality for orbifolds]\label{qserre}
Let $\bt^\vee(z)=\bc((q^*F)^{inv})\bt(z)+(1-\bc((q^*F)^{inv}))z$. Then we have 
$$\D_{(\bc^\vee,F^\vee)}(\bt^\vee)=\exp\left(-s_0\text{rank}\,F\<\bpsi\>_{1,1,0}\right)\D_{(\bc, F)}(\bt).$$
\end{theorem}

\begin{proof}
One may prove this result by comparing the formulas for $\D_{(\bc, F)}$ and $\D_{(\bc^\vee,F^\vee)}$ given by Theorem \ref{qrr}. We proceed differently by comparing the differential equation (\ref{qrrd}) for $(\bc, F)$ and $(\bc^\vee, F^\vee)$. The equation satisfied by $\D_{(\bc, F)}$ is 
\begin{equation}\label{usual_eq}
\frac{\partial\D_{(\bc,F)}}{\partial s_k}=\left[\left(\sum_{m+h=k+1; m,h\geq 0}\frac{(A_m)_hz^{m-1}}{m!}+\frac{ch_k(q^*F^{inv})}{2}\right)^\wedge+C_k\right]\D_{(\bc,F)}.
\end{equation}
We write the equation satisfied by $\D_{(\bc^\vee, F^\vee)}$ as
\begin{equation}\label{veeeq1}
(-1)^{k+1}\frac{\partial\D_{(\bc^\vee, F^\vee)}}{\partial s_k}=\left[\left(\sum_{m+h=k+1; m,h\geq 0}\frac{(A_m^\vee)_hz^{m-1}}{m!}+\frac{ch_k(q^*F^{\vee inv})}{2}\right)^\wedge+C_k^\vee\right]\D_{(\bc^\vee, F^\vee)}.
\end{equation}
For a fixed $i\in \sI$, the first term on the right-hand side of (\ref{veeeq1}) is
$$
\sum_{\overset{m+h=k+1}{m,h\geq 0}}\frac{1}{m!}\sum_{0\leq l\leq r_i-1}ch_h(F_i^{\vee (l)})B_m(l/r_i)z^{m-1}.
$$
We now analyze this for each fixed $m,h$. Using $F_i^{\vee (l)}=F_i^{(r_i-l)\vee}$ for $0<l<r_i$ and $F_i^{\vee (0)}=F_i^{(0)\vee}$, we can write this as
\begin{equation*}
\begin{split}
&\frac{1}{m!}\sum_{1\leq l\leq r_i-1}ch_h(F_i^{(r_i-l)\vee})B_m(l/r_i)z^{m-1}+\frac{1}{m!}ch_h(F_i^{(0)\vee})B_mz^{m-1}\\
=&(-1)^{m+h}\frac{1}{m!}\sum_{1\leq l\leq r_i-1}ch_h(F_i^{(r_i-l)})B_m(\frac{r_i-l}{r_i})z^{m-1}+(-1)^h\frac{1}{m!}ch_h(F_i^{(0)})B_mz^{m-1}\\
=&(-1)^{m+h}\frac{1}{m!}\sum_{1\leq l\leq r_i-1}ch_h(F_i^{(l)})B_m(\frac{l}{r_i})z^{m-1}+(-1)^h\frac{1}{m!}ch_h(F_i^{(0)})B_mz^{m-1}.
\end{split}
\end{equation*}
For $m\neq 1$, since $(-1)^mB_m=B_m$, this sum is 
$$(-1)^{k+1}\frac{1}{m!}\sum_{0\leq l\leq r_i-1}ch_h(F_i^{(l)})B_m(\frac{l}{r_i})z^{m-1},$$ where we use $m+h=k+1$. For $m=1$ and $h=k$, we have
$$(-1)^kch_k(F_i^{(0)})B_1=\frac{-1}{2}(-1)^kch_k(F_i^{(0)}),$$ which cancels with the term $ch_k(F_i^{\vee (0)})/2$.

Therefore we conclude that (\ref{veeeq1}) is 
\begin{equation}\label{veeeq}
\frac{\partial\D_{(\bc^\vee, F^\vee)}}{\partial s_k}=\left[\left(\sum_{m+h=k+1; m,h\geq 0}\frac{(A_m)_hz^{m-1}}{m!}+\frac{ch_k(q^*F^{inv})}{2}\right)^\wedge+C_k^\vee\right]\D_{(\bc^\vee, F^\vee)}.
\end{equation}
By Lemma \ref{number_C_k}, $C_k^\vee=0$ for $k\geq 1$, and $$C_0^\vee=\frac{1}{2}\text{rank}\,F^\vee\<\bpsi\>_{1,1,0}-\<c_1(F^\vee)\>_{1,1,0}=\frac{1}{2}\text{rank}\,F\<\bpsi\>_{1,1,0}+\<c_1(F)\>_{1,1,0}.$$
The result follows by comparing (\ref{veeeq}) with (\ref{usual_eq}),
\end{proof}

\subsection{Euler class}\label{qserre_for_euler}
We consider the case of twisting by a $\com^*$-equivariant Euler class $e(\cdot)$, where $\com^*$ acts on $F$ by scaling the fibers. Let the dual bundle $F^\vee$ be equipped with the dual $\com^*$-action and let $e^{-1}(\cdot)$ be the inverse $\com^*$-equivariant Euler class. If $\rho_j$ are the Chern roots of $F$, then $e^{-1}(F^\vee)=\prod_j (-\lambda-\rho_j)^{-1}$. The main result of this section, Proposition \ref{eulerserre}, is a relation between $(e,F)$-twisted invariants and $(e^{-1},F)$-twisted invariants. Note that this is not a special case of Theorem \ref{qserre} since $e^{-1}\neq e^\vee$.

Let $F$ be a vector bundle on $\X$. For a component $\X_i$ of $I\X$ we define  the age of the bundle $F$ on $\X_i$ to be
$$age(F_i):=\sum_{1\leq l\leq r_i-1} \frac{l}{r_i}rank F_i^{(l)}.$$ The bundle $F_i^{mov}$ is defined to be $\oplus_{1\leq l\leq r_i-1} F_i^{(l)}$. Let $M: H^*(I\X)\to H^*(I\X)$ be defined as follows: on $H^*(\X_i)$, $M$ is the multiplication by the number $(-1)^{\frac{1}{2}rank F_i^{mov}-age(F_i)}$. Put $$\bt^*(z)=z+(-1)^{\frac{1}{2}rank (q^*F)^{inv}}M e((q^*F)^{inv})(\bt(z)-1z),$$ and define a change $\diamond: Q^d\mapsto Q^d(-1)^{\<ch_1(F),d\>}$ in the Novikov ring. 

\begin{proposition}\label{eulerserre}
We have 
\begin{equation*}
\begin{split}
&\quad \exp\left(\frac{\pi\sqrt{-1}}{2}\text{rank}F\<\bpsi\>_{1,1,0}+\pi\sqrt{-1}\<c_1(F)\>_{1,1,0}\right)\D_{(e^{-1},F^\vee)}(\bt^*,Q)\\
&=\exp\left(-\ln \lambda \text{rank}F\<\bpsi\>_{1,1,0}\right)\D_{(e,F)}(\bt,\diamond Q).
\end{split}
\end{equation*}
\end{proposition}

\begin{proof}
If we write $e^{-1}(\cdot)=\exp{(\sum_{k\geq 0}s_k^*ch_k(\cdot))}$ and $e(\cdot)=\exp{(\sum_{k\geq 0} s_kch_k(\cdot))}$, then we find that $s_k^*=(-1)^{k+1}s_k$ for $k>0$ and $s_0^*=-s_0-\pi\sqrt{-1}$. The proof of Theorem \ref{qserre} shows that $\D_{s^*}$ satisfies the differential equation
\begin{equation}\label{sstar}
\frac{\partial\D_{s^*}}{\partial s_k}=\left[\left(\sum_{m+h=k+1; m,h\geq 0}\frac{(A_m)_hz^{m-1}}{m!}+\frac{ch_k(q^*F^{inv})}{2}\right)^\wedge+C_k^\vee\right]\D_{s^*}.
\end{equation}
Also, $$\D_{s^*}|_{s_k=0}=\D_{s^*}|_{s_0^*=-\pi\sqrt{-1}, s_k^*=0 \text{ for }k>0}.$$
By Theorem \ref{qrr}, we have
\begin{equation*}
\begin{split}
&\exp\left(\frac{\pi\sqrt{-1}}{2}\text{rank}F^\vee\<\bpsi\>_{1,1,0}-\pi\sqrt{-1}\<c_1(F^\vee)\>_{1,1,0}\right)\D_{s^*}\vert_{\overset{s_0^*=-\pi\sqrt{-1}}{s_k^*=0 \text{ for }k>0}}\\
=&\exp(-\pi\sqrt{-1}[(A_1)_0+\frac{1}{2}ch_0(q^*F^{inv})]^\wedge)\exp(-\pi\sqrt{-1}((A_0)_1/z)^\wedge)\D_0.
\end{split}
\end{equation*}
For a fixed $i\in \sI$, we have
\begin{equation*}
\begin{split}
&((A_1)_0+\frac{1}{2}ch_0(q^*F^{inv}))|_{\X_i}=\sum_{0<l<r_i}ch_0(F_i^{(l)})B_1(l/r_i),\\ 
&(A_0)_1/z|_{\X_i}=ch_1(q^*F|_{\X_i})/z.
\end{split}
\end{equation*}
The operator $$\exp(-\pi\sqrt{-1}((A_0)_1/z)^\wedge)$$ can be computed directly using Appendix \ref{quantized_operator} and is seen to yield the change $\diamond$ via the divisor flow. The operator $\exp(-\pi\sqrt{-1}[(A_1)_0+\frac{1}{2}ch_0(q^*F^{inv})]^\wedge)$ may be computed using Appendix \ref{quantized_operator}, one sees that it yields the operator $M$.

Solving the equation (\ref{sstar}) and using the expression of $\D_{s^*}|_{s_k=0}$ yields the desired formula. 
\end{proof}

\section{Proof of Theorem \ref{qrr}}\label{pfqrr}
In this section we prove Theorem \ref{qrr}. The proof is rather lengthy and somewhat unpleasant, however the idea (which we borrowed from \cite{CG}) of the proof is quite simple.

\subsection{Overview}\label{proof_overview}
For the convenience of what follows, we introduce a new notation.
\begin{definition}
Let $a_j\in H^{p_j}(\X_{i_j},\com), j=1,...,n$ be cohomology classes, $A\in H^*(\Mbar_{g,n}(\X,d; i_1,...,i_n))$, and $k_1,...,k_n$ nonnegative integers. Define $$\<a_1\bpsi^{k_1},...,a_n\bpsi^{k_n};A\>_{g,n,d}:=\int_{[\Mbar_{g,n}(\X,d; i_1,...,i_n)]^w}(ev_1^*a_1)\bpsi_1^{k_1}\cup ...\cup(ev_n^*a_n)\bpsi_n^{k_n}\cup A.$$
\end{definition}
Let us explain the structure of the proof. As explained in Section \ref{oqrr}, the twisted descendant potentials $\D_s$ are viewed as a family of asymptotic elements depending on variables $s=(s_0, s_1,...)$. We know that
$$\D_s|_{s_0=s_1=...=0}=\D_\X.$$
To prove Theorem \ref{qrr}, we find a system of differential equations in $s_k$ satisfied by $\D_s$, and solve the initial value problem with the initial condition given by above. Such a system of differential equations is found by doing the naive thing: compute $\partial \D_s/\partial s_k$. A direct computation yields  
$$\D_s^{-1}\frac{\partial \D_s}{\partial s_k}$$
\begin{equation}\label{dds_k}
=\sum_{g,n,d}\frac{Q^d\hbar^{g-1}}{n!}\<\bt(z),...,\bt(z);\frac{\partial}{\partial s_k}\bc(F_{g,n,d})\>_{g,n,d}
+\sum_{g,n,d}\frac{Q^d\hbar^{g-1}}{(n-1)!}\<\frac{\partial}{\partial s_k}\bt(z),...,\bt(z);\bc(F_{g,n,d})\>_{g,n,d}.
\end{equation}

The second term in (\ref{dds_k}), called the derivative term, is equal to 
\begin{equation}\label{der}
-\frac{1}{2}\sum_{g,n,d}\frac{Q^d\hbar^{g-1}}{(n-1)!}\<ch_k((q^*F)^{inv})(\bt(z)-1z),...,\bt(z);\bc(F_{g,n,d})\>_{g,n,d}.
\end{equation}
This can be seen from 
\begin{equation*}
\begin{split}
&\frac{\partial}{\partial s_k}\bt(z)=\frac{\partial}{\partial s_k}(\bc((q^*F)^{inv})^{-1/2}\bq(z)+1z)\\
&=-\frac{1}{2}\bc((q^*F)^{inv})^{-1/2}ch_k((q^*F)^{inv})\bq(z)=-\frac{1}{2}ch_k((q^*F)^{inv})(\bt(z)-1z).
\end{split}
\end{equation*}

Since $$\frac{\partial}{\partial s_k}\bc(F_{g,n,d})=\bc(F_{g,n,d})ch_k(F_{g,n,d}),$$ the first term in (\ref{dds_k}) is equal to 
\begin{equation}\label{chern_term}
\sum_{g,n,d}\frac{Q^d\hbar^{g-1}}{n!}\<\bt(z),...,\bt(z);\bc(F_{g,n,d})ch_k(F_{g,n,d})\>_{g,n,d}.
\end{equation}
The Chern character $ch_k(F_{g,n,d})$ appearing in (\ref{chern_term}) will be computed by applying Grothendieck-Riemann-Roch formula. The result is then combined with (\ref{der}) to obtained the following differential equation, written using Givental's formalism: 
\begin{equation}\label{qrrd}
\frac{\partial\D_s}{\partial s_k}=\left[\left(\sum_{m+h=k+1; m,h\geq 0}\frac{(A_m)_hz^{m-1}}{m!}+\frac{ch_k((q^*F)^{inv})}{2}\right)^\wedge+C_k\right]\D_s.
\end{equation}
Here we define $$C_k:=-\int_{[\Mbar_{1,1}(\X,0)']^w}\sum_{a+b=k+1; a, b\geq 0} ev^*ch_a(F)(Td^\vee(L_1))_b\bc(F_{1,1,0}).$$ The term $(-)_b$ means the degree $2b$ part of a cohomology class, and $Td^\vee$ is the dual Todd class defined by the property that $Td^\vee(L^\vee)=Td (L)$ for any line bundle $L$. Recall that the superscript $^\wedge$ indicates the quantization, discussed in Section \ref{quantization}.

The proof of Theorem \ref{qrr} will be completed in the next few sections. In the next Section we derive Theorem \ref{qrr} from (\ref{qrrd}). The computation of $ch_k(F_{g,n,d})$ by applying Grothendieck-Riemann-Roch formula will be presented in Section \ref{grrcal}. In Section \ref{find_de} we derive the equation (\ref{qrrd}) from these computations.

\subsection{From (\ref{qrrd}) to (\ref{qrr})}
We first derive Theorem \ref{qrr} from (\ref{qrrd}). 
\begin{lemma}\label{number_C_k}
$C_k=0$ for $k\geq 1$. $C_0=\frac{1}{2}\text{rank}\,F\<\bpsi\>_{1,1,0}-\<c_1(F)\>_{1,1,0}$.
\end{lemma}
\begin{proof}
The virtual complex dimension of $\Mbar_{1,1}(\X,0)'$ is $1$ (note that the marked point is non-stacky). The integrand involved in $C_k$ is of degree at least $2(k+1)$. So $C_k=0, k\geq 1$ for dimension reason. The degree $2$ part of the integrand of $C_0$ is $(ev^*ch_0(F)(Td^\vee(L_1))_1+ev^*ch_1(F))\bc(F_{1,1,0})_0$, where $\bc(F_{1,1,0})_0=\exp(s_0ch_0(F_{1,1,0}))$ denotes the degree-$0$ part of $\bc(F_{1,1,0})$. By Riemann-Roch, we find that the virtual rank of $F_{1,1,0}$ is $0$, thus $ch_0(F_{1,1,0})=0$ and $\bc(F_{1,1,0})_0=1$. We conclude by observing that $(Td^\vee(L_1))_1=-\frac{1}{2}\bpsi$.
\end{proof}
\begin{remark}
Our proof of Lemma \ref{number_C_k} uses a dimension argument and is valid in {\em non-equivariant} Gromov-Witten theory. The exact evaluation of $C_k$ in equivariant Gromov-Witten theory requires an explicit description of the moduli stack $\Mbar_{1,1}(\X,0)'$ and its virtual class. Such a description is not known for Deligne-Mumford stacks $\X$, thus an exact evaluation of $C_k$ in equivariant Gromov-Witten theory remains unknown. If the torus acts with isolated fixed points, virtual localization formula yields a calculation of $C_k$. We will not discuss it here.
\end{remark}
For simplicity, write 
$$\alpha_k:=\left(\sum_{m>0} \frac{(A_m)_{k+1-m}z^{m-1}}{m!}+\frac{ch_k(q^* F^{inv})}{2}\right)^\wedge,\quad
\beta_k:=\left(\frac{(A_0)_{k+1}}{z}\right)^\wedge.$$
As explained in \cite{C}, Example 1.3.4.1, the cocycle $\C(\sum_j s_j\alpha_j,\beta_k)$ is equal to 
\begin{equation}\label{cocycle_GRR}
\begin{split}
&\C\left(\sum_{j\geq 0}\frac{s_j (A_2)_{j-1}z}{2}, \frac{(A_0)_{k+1}}{z}\right)\\
=&-\frac{1}{2}\text{str}\left((A_0)_{k+1}\sum_{j\geq 0}\frac{s_j (A_2)_{j-1}}{2}\right)\\
=&-\frac{1}{2}\int_{II\X}e(T_{II\X})\wedge (Iq)^*(A_0)_{k+1}\wedge\sum_{j>0}s_j\frac{(Iq)^*(A_2)_{j-1}}{2}\quad \text{by Appendix \ref{cocycle_calculation}}\\
=&0 \quad \text{since the degrees of integrands exceed the dimension of }II\X.
\end{split}
\end{equation}

Solving (\ref{qrrd}), we find
$$\D_s=\exp(\sum_k s_k\alpha_k)\exp(\sum_k s_k\beta_k)\exp(s_0C_0)\D_0,$$
which gives Theorem \ref{qrr}.

\begin{remark}
Our derivation of Theorem \ref{qrr} from (\ref{qrrd}) uses the exact values of $C_k$ and the cocycles, and is valid for {\em non-equivariant} Gromov-Witten theory. In this paper we only consider non-equivariant Gromov-Witten theory. Note however that (\ref{qrrd}) is valid in full generality. 
\end{remark}

\subsection{GRR Calculation}\label{grrcal}
In this Section we compute $ch_k(F_{g,n,d})\cap [\Mbar_{g,n}(\X,d)]^{vir}$ by applying Grothendieck-Riemann-Roch formula. For technical reasons we proceed as follows. The construction in \cite{AGOT} using Hilbert functors for Deligne-Mumford stacks provides a family of orbicurves $$\U\to \M$$  over a smooth base stack $\M$ and an embedding $$\Mbar_{g,n}(\X,d)\to \M$$ satisfying the following 
\begin{property}\label{property_of_nice_family_of_curves}
\hfill
\begin{enumerate}
\item
the family $\U\to\M$ pulls back to the universal family over $\Mbar_{g,n}(\X,d)$,
\item
the vector bundle $E=ev_{n+1}^*F$ extends to a vector bundle over $\U$,
\item
the Kodaira-Spencer map $T_m\M\to Ext^1(\sO_{\U_m}, \sO_{\U_m})$ is surjective for all $m\in \M$.
\end{enumerate}
\end{property}

Details can be found in \cite{AGOT}, Proposition 3.1.1. We check that Grothendieck-Riemann-Roch formula (Corollary \ref{rr}) can be applied to $\U\to\M$. First note that Property \ref{property_of_nice_family_of_curves} and the smoothness of $\M$ imply that $\U$ is a smooth Deligne-Mumford stack. By the construction in \cite{AGOT}, $\U\to \M$ factors as $$\U\to \bar{\mathbb{A}}\times \M\to \M,$$ where $\bar{\mathbb{A}}$ is smooth, $\U\to \bar{\mathbb{A}}\times \M$ is a regular embedding, and $\bar{\mathbb{A}}\times \M\to \M$ is the projection. Therefore $\U\to \M$ is a local complete intersection (lci) morphism\footnote{The notion of a lci morphism for Deligne-Mumford stacks is the same as that for schemes (\cite{Fu}, Appendix B.7.6). }. Moreover, since the relative tangent bundle of $\bar{\mathbb{A}}\times \M\to \M$ is just the tangent bundle of $\bar{\mathbb{A}}$ pulled back to $\bar{\mathbb{A}}\times \M$, it follows that the lci virtual tangent bundle of $\U\to \M$ coincides with its relative tangent bundle. 

We can compute $ch(f_*ev_{n+1}^*F)\cap[\Mbar_{g,n}(\X,d)]^{vir}$ by applying Corollary \ref{rr} to the morphism $\U\to \M$ then capping with $[\Mbar_{g,n}(\X,d)]^{vir}$. Therefore for the rest of this section, we assume Property \ref{property_of_nice_family_of_curves}. To avoid introducing cumbersome new notations, we will express our computations {\em as if} they were done for the morphism $\Mbar_{g,n+1}(\X,d)'\to \Mbar_{g,n}(\X,d)$. Namely we assume that the moduli stack $\Mbar_{g,n}(\X, d)$ is smooth and its universal family has everywhere-surjective Kodaira-Spencer map.

The Grothendieck-Riemann-Roch calculation we need is done individually for each component $\Mbar_{g,n+1}(\X,d;i_1,...,i_n,0)$. We begin with an analysis of the components of the inertia stacks $I\Mbar_{g,n+1}(\X,d;i_1,...,i_n,0)$ required for this calculation. There are three types of components of the inertia stack $I\Mbar_{g,n+1}(\X,d;i_1,...,i_n,0)$ that are mapped to $\Mbar_{g,n}(\X,d;i_1,...,i_n)$:
\begin{enumerate}
\item the main stratum $\Mbar_{g,n+1}(\X,d;i_1,...,i_n,0)$,
\item the divisors of marked points $\D_{j,(i_1,...,i_n)}$, and
\item the locus of nodes $\Z_{r,(i_1,...,i_n)}$.
\end{enumerate}
In the rest of this Section we work out contributions from each of them to GRR formula of $ch(f_*ev_{n+1}^*F)\cap[\Mbar_{g,n}(\X,d)]^{vir}$. 


\subsubsection{Main stratum}
The computation on the main stratum does not depend on $(i_1,...,i_n)$. To simplify notation, we describe it for $f:\Mbar_{g,n+1}(\X,d)'\to\Mbar_{g,n}(\X,d)$.

The restrictions of $\tch(E)$ and $\ttd(T_f)$ to $\Mbar_{g,n+1}(\X,d)'$ are $ch(E)$ and $Td(T_f)$ respectively. To compute $Td(T_f)=Td^\vee(\Omega_f)$, we use the following Lemma.
\begin{lemma}\label{cot}
There are exact sequences of sheaves
$$0\to \Omega_f\to \omega_f\to i_*\sO_\Z\to 0,$$
$$0\to \omega_f\to L_{n+1}\to \oplus_j s_{j*}\sO_{\D_j}\to 0,$$
where 
\begin{itemize}
\item $L_{n+1}$ is the tautological line bundle on $\Mbar_{g,n+1}(\X,d)'$ corresponding to the $(n+1)$-st marked point.
\item $s_j:\D_j\to \Mbar_{g,n+1}(\X,d)'$ are inclusions of the divisors of marked points.
\item $i: \Z\to \Mbar_{g,n+1}(\X,d)'$ is the inclusion of the locus of the nodes.
\end{itemize}
\end{lemma}
\begin{proof}
We prove the first exact sequence. The second sequence can be proved by a similar argument. Away from $\Z$, two sheaves $\Omega_f$ and $\omega_f$ are the same. Consider a family $S\leftarrow \C\overset{f}{\to}\X$ of orbifold stable maps with $S={\rm Spec}R$ such that the fiber of $\C/S$ over a point of $S$ is a nodal orbicurve. \'Etale-locally near a node\footnote{We use the condition on Kodaira-Spencer map to give this description.}, we may write $\C$ as the quotient $[U/\mu_r]$ where $U$ is the nodal curve ${\rm Spec}(R[z,w]/(zw-t))$ and $\mu_r$ acts on $U$ via $$(z,w)\mapsto (\zeta_rz,\zeta_r^{-1}w).$$

On this neighborhood, the dualizing sheaf $\omega_f$ corresponds to the $\mu_r$-equivariant sheaf $\omega_U$ with invariant generator $\frac{dz\wedge dw}{d(zw)}$. The sheaf $\Omega_f$ of K\"ahler differentials corresponds to the $\mu_r$-equivariant sheaf $\Omega_U$ with generators $dz, dw$ and a relation $wdz+zdw=0$. There is an equivariant inclusion $\Omega_f\hookrightarrow \omega_f$ defined by $$dz\mapsto z\frac{dz\wedge dw}{d(zw)},\quad dw\mapsto -w\frac{dz\wedge dw}{d(zw)}.$$ The cokernel corresponds to the $\mu_r$-equivariant sheaf generated by $\frac{dz\wedge dw}{d(zw)}$ with coefficients in $\sO_S$. This sheaf is identified with $i_*\sO_\Z$, proving the first exact sequence.
\end{proof}

Therefore we have
$$Td^\vee(\Omega_f)=Td^\vee(L_{n+1})Td^\vee(-i_*\sO_\Z))\prod_j Td^\vee(-s_{j*}\sO_{\D_j}).$$
Note that the $\D_j$'s and $\Z$ are disjoint, and the restrictions of $L_{n+1}$ to them are trivial. So we have 
\begin{equation*}
\begin{split}
&(Td^\vee(-s_{j_1*}\sO_{\D_{j_1}})-1)(Td^\vee(-s_{j_2*}\sO_{\D_{j_2}})-1)= 0\quad \text{for } 1\leq j_1< j_2\leq n,\\ 
&(Td^\vee(-s_{j*}\sO_{\D_j})-1)(Td^\vee(L_{n+1})-1)=0\quad \text{for } 1\leq j\leq n,\\
&(Td^\vee(-s_{j*}\sO_{\D_j})-1)(Td^\vee(-i_*\sO_\Z)-1)=0\quad \text{for } 1\leq j\leq n,\\
&(Td^\vee(-i_*\sO_\Z)-1)(Td^\vee(L_{n+1})-1)=0.\\
\end{split}
\end{equation*}
Equivalently,
\begin{equation*}
\begin{split}
&Td^\vee(-s_{j_1*}\sO_{\D_{j_1}}-s_{j_2*}\sO_{\D_{j_2}})-1=(Td^\vee(-s_{j_1*}\sO_{\D_{j_1}})-1)+(Td^\vee(-s_{j_2*}\sO_{\D_{j_2}})-1)\quad \text{for } 1\leq j_1< j_2\leq n,\\ 
&Td^\vee(-s_{j*}\sO_{\D_j}+L_{n+1})-1=(Td^\vee(-s_{j*}\sO_{\D_j})-1)+(Td^\vee(L_{n+1})-1)\quad \text{for } 1\leq j\leq n,\\
&Td^\vee(-s_{j*}\sO_{\D_j}-i_*\sO_\Z)-1=(Td^\vee(-s_{j*}\sO_{\D_j})-1)+(Td^\vee(-i_*\sO_\Z)-1)\quad \text{for } 1\leq j\leq n,\\
&Td^\vee(-i_*\sO_\Z+L_{n+1})-1=(Td^\vee(-i_*\sO_\Z)-1)+(Td^\vee(L_{n+1})-1).\\
\end{split}
\end{equation*}
Using these equations repeatedly, we find 
\begin{equation*}
\begin{split}
Td^\vee(\Omega_f)-1&=Td^\vee\left(L_{n+1}+\sum_j (-s_{j*}\sO_{\D_j})-i_*\sO_\Z\right)-1\\
&=(Td^\vee(L_{n+1})-1)+\sum_j (Td^\vee(s_{j*}\sO_{\D_j})^{-1}-1)+(Td^\vee(i_*\sO_\Z)^{-1}-1).
\end{split}
\end{equation*}
Hence the contribution from the main stratum is
\begin{equation*}
\begin{split}
\{f_*(ch(E)Td^\vee(L_{n+1}))&+\sum_j f_*(ch(E)(Td^\vee(s_{j*}\sO_{\D_j})^{-1}-1))\\
&+f_*(ch(E)(Td^\vee(i_*\sO_\Z)^{-1}-1))\}\cap[\Mbar_{g,n}(\X,d)]^{vir}.
\end{split}
\end{equation*}

The term $Td^\vee(s_{j*}\sO_{\D_j})^{-1}-1$ is computed as follows:
Consider the exact sequence 
\begin{equation}\label{divseq}
0\to\sO(-\D_j)\to\sO\to s_{j*}\sO_{\D_j}\to 0.
\end{equation}
Note that $s_j^*(-\D_j)=c_1(N_j^\vee)$ with $N_j^\vee$ the conormal bundle of $\D_j\to \Mbar_{g,n+1}(\X,d)'$. It follows that 
\begin{equation*}
\begin{split}
Td^\vee(s_{j*}\sO_{\D_j})^{-1}-1 &= Td^\vee(\sO(-\D_j))-1=\sum_{r\geq 1}\frac{B_r}{r!}(-\D_j)^r\\
&= -s_{j*}\sum_{r\geq 1}\frac{B_r}{r!}(c_1(N_j^\vee))^{r-1}=-s_{j*}\left[\frac{Td^\vee(N_j^\vee)}{c_1(N_j^\vee)}\right]_+.
\end{split}
\end{equation*}
Here and henceforth the symbol $[\cdot]_+$ denotes power series truncation, which removes terms containing negative powers of cohomology classes.

The term $Td^\vee(i_*\sO_\Z)^{-1}-1$ is computed as follows: Let $\phi: \tilde{\Z}\to \Z$ be the 
double cover of $\Z$ consisting of nodes and choices of a branch at each node. $\tilde{\Z}$ is a disjoint union of open-and-closed substacks of the form $\Mbar_{g-1,n+\{+,-\}}(\X,d)\times_{I\X\times I\X} I\X$ or of the form $\Mbar_+\times_{I\X}\Mbar_-$, where $\Mbar_\pm=\Mbar_{g_\pm,n_\pm+1}(\X,d_\pm)$ such that $g_++g_-=g, n_++n_-=n, d_++d_-=d$ is an ordered splitting of $g,n,d$. This follows from the fact that $\Z$ is the universal gerbe of nodes over $f(\Z)$ (c.f. \cite{AGV}, Proposition 5.2.1).

Let $L_+$ be the line bundle on $\Mbar_+$ whose fiber at an orbifold stable map is the cotangent space\footnote{This is {\em not} the cotangent space on the coarse curve.} at the marked point of gluing. The line bundle $L_-$ on $\Mbar_-$ is similarly defined. On $\Mbar_{g-1,n+\{+,-\}}(\X,d)$ the cotangent line bundles at marked points $+$ and $-$ are also denoted by $L_+$ and $L_-$.

By \cite{M} Lemma 5.1, there is a polynomial $P$ such that $$Td^\vee(i_*\sO_\Z)^{-1}-1=i_*P(c_1(N),c_2(N)),$$ where $N$ is the normal bundle of $\Z\subset\Mbar_{g,n+1}(\X,d)'$. Thus we have 
$$Td^\vee(i_*\sO_\Z)^{-1}-1=\frac{1}{2}i_*\phi_*P(c_1(\phi^*N),c_2(\phi^*N)).$$ Denote $\iota=i\circ\phi$. Using $\phi^*N=L_+^\vee\oplus L_-^\vee$ and the expression of $P$ in \cite{M}, page 303, we find
\begin{equation}\label{nodal_term_main_comp}
\begin{split}
Td^\vee(i_*\sO_\Z)^{-1}-1&=\frac{1}{2}\iota_*\left(\sum_{s\geq 2}\frac{B_s}{s!}\sum_{a+b=s-2}(-1)^a\psi_+^a\psi_-^b\right)\\
&=\frac{1}{2}\iota_*\left(\frac{1}{\psi_++\psi_-}\left(\frac{1}{e^{\psi_+}-1}-\frac{1}{\psi_+}+\frac{1}{2}+\frac{1}{e^{\psi_-}-1}-\frac{1}{\psi_-}+\frac{1}{2}\right)\right)\\
&=\frac{1}{2}\iota_*\left[\frac{1}{\psi_++\psi_-}\left(\frac{Td^\vee(L_+)}{\psi_+}+\frac{Td^\vee(L_-)}{\psi_-}\right)\right]_+.
\end{split}
\end{equation}
Here $\psi_\pm=c_1(L_\pm)$.

Therefore the contribution from the main stratum is
\begin{equation*}
\begin{split}
&\quad f_*(ch(E)Td^\vee(L_{n+1}))\cap[\Mbar_{g,n}(\X,d)]^{vir}\\
&-\sum_j f_*s_{j*}\left(ch(s_j^*E)\left[\frac{Td^\vee(N_j^\vee)}{c_1(N_j^\vee)}\right]_+\right)\cap[\Mbar_{g,n}(\X,d)]^{vir}\\
&+\frac{1}{2}(f\circ\iota)_*\left(ch(\iota^*E)\left[\frac{1}{\psi_++\psi_-}\left(\frac{Td^\vee(L_+)}{\psi_+}+\frac{Td^\vee(L_-)}{\psi_-}\right)\right]_+\right)\cap[\Mbar_{g,n}(\X,d)]^{vir}.
\end{split}
\end{equation*}

The contribution to $ch(f_{(i_1,...,i_n)*}E)\cap[\Mbar_{g,n}(\X,d;i_1,...,i_n)]^{vir}$ from the main stratum can be found by restricting the above to $\Mbar_{g,n}(\X,d;i_1,...,i_n)$. It is the sum of the following three terms:
\begin{equation}\label{co0}
f_{(i_1,...,i_n)*}(ch(E)Td^\vee(L_{n+1}))\cap[\Mbar_{g,n}(\X,d;i_1,...,i_n)]^{vir},
\end{equation}
\begin{equation}\label{co1}
-\sum_j f_{(i_1,...,i_n)*}s_{j,(i_1,...,i_n)*}\left(ch(s_{j,(i_1,...,i_n)}^*E)\left[\frac{Td^\vee(N_j^\vee)}{c_1(N_j^\vee)}\right]_+\right)\cap[\Mbar_{g,n}(\X,d;i_1,...,i_n)]^{vir},
\end{equation}
\begin{equation}\label{co2}
\frac{1}{2}(f_{(i_1,...,i_n)}\circ\iota_{(i_1,...,i_n)})_*\left(ch(\iota_{(i_1,...,i_n)}^*E)\left[\frac{1}{\psi_++\psi_-}\left(\frac{Td^\vee(L_+)}{\psi_+}+\frac{Td^\vee(L_-)}{\psi_-}\right)\right]_+\right)\cap[\Mbar_{g,n}(\X,d;i_1,...,i_n)]^{vir}.
\end{equation}
Here the subscript $_{(i_1,...,i_n)}$ indicates the restriction to $\Mbar_{g,n}(\X,d;i_1,...,i_n)$. We call (\ref{co0}) the codim-0 term, (\ref{co1}) the codim-1 term, and (\ref{co2}) the codim-2 term.
\begin{remark}
Consider the stack $\Mbar_+\times_{I\X}\Mbar_-$ parametrizing maps whose domains consist of two parts separated by a distinguished node\footnote{By definition, a section of the gerbe at the distinguished node is part of the data in this moduli problem.}. If the order of the automorphism group of this node is $r$, then $r\psi_+=\bpsi_+$ is the first Chern class of the line bundle whose fiber is the cotangent line of the {\em coarse curve} at the marked point of gluing. Similarly $r\psi_-=\bpsi_-$. The same statements hold for $L_\pm$ on $\Mbar_{g-1,n+\{+,-\}}(\X,d)$.
\end{remark}

\subsubsection{Marked points}
We compute the contribution from the divisors formed by marked points. Let $s_{j,(i_1,...,i_n)}:\D_{j,(i_1,...,i_n)}\to \Mbar_{g,n+1}(\X,d;i_1,...,i_n,0)$ be the divisor of the $j$-th marked point. We know that $\D_{j,(i_1,...,i_n)}\simeq \Mbar_{g,n}(\X,d;i_1,...,i_n)\times B\mu_{r_{i_j}}$ and the diagram  
\begin{equation*}
\begin{CD}
\Mbar_{g,n}(\X,d;i_1,...,i_n)\times B\mu_{r_{i_j}} @>  >>  \X \\
@V{ }VV  \\
\Mbar_{g,n}(\X,d;i_1,...,i_n).
\end{CD}
\end{equation*}
defined by the restriction of the universal orbifold stable map is equivalent to the evaluation map $$ev_j: \Mbar_{g,n}(\X,d;i_1,...,i_n)\to \X_{i_j}.$$ Also, the generator $\mathfrak{u}_{r_{i_j}}\in \mu_{r_{i_j}}$ acts on the conormal bundle $N_j^\vee$ with eigenvalue $\zeta_{r_{i_j}}^{-1}$.

The locus $\D_{j,(i_1,...,i_n)}$ contributes components of $I\Mbar_{g,n+1}(\X,d;i_1,...,i_n,0)$ which are mapped to $\Mbar_{g,n}(\X,d;i_1,...,i_n)$. These components are
$$\Mbar_{g,n}(\X,d;i_1,...,i_n)\times (IB\mu_{r_{i_j}}\setminus B\mu_{r_{i_j}})=:\coprod_{1\leq l\leq r_{i_j}-1}\D_{j,(i_1,...,i_n)}(l)$$
where $\D_{j,(i_1,...,i_n)}(l)$ is defined as follows. The inertia stack $IB\mu_{r_{i_j}}$ can be described as $$IB\mu_{r_{i_j}}=\coprod_{0\leq k\leq r_{i_j}-1}[\text{Spec}\,\com/C(\mathfrak{u}_{r_{i_j}}^k)].$$ Define $$\D_{j,(i_1,...,i_n)}(l):=\Mbar_{g,n}(\X,d;i_1,...,i_n)\times [\text{Spec}\,\com/C(\mathfrak{u}_{r_{i_j}}^l)]\simeq \Mbar_{g,n}(\X,d;i_1,...,i_n)\times B\mu_{r_{i_j}}.$$ These components arise in the following way. The auotomorphism group of an object of $\D_{j,(i_1,...,i_n)}$ splits as a product $Aut\times  \mu_{r_{i_j}}$ where the first factor $Aut$ is the automorphism group of the correponding object in $\Mbar_{g,n}(\X,d;i_1,...,i_n)$. The components $\D_{j,(i_1,...,i_n)}(l), 1\leq l\leq r_{i_j}-1$ correspond to taking the identity element of the factor $Aut$ and elemenets $\mathfrak{u}_{r_{i_j}}^l, 1\leq l\leq r_{i_j}-1$ in the second factor $\mu_{r_{i_j}}$.

By Lemma \ref{cot} and the exact sequence (\ref{divseq}), we see that the pullback of $T_{f_{(i_1,...,i_n)}}$ to $\D_{j,(i_1,...,i_n)}(l)$ has trivial invariant part, and the moving part is the pullback of $N_j$ to $\D_{j,(i_1,...,i_n)}(l)$.

The restriction $E|_{\D_{j,(i_1,...,i_n)}}$ is decomposed into a direct sum $\oplus_{0\leq k\leq r_{i_j}-1} E^{(k)}$ of $\mathfrak{u}_{r_{i_j}}$-eigenbundles, where $E^{(k)}$ has $\mathfrak{u}_{r_{i_j}}$-eigenvalue $\zeta_{r_{i_j}}^k$ and $\zeta_{r_{i_j}}=\exp(2\pi\sqrt{-1}\frac{1}{r_{i_j}})$. Let $P_l:\D_{j,(i_1,...,i_n)}(l)\to\D_{j,(i_1,...,i_n)}$ be the projection. Then we have
$$ch(\rho(P_l^*E|_{\D_{j,(i_1,...,i_n)}}))=\sum_{0\leq k\leq r_{i_j}-1} \zeta_{r_{i_j}}^{kl}ch(P_l^*(E^{(k)})).$$
So the contribution from $\D_{j,(i_1,...,i_n)}(l)$ is 
$$(f_{(i_1,...,i_n)}\circ s_{j,(i_1,...,i_n)}\circ P_l)_*\left(\frac{\sum_{0\leq k\leq r_{i_j}-1} \zeta_{r_{i_j}}^{kl}ch(P_l^*(E^{(k)}))}{1-\zeta_{r_{i_j}}^{-l}ch(P_l^*N_j^\vee)}\right)\cap[\Mbar_{g,n}(\X,d;i_1,...,i_n)]^{vir}.$$

Let $\gamma_l:\Mbar_{g,n}(\X,d;i_1,...,i_n)\to \D_{j,(i_1,...,i_n)}(l)\simeq \Mbar_{g,n}(\X,d;i_1,...,i_n)\times B\mu_{r_{i_j}}$ be the map such that the map to the first factor is the identity and the map to the second factor corresponds to the trivial $\mu_{r_{i_j}}$-bundle. We have $$\gamma_{l*}\gamma_l^*=r_{i_j}\cdot id \quad \text{and } f_{(i_1,...,i_n)}\circ s_{j,(i_1,...,i_n)}\circ P_l\circ \gamma_l=id.$$ Hence we can write the contribution from $\D_{j,(i_1,...,i_n)}(l)$ as
$$\frac{1}{r_{i_j}}\left(\frac{\sum_{0\leq k\leq r_{i_j}-1} \zeta_{r_{i_j}}^{kl}ch(\gamma_l^*P_l^*(E^{(k)}))}{1-\zeta_{r_{i_j}}^{-l}ch(\gamma_l^*P_l^*N_j^\vee)}\right)\cap[\Mbar_{g,n}(\X,d;i_1,...,i_n)]^{vir}.$$ The following Lemma is straightforward.
\begin{lemma}\label{bundle}
\hfill
\begin{enumerate}

\item For $E=ev_{n+1}^*F$, we have $\gamma_l^*P_l^*(E^{(k)})=ev_j^*(F_{i_j}^{(k)})$.

\item $\gamma_l^*P_l^*N_j^\vee=L_j$.
\end{enumerate}
\end{lemma}
\begin{proof}
The second statement follows from the definition. We prove the first statement. Let $S\to \Mbar_{g,n}(\X,d)$ be a morphism and $S\leftarrow\C\to \X$ the corresponding orbifold stable map. Restricting to the divisor of the $j$-th marked point yields morphisms $$S\overset{p}{\leftarrow} S\times B\mu_{r_{i_j}}\overset{\rho}{\to} \X.$$ By the description of the inertia stack $I\X$ in Remark \ref{inerem} (i), these morphisms correspond to a morphism $\tilde{\rho}:S\to \X_{i_j}$. Consider the component $B\mu_{r_{i_j}}\simeq [\text{Spec}\,\com/C(\mathfrak{u}_{r_{i_j}}^l)]\subset IB\mu_{r_{i_j}}$ and let $\pi_l:S\times [\text{Spec}\,\com/C(\mathfrak{u}_{r_{i_j}}^l)]\to S\times B\mu_{r_{i_j}}$ be the projection. Let $\gamma :S\to S\times [\text{Spec}\,\com/C(\mathfrak{u}_{r_{i_j}}^l)]$ be the section of $p\circ \pi_l$ such that the map to the first factor is the identity and the map to the second factor corresponds to the trivial $\mu_{r_{i_j}}$-bundle. Let $(\rho^*F)^{(k)}$ be the eigen sub-bundle of $\rho^*F$ on which $\mathfrak{u}_{r_{i_j}}$ acts with eigenvalue $\zeta_{r_{i_j}}^k$. To prove the first statement it suffices to prove $$\gamma^*\pi_l^*((\rho^*F)^{(k)})=\tilde{\rho}^*(F_{i_j}^{(k)}).$$ This is obtained immediately from (\ref{eigen-bdle-description}) by pulling back via $\tilde{\rho}$. 


\end{proof}
Therefore the contribution from $\D_{j,(i_1,...,i_n)}(l)$ can be written as
$$\frac{1}{r_{i_j}}\left(\frac{\sum_{0\leq k\leq r_{i_j}-1} \zeta_{r_{i_j}}^{kl}ch(ev_j^*(F_{i_j}^{(k)}))}{1-\zeta_{r_{i_j}}^{-l}ch(L_j)}\right)\cap[\Mbar_{g,n}(\X,d;i_1,...,i_n)]^{vir}.$$

The contributions from $\D_{j,(i_1,...,i_n)}(1),...,\D_{j,(i_1,...,i_n)}(r_{i_j}-1)$ add up to
\begin{equation*}
\begin{split}
&\quad \sum_{1\leq l\leq r_{i_j}-1}\frac{1}{r_{i_j}}\left(\frac{\sum_{0\leq k\leq r_{i_j}-1} \zeta_{r_{i_j}}^{kl}ch(ev_j^*(F_{i_j}^{(k)}))}{1-\zeta_{r_{i_j}}^{-l}ch(L_j)}\right)\cap[\Mbar_{g,n}(\X,d;i_1,...,i_n)]^{vir}\\
&=\frac{1}{r_{i_j}}\sum_{0\leq k\leq r_{i_j}-1} ch(ev_j^*(F_{i_j}^{(k)})) \sum_{1\leq l\leq r_{i_j}-1} \frac{\zeta_{r_{i_j}}^{kl}}{1-\zeta_{r_{i_j}}^{-l}e^{c_1(L_j)}}\cap[\Mbar_{g,n}(\X,d;i_1,...,i_n)]^{vir}.
\end{split}
\end{equation*}
For each $k$ with $0\leq k < r_{i_j}$ we have
$$\sum_{1\leq l\leq r_{i_j}-1} \frac{\zeta_{r_{i_j}}^{kl}}{1-\zeta_{r_{i_j}}^{-l}e^{c_1(L_j)}}
=\frac{r_{i_j}e^{kc_1(L_j)}}{1-e^{r_{i_j}c_1(L_j)}}-\frac{1}{1-e^{c_1(L_j)}}.$$
Using $\gamma_{0*}\gamma_0^*=r_{i_j}\cdot id$ we can rewrite the part of codim-1 term (\ref{co1}) that comes from $\D_{j,(i_1,...,i_n)}$ as
\begin{equation*}
\begin{split}
&\quad -\frac{1}{r_{i_j}}\sum_{0\leq k\leq r_{i_j}-1} ch(ev_j^*(F_{i_j}^{(k)})) \sum_{n\geq 1}\frac{B_n}{n!}c_1(L_j)^{n-1}\cap[\Mbar_{g,n}(\X,d;i_1,...,i_n)]^{vir}\\
&=-\frac{1}{r_{i_j}}\sum_{0\leq k\leq r_{i_j}-1} ch(ev_j^*(F_{i_j}^{(k)}))\frac{1}{c_1(L_j)}\left(\frac{c_1(L_j)}{e^{c_1(L_j)}-1}-1\right)\cap[\Mbar_{g,n}(\X,d;i_1,...,i_n)]^{vir}.
\end{split}
\end{equation*}
Here we use $\sum_{n\geq 1}\frac{B_n}{n!}x^{n-1}=\frac{1}{x}(\frac{x}{e^x-1}-1)$. Combining this part of codim-1 term and contributions from $\D_{j,(i_1,...,i_n)}(1)$,..., $\D_{j,(i_1,...,i_n)}(r_{i_j}-1)$, we find that their sum is equal to
$$\frac{1}{r_{i_j}}\sum_{0\leq k\leq r_{i_j}-1} ch(ev_j^*(F_{i_j}^{(k)})) \left(\frac{r_{i_j}e^{kc_1(L_j)}}{1-e^{r_{i_j}c_1(L_j)}}+\frac{1}{c_1(L_j)}\right)\cap[\Mbar_{g,n}(\X,d;i_1,...,i_n)]^{vir}.$$
Using the definition of Bernoulli polynomials, we see that this is
$$-\sum_{0\leq k\leq r_{i_j}-1} ch(ev_j^*(F_{i_j}^{(k)}))\sum_{n\geq 1}\frac{B_n(k/r_{i_j})}{n!}(r_{i_j}c_1(L_j))^{n-1}\cap[\Mbar_{g,n}(\X,d;i_1,...,i_n)]^{vir}$$
\begin{equation}\label{mark}
=-\sum_{n\geq 1}\frac{\sum_{0\leq k\leq r_{i_j}-1} ch(ev_j^*(F_{i_j}^{(k)}))B_n(k/r_{i_j})}{n!}\bpsi_j^{n-1}\cap[\Mbar_{g,n}(\X,d;i_1,...,i_n)]^{vir}.
\end{equation}
Here we also use $\bpsi_j=r_{i_j}\psi_j$, which follows from the fact that $L_j^{\otimes r_{i_j}}=\pi_n^*L_j$.

\subsubsection{Nodes}
We proceed to compute the contributions from the locus of nodes in a similar fashion. Let $$\phi_{r,(i_1,...,i_n)}:\tilde{\Z}_{r,(i_1,...,i_n)}\to \Z_{r,(i_1,...,i_n)}$$ be the double covering of $\Z_{r,(i_1,...,i_n)}$ consisting of nodes and choices of a branch at each node and $$\iota_{r,(i_1,...,i_n)}: \tilde{\Z}_{r,(i_1,...,i_n)}\to \Mbar_{g,n+1}(\X,d;i_1,...,i_n,0)$$ be $\phi_{r,(i_1,...,i_n)}$ followed by the inclusion. 

The components $\tilde{\Z}_{r,(i_1,...,i_n)}(1)$,..., $\tilde{\Z}_{r,(i_1,...,i_n)}(r-1)$ of $I\tilde{\Z}_{r,(i_1,...,i_n)}$ which are mapped to $\Mbar_{g,n}(\X,d;i_1,...,i_n)$ can be defined similarly as $\D_{j,(i_1,...,i_n)}(l)$. Since $\tilde{\Z}$ can be identified with a disjoint union of the stack $$\Mbar_{g-1,n+\{+,-\}}(\X,d)\times_{I\X\times I\X} I\X$$ and stacks of the form $\Mbar_+\times_{I\X}\Mbar_-$, each $\tilde{\Z}_{r,(i_1,...,i_n)}(l)$ is isomorphic to $\tilde{\Z}_{r,(i_1,...,i_n)}$. Let $$P_l:\tilde{\Z}_{r,(i_1,...,i_n)}(l)\to \tilde{\Z}_{r,(i_1,...,i_n)}$$ be the projection, and $\gamma_l$ an inverse of $P_l$. Note that $\gamma_{l*}\gamma_l^*=id$.

By the Koszul complex $$0\to\sO(L_+\otimes L_-)\to\sO(L_+)\oplus\sO(L_-)\to\sO\to\sO_{\tilde{\Z}}\to 0 $$ and Lemma \ref{cot}, we see that the invariant part of the pullback of $T_{f_{(i_1,...,i_n)}}$ to $\tilde{\Z}_{r,(i_1,...,i_n)}(l)$ is the sum of a trivial bundle $\sO$, and $-\sO$, and $-(L_+\otimes L_-)^\vee$. The moving part is $(L_+^\vee\oplus L_-^\vee)$. 

The contribution from $\tilde{\Z}_{r,(i_1,...,i_n)}(l)$ is
$$\frac{1}{2}(f_{(i_1,...,i_n)}\circ\iota_{r,(i_1,...,i_n)})_*P_{l*}\gamma_{l*}\left(\frac{\sum_k \zeta_r^{lk}ch(\gamma_l^*P_l^*(E^{(k)}))Td(-(L_+\otimes L_-)^\vee)}{1-\zeta_r^{-l}ch(L_+)-\zeta_r^l ch(L_-)+ch(\wedge^2(L_+\oplus L_-))}\right).$$
Put $\psi_+=c_1(L_+), \psi_-=c_1(L_-)$. Note that $Td(-(L_+\otimes L_-)^\vee)=Td^\vee(-(L_+\otimes L_-))=1/Td^\vee(L_+\otimes L_-)$. We have
\begin{equation*}
\begin{split}
&\frac{Td(-(L_+\otimes L_-)^\vee)}{1-\zeta_r^{-l}ch(L_+)-\zeta_r^l ch(L_-)+ch(\wedge^2(L_+\oplus L_-))}\\
=&\frac{e^{\psi_++\psi_-}-1}{(\psi_++\psi_-)(1-\zeta_r^{-l}e^{\psi_+}-\zeta_r^l e^{\psi_-}+e^{\psi_++\psi_-})}\\
=&\frac{e^{\psi_++\psi_-}-1}{(\psi_++\psi_-)(1-\zeta_r^{-l}e^{\psi_+})(1-\zeta_r^l e^{\psi_-})}=\frac{1}{\psi_++\psi_-}\left(1+\frac{1}{\zeta_r^{-l}e^{\psi_+}-1}+\frac{1}{\zeta_r^l e^{\psi_-}-1}\right).
\end{split}
\end{equation*}

Also, for $0< k< r$, 
$$\sum_{l=1}^{r-1}\frac{\zeta_r^{kl}}{\zeta_r^{-l}e^x-1}=\frac{re^{kx}}{e^{rx}-1}-\frac{1}{e^x-1},\quad
\sum_{l=1}^{r-1}\frac{\zeta_r^{kl}}{\zeta_r^l e^x-1}=\frac{re^{(r-k)x}}{e^{rx}-1}-\frac{1}{e^x-1}.$$ 
And 
\[\sum_{l=1}^{r-1}\frac{1}{\zeta_r^{-l}e^x-1}=\sum_{l=1}^{r-1}\frac{1}{\zeta_r^l e^x-1}=\frac{r}{e^{rx}-1}-\frac{1}{e^x-1},\quad \sum_{l=1}^{r-1}\zeta_r^{kl}=\left\{
 \begin{array}{rr}
 -1, \quad &\mbox{$k\neq 0$}\\
 r-1,\quad &\mbox{$k=0$}.
 \end{array}\right.   \]

Therefore 
\begin{equation}\label{exponential_sums}
\begin{split}
 &\sum_{l=1}^{r-1}\left(\zeta_r^{kl}+\frac{\zeta_r^{kl}}{\zeta_r^{-l}e^{\psi_+}-1}+\frac{\zeta_r^{kl}}{\zeta_r^l e^{\psi_-}-1}\right)\\
=&\left\{
 \begin{array}{rr}
 \frac{re^{k\psi_+}}{e^{r\psi_+}-1}-\frac{1}{e^{\psi_+}-1}+\frac{re^{(r-k)\psi_-}}{e^{r\psi_-}-1}-\frac{1}{e^{\psi_-}-1} -1, \quad &\mbox{$k\neq 0$}\\
 \frac{r}{e^{r\psi_+}-1}-\frac{1}{e^{\psi_+}-1}+\frac{r}{e^{r\psi_-}-1}-\frac{1}{e^{\psi_-}-1}+ r-1,\quad &\mbox{$k=0$}.
 \end{array}\right.   
\end{split}
\end{equation}

We have $\gamma_l^*P_l^*E^{(k)}=ev_{node}^*(F^{(k)})$, which is similar to Lemma \ref{bundle} (see Appendix \ref{vbdle} for the definition of $ev_{node}$).  What we need to do now is to combine the part of codim-2 term (\ref{co2}) from $\tilde{\Z}_{r,(i_1,...,i_n)}$ and contributions of $\tilde{\Z}_{r,(i_1,...,i_n)}(1),...,\tilde{\Z}_{r,(i_1,...,i_n)}(r-1)$. First note that the term $ch(\iota_{(i_1,...,i_n)}^*E)$ in (\ref{co2}) breaks into a sum of terms $ch(ev_{node}^*((q^*F)^{(k)}))$ for $0\leq k <r$. The term in (\ref{co2}) corresponding to $k$ is (the pushforward of) $ch(ev_{node}^*((q^*F)^{(k)}))$ multiplied by $\frac{1}{\psi_++\psi_-}$ and $$\frac{1}{e^{\psi_+}-1}-\frac{1}{\psi_+}+\frac{1}{2}+\frac{1}{e^{\psi_-}-1}-\frac{1}{\psi_-}+\frac{1}{2},$$
see (\ref{nodal_term_main_comp}). Adding this to (\ref{exponential_sums}), we get for $k\neq 0$
\begin{equation*}
\begin{split}
&\quad \frac{re^{k\psi_+}}{e^{r\psi_+}-1}-\frac{1}{\psi_+}+\frac{re^{(r-k)\psi_-}}{e^{r\psi_-}-1}-\frac{1}{\psi_-}\\
&=r\left(\frac{e^{\frac{k}{r}r\psi_+}}{e^{r\psi_+}-1}-\frac{1}{r\psi_+}+\frac{e^{\frac{r-k}{r}r\psi_-}}{e^{r\psi_-}-1}-\frac{1}{r\psi_-}\right)\\
&=r\sum_{n\geq 1}\left(\frac{B_n(k/r)}{n!}(r\psi_+)^{n-1}+\frac{B_n(1-k/r)}{n!}(r\psi_-)^{n-1}\right),
\end{split}
\end{equation*}

and for $k=0$
\begin{equation*}
\begin{split}
&\quad \frac{r}{e^{r\psi_+}-1}-\frac{1}{\psi_+}+\frac{r}{e^{r\psi_-}-1}-\frac{1}{\psi_-}+ r\\
&=r\left(\frac{1}{e^{r\psi_+}-1}-\frac{1}{r\psi_+}+\frac{1}{2}+\frac{1}{e^{r\psi_-}-1}-\frac{1}{r\psi_-}+\frac{1}{2}\right)\\
&=r\sum_{n\geq 2}\left(\frac{B_n}{n!}(r\psi_+)^{n-1}+\frac{B_n}{n!}(r\psi_-)^{n-1}\right)\\
&=r\sum_{n\geq 1}\left(\frac{B_n}{n!}(r\psi_+)^{n-1}+\frac{B_n(1)}{n!}(r\psi_-)^{n-1}\right).
\end{split}
\end{equation*}

By these calculation it follows that the combined contribution is the following expression capped with the virtual class:
 
\begin{equation*}
\begin{split}
&{\footnotesize
\frac{r^2}{2}(f_{(i_1,...,i_n)}\circ\iota_{r,(i_1,...,i_n)})_*\sum_{n\geq 1}\frac{1}{n!}\frac{1}{r\psi_++r\psi_-}\sum_{0\leq l<r} ch(ev_{node}^*((q^*F)^{(l)}))(B_n(\frac{l}{r})(r\psi_+)^{n-1}+B_n(1-\frac{l}{r})(r\psi_-)^{n-1})}\\
&=\frac{r^2}{2}(f_{(i_1,...,i_n)}\circ\iota_{r,(i_1,...,i_n)})_*\sum_{n\geq 1}\frac{1}{n!}\left[\sum_l ch(ev_{node}^*((q^*F)^{(l)}))B_n(l/r)\right]\sum_{a+b=n-2}(-1)^b(r\psi_+)^a(r\psi_-)^b\\
&=\frac{r^2}{2}(f\circ\iota)_*\sum_{n\geq 2}\frac{1}{n!}\left[\sum_l ch(ev_{node}^*((q^*F)^{(l)}))B_n(l/r)\right]\frac{(\bpsi_+)^{n-1}+(-1)^n(\bpsi_-)^{n-1}}{\bpsi_++\bpsi_-}.
\end{split}
\end{equation*}
Here we use $r\psi_\pm=\bpsi_\pm$. Note that we rewrite $\frac{1}{\psi_++\psi_-}$ as $r\cdot\frac{1}{r\psi_++r\psi_-}$, which gives a factor of $r$.

Combining all together, we find
\begin{equation}\label{chern_char}
\begin{split}
&\quad ch(f_*ev_{n+1}^*F)\cap[\Mbar_{g,n}(\X,d)]^{vir}\\
&=f_*(ch(ev^*F)Td^\vee(L_{n+1}))\cap[\Mbar_{g,n}(\X,d)]^{vir}\\
&\quad -\sum_{i=1}^n\sum_{m\geq 1}\frac{ev_i^*A_m}{m!}(\bpsi_i)^{m-1}\cap[\Mbar_{g,n}(\X,d)]^{vir}\\
&\quad +\frac{1}{2}(f\circ\iota)_*\sum_{m\geq 2}\frac{1}{m!}r_{node}^2(ev_{node}^*A_m)\left(\frac{(\bpsi_+)^{m-1}+(-1)^m(\bpsi_-)^{m-1}}{\bpsi_++\bpsi_-}\right)\cap[\Mbar_{g,n}(\X,d)]^{vir}.
\end{split}
\end{equation}

\subsection{Finding the differential equation}\label{find_de}
We begin with the following splitting property of the virtual fundamental classes, which will be used in the calculations. Let $\mathfrak{M}^{tw}_{g, n}$ be the (Artin) stack of twisted curves of genus $g$ with $n$ marked gerbes (not trivialized). First consider the case of separating nodes. Let $$\mathfrak{D}^{tw}(g_+;n_+|g_-; n_-):=\coprod_{\{1,...,n\}=A\cup B, |A|=n_+, |B|=n_-}\mathfrak{D}^{tw}(g_+; A|g_-; B),$$
where the right-hand side is defined as in \cite{AGV2}, Section 5.1. There is a natural forgetful map $\Mbar_{g,n}(\X, d)\to \mathfrak{M}^{tw}_{g,n}$ and a natural gluing map $gl: \mathfrak{D}^{tw}(g_+;n_+|g_-; n_-)\to \mathfrak{M}^{tw}_{g,n}$ as defined in \cite{AGV2}, Proposition 5.1.3. Consider the cartesian diagram formed by these maps:
$$\begin{CD}
\mathfrak{D}_{g,n}(\X)@> >>\Mbar_{g,n}(\X,d)\\
@V{}VV @V{}VV\\
\mathfrak{D}^{tw}(g_+;n_+|g_-; n_-)@>gl >> \mathfrak{M}^{tw}_{g,n}.
\end{CD}$$
There is a natural map $$\mathfrak{g}: \bigcup_{d=d_++d_-}\Mbar_{g_+,n_++1}(\X,d_+)\times_{I\X}\Mbar_{g_-,n_-+1}(\X,d_-)\to \mathfrak{D}_{g,n}(\X).$$ This is the universal gerbe over the distinguished node (see \cite{AGV}, Proposition 5.2.1).

Similarly, for non-separating nodes we write $gl$ for the map obtained by gluing the last two marked points. There is a similar cartesian diagram and a similar map $\mathfrak{g}$, which we do not describe explicitly. 

\begin{proposition}\label{spvirclass}
Let $\Mbar_{g_+,n_++1}(\X,d_+)\times_{I\X}\Mbar_{g_-,n_-+1}(\X,d_-)\subset \tilde{\Z}_r\overset{\iota_r}{\to}\Mbar_{g,n}(\X,d)$.
Consider the diagram of gluing,
$$\begin{CD}
\Mbar_{g_+,n_++1}(\X,d_+)\times_{I\X}\Mbar_{g_-,n_-+1}(\X,d_-)@> >>I\X\\
@V{}VV @V{\delta}VV\\
\Mbar_{g_+,n_++1}(\X,d_+)\times \Mbar_{g_-,n_-+1}(\X,d_-)@>ev_+\times \check{ev}_->> I\X\times I\X.
\end{CD}$$
Here $\delta: I\X\to I\X\times I\X$ is the diagonal map, and $\check{ev}_-$ is the composite $$\Mbar_{g_-, n_-+1}(\X, d_-)\overset{ev_-}{\to}I\X\overset{I}{\to}I\X.$$
Then 
$$\sum_{d_++d_-=d}\delta^!([\Mbar_{g_+,n_++1}(\X,d_+)]^w\times [\Mbar_{g_-,n_-+1}(\X,d_-)]^w)=r^2\mathfrak{g}^*(gl^![\Mbar_{g,n}(\X,d)]^w).$$
Similarly, for $\Mbar_{g-1,n+\{+,-\}}(\X,d)\times_{I\X\times I\X}I\X\subset \tilde{\Z}_r\overset{\iota_r}{\to}\Mbar_{g,n}(\X,d)$, we have $$\delta^![\Mbar_{g-1,n+\{+,-\}}(\X,d)]^w=r^2\mathfrak{g}^*(gl^![\Mbar_{g,n}(\X,d)]^w).$$
\end{proposition}
This Proposition is more general than Proposition 5.3.1 of \cite{AGV}. The proof of this Proposition is the same as that of \cite{AGV2}, Proposition 5.3.1, with the straightforward adjustment for weighted virtual classes. In particular, the factor $r^2$ arises since when a stacky node of order $r$ is split into two stacky marked points, each marked point should receive a factor of $r$ in order to get the weighted virtual class. (Note that $r$ should be interpreted as a locally constant function.) 

We now process the term (\ref{chern_term}). According to the GRR calculation (\ref{chern_char}), (\ref{chern_term}) splits into three parts:

{\em Codim-1:}
\begin{equation}\label{codim1}
-\sum_{g,n,d}\frac{Q^d\hbar^{g-1}}{(n-1)!}\<\left(\sum_{m\geq 1}\frac{A_m}{m!}(\bpsi)^{m-1}\right)_k\bt,\bt,...,\bt;\bc(F_{g,n,d})\>_{g,n,d}.
\end{equation}

{\em Codim-2:}
\begin{equation}\label{codim2}
\frac{1}{2}\sum_{g,n,d}\frac{Q^d\hbar^{g-1}}{n!}\<\bt,...,\bt;\left[(f\circ\iota)_*\sum_{m\geq 2}\frac{1}{m!}r_{node}^2ev_{node}^*A_m\frac{\bpsi_+^{m-1}+(-1)^m\bpsi_-^{m-1}}{\bpsi_++\bpsi_-}\right]_k\bc(F_{g,n,d})\>_{g,n,d}.
\end{equation}

{\em Codim-0:}
\begin{equation*}
\begin{split}
&\quad \sum_{g,n,d}\frac{Q^d\hbar^{g-1}}{n!}\<\bt,...,\bt;(f_*(ch(ev^*F)Td^\vee(L_{n+1})))_k\bc(F_{g,n,d})\>_{g,n,d}\\
&=\sum_{g,n,d}\frac{Q^d\hbar^{g-1}}{n!}\<f^*\bt,...,f^*\bt,(ch(ev^*F)Td^\vee(L_{n+1}))_{k+1};\bc(F_{g,n+1,d})\>_{g,n+1,d}^\bullet\\
&=\sum_{g,n,d}\frac{Q^d\hbar^{g-1}}{n!}\<\bt,...,\bt,(ch(F)Td^\vee(L_{n+1}))_{k+1};\bc(F_{g,n+1,d})\>_{g,n+1,d}^\bullet\\
&\quad -\sum_{g,n,d}\frac{Q^d\hbar^{g-1}}{(n-1)!}\<ch_{k+1}(F)\orbcup\left[\frac{\bt(\bpsi)}{\bpsi}\right]_+,\bt,...,\bt;\bc(F_{g,n,d})\>_{g,n,d},
\end{split}
\end{equation*}
where we have used Lemma \ref{split}. This is equal to the sum of the following four terms:
\begin{equation}\label{3}
\sum_{g,n,d}\frac{Q^d\hbar^{g-1}}{(n-1)!}\<\bt,...,\bt,(ch(F)Td^\vee(L))_{k+1};\bc(F_{g,n,d})\>_{g,n,d}^\bullet
\end{equation}
\begin{equation}\label{4}
-\sum_{g,n,d}\frac{Q^d\hbar^{g-1}}{(n-1)!}\<ch_{k+1}(F)\orbcup\left[\frac{\bt(\bpsi)}{\bpsi}\right]_+,\bt,...,\bt;\bc(F_{g,n,d})\>_{g,n,d}
\end{equation}
\begin{equation}\label{5}
-\frac{1}{2\hbar}\<\bt,\bt,(ch(F)Td^\vee(L))_{k+1};\bc(F_{0,3,0})\>_{0,3,0}^\bullet
\end{equation}
\begin{equation}\label{6}
-\<(ch(F)Td^\vee(L))_{k+1};\bc(F_{1,1,0})\>_{1,1,0}^\bullet.
\end{equation}
Here $\<...\>_{...}^\bullet$ denotes invariants defined from moduli spaces of maps with the last marked point untwisted, and we use the property $$f^*\bt(\bpsi_j)=\bt(\bpsi_j)-s_{j*}\left[\frac{\bt(\bpsi_j)}{\bpsi_j}\right]_+.$$ Since $\bpsi_j$ are pulled back from $\Mbar_{g,n}(X,d)$, this follows from the case of schemes (see for instance \cite{CG}).

Observe that, by Lemma \ref{ocup}, $$ch_{k+1}(F)\orbcup\left[\frac{\bt(\bpsi)}{\bpsi}\right]_+=ch_{k+1}(q^*F)\left[\frac{\bt(\bpsi)}{\bpsi}\right]_+=ch_{k+1}(q^*F)\left(\frac{\bt(\bpsi)-t_0}{\bpsi}\right).$$
On a component $\X_i$, we have
\begin{equation*}
\begin{split} 
&\quad \sum_{m\geq 1}\frac{A_m}{m!}z^{m-1}|_{\X_i}=\sum_{m\geq 1}\sum_{0\leq l\leq r_i-1}\frac{ch(F_i^{(l)})B_m(l/r_i)}{m!}z^{m-1}\\
&=\sum_{0\leq l\leq r_i-1}ch(F_i^{(l)})\left(\sum_{m\geq 1}\frac{B_m(l/r_i)}{m!}z^{m-1}\right)=\sum_{0\leq l\leq r_i-1}ch(F_i^{(l)})\left(\frac{e^{\frac{l}{r_i}z}}{e^z-1}-\frac{1}{z}\right),
\end{split}
\end{equation*}
and $ch_{k+1}(q^*F)|_{\X_i}=\sum_{0\leq l\leq r_i-1}ch_{k+1}(F_i^{(l)})$.

For each $l$ we have
$$\left(\frac{e^{\frac{l}{r_i}z}}{e^z-1}-\frac{1}{z}\right)\bt(z)+\frac{\bt(z)-t_0}{z}=\left(\frac{e^{\frac{l}{r_i}z}}{e^z-1}\bt(z)-\frac{t_0}{z}\right)=\left[\frac{e^{\frac{l}{r_i}z}}{e^z-1}\bt(z)\right]_+.$$
Hence
$$\left.\left(\sum_{m\geq 1}\frac{A_m}{m!}\bpsi^{m-1}\right)_k\right|_{\X_i}\bt(\bpsi)+ch_{k+1}(q^*F)|_{\X_i}\left[\frac{\bt(\bpsi)}{\bpsi}\right]_+=\sum_{0\leq l\leq r_i-1}\left[\left(\frac{ch(F_i^{(l)})e^{\frac{l}{r_i}\bpsi}}{e^{\bpsi}-1}\right)_k\bt(\bpsi)\right]_+.$$
Therefore the sum of (\ref{codim1}) and (\ref{4}) is
$$-\sum_{g,n,d}\frac{Q^d\hbar^{g-1}}{(n-1)!}\<\left[\left(\frac{\sum_l ch(F_i^{(l)})e^{\frac{l}{r_i}\bpsi}}{e^z-1}\right)_k\bt\right]_+,\bt,...,\bt;\bc(F_{g,n,d})\>_{g,n,d}.$$
Also, $(ch(F)Td^\vee(L))_{k+1}=[(ch(F)\frac{Td^\vee(L)}{\psi})_k\psi]_+=[(\frac{ch(F)}{e^\psi-1})_k1\psi]_+$.
Hence the sum of (\ref{codim1}), (\ref{3}) and (\ref{4}) is
$$-\sum_{g,n,d}\frac{Q^d\hbar^{g-1}}{(n-1)!}\<\left[\left(\frac{\sum_l ch(F_i^{(l)})e^{\frac{l}{r_i}\bpsi}}{e^z-1}\right)_k(\bt-1\bpsi)\right]_+,\bt,...,\bt;\bc(F_{g,n,d})\>_{g,n,d}.$$
Adding (\ref{der}) to this, we get
\begin{equation}\label{pqterm}
-\sum_{g,n,d}\frac{Q^d\hbar^{g-1}}{(n-1)!}\<\left[\left(\left(\frac{\sum_l ch(F_i^{(l)})e^{\frac{l}{r_i}\bpsi}}{e^z-1}\right)_k+\frac{ch_k((q^*F)^{inv})}{2}\right)\bq(\bpsi)\right]_+,\bt,...,\bt;\bc(F_{g,n,d})\>_{g,n,d}.\end{equation}
The restriction to $\X_i$ of the operator 
\begin{equation}\label{operator}
\sum_{\overset{m+h=k+1}{m,h\geq 0}}\frac{(A_m)_hz^{m-1}}{m!}+\frac{ch_k((q^*F)^{inv})}{2}
\end{equation} 
is
$$
\sum_{\overset{m+h=k+1}{m,h\geq 0}}\frac{(A_m)_hz^{m-1}}{m!}|_{\X_i}+\frac{ch_k((q^*F)^{inv}|_{\X_i})}{2}$$
\begin{equation}\label{operator_1}
=\left(\frac{\sum_l ch(F_i^{(l)})e^{\frac{l}{r_i}z}}{e^z-1}\right)_k+\frac{ch_k((q^*F)^{inv}|_{\X_i})}{2}.
\end{equation}
Note that the operator (\ref{operator}) is infinitesimally symplectic by Corollary \ref{infsymp}.

By \cite{C}, Example 1.3.3.1 (see Appendix \ref{quantized_operator}), the quantization of the $pq$-terms of the quadratic Hamiltonian of (\ref{operator}) applied to $\D_s$ gives (\ref{pqterm}).

It is straightforward to check that the $q^2$-term of the Hamiltonian of the operator (\ref{operator}) only comes from $(A_0)_{k+1}/z=ch_{k+1}(q^*F)/z$. Using statements from Remark \ref{spcase} we can calculate (\ref{5}) directly, the answer is $$-\frac{1}{2\hbar}\int_{I\X}t_0\wedge t_0\wedge ch_{k+1}(q^*F)\wedge \bc((q^*F)^{inv}).$$Then by Appendix \ref{quantized_operator} the quantization of the $q^2$-term yields exactly (\ref{5}).

Now we handle the codim-2 terms (\ref{codim2}), following the approach of \cite{C}. Pulling back to $\Mbar_+\times\Mbar_-$ and $\Mbar_{g-1,n+\{+,-\},d}$ and using Lemma \ref{split} and Proposition \ref{spvirclass}, we express (\ref{codim2}) as 
\begin{equation}\label{codim2_expanded}
\begin{split}
\frac{\hbar}{2}\sum_{g,n,d}\sum_{g_1+g_2=g}\sum_{n_1+n_2=n}\sum_{d_1+d_2=d}\frac{Q^{d_1+d_2}\hbar^{g_1-1+g_2-1}}{n_1!n_2!}\sum_{a,b, c}\<\bt,...,\bt,\frac{\sO_{a,b, c}'\bpsi_+^a}{\sqrt{\bc((q^*F)^{inv})}};\bc(F_{g_1,n_1+1,d_1})\>_{g_1,n_1+1,d_1}\\
\times\<\frac{\sO_{a,b,c,}''\bpsi_-^b}{\sqrt{\bc((q^*F)^{inv})}},\bt,...,\bt;\bc(F_{g_2,n_2+1,d_2})\>_{g_2,n_2+1,d_2}\\
+\frac{\hbar}{2}\sum_{g,n,d}\frac{Q^d\hbar^{g-1-1}}{n!}\sum_{a,b, c}\<\bt,...,\bt,\frac{\sO_{a,b, c}'\bpsi_+^a}{\sqrt{\bc((q^*F)^{inv})}},\frac{\sO_{a,b,c,}''\bpsi_-^b}{\sqrt{\bc((q^*F)^{inv})}}; \bc(F_{g-1,n+2,d})\>_{g-1,n+2,d}.
\end{split}
\end{equation}
Here\footnote{Note that the term $\sum_{m\geq 2}\frac{A_m}{m!}\frac{\bpsi_+^{m-1}+(-1)^m\bpsi_-^{m-1}}{\bpsi_++\bpsi_-}$ belongs to $End(H^*(I\X))[[\bpsi_+, \bpsi_-]]$, which is identified with $H^*(I\X)[[\bpsi_+]]\otimes H^*(I\X)[[\bpsi_-]]$ using the pairing on $H^*(I\X)$.} $$\sum_{a,b}\sO_{a,b}\bpsi_+^a\bpsi_-^b=\left(\sum_{m\geq 2}\frac{A_m}{m!}\frac{\bpsi_+^{m-1}+(-1)^m\bpsi_-^{m-1}}{\bpsi_++\bpsi_-}\right)_{k-1}\cdot (g^{\alpha\beta}\phi_\alpha\otimes \phi_\beta)\in H^*(I\X)[[\bpsi_+]]\otimes H^*(I\X)[[\bpsi_-]],$$
$g^{\alpha\beta}$ is the matrix entry of the inverse of the matrix $(g_{\alpha\beta})$ with $g_{\alpha\beta}=(\phi_\alpha,\phi_\beta)_{orb}$, and we write $\sO_{a,b}\in H^*(I\X)\otimes H^*(I\X)$ in its K\"unneth decomposition:
$$\sO_{a,b}=\sum_c \sO_{a,b, c}'\otimes \sO_{a,b,c}'', \quad \sO_{a,b,c}', \sO_{a,b,c}''\in H^*(I\X).$$

Due to twisted dilaton shift, we have $$\frac{\partial}{\partial q_k^\alpha}=\frac{1}{\sqrt{\bc((q^*F)^{inv})}}\frac{\partial}{\partial t_k^\alpha}.$$
Comparing this with \cite{C}, Example 1.3.3.1 (see Appendix \ref{quantized_operator}), we find that (\ref{codim2_expanded}) coincides with the quantization of the $p^2$-terms of the Hamiltonian of $\sum_{m, h\geq 0, m+h=k+1}\frac{(A_m)_h}{m!}z^{m-1}+ch_k((q^*F)^{inv})/2$ applied to $\D_s$ (note that the Hamiltonian of $(A_0)_{k+1}/z+(A_1)_k+ch_k((q^*F)^{inv})/2$ has no $p^2$-terms).


Putting the above together, we just proved 
$$\frac{\partial\D_s}{\partial s_k}=\left[\left(\sum_{m+h=k+1; m,h\geq 0}\frac{(A_m)_hz^{m-1}}{m!}+\frac{ch_k((q^*F)^{inv})}{2}\right)^\wedge+(\ref{6})\right]\D_s.$$
Note that (\ref{6}) is equal to $C_k$ defined in Section \ref{proof_overview}. This concludes the proof of (\ref{qrrd}).

\appendix

\section{A Grothendieck-Riemann-Roch formula for Stacks}\label{grr}
Let $\X$ and $\Y$ be Deligne-Mumford stacks with quasi-projective coarse moduli spaces.
Let $f:\X\to\Y$ be a proper morphism of Deligne-Mumford stacks. Assume that $f$ factors as 
\begin{equation}\label{lciq} f=g\circ i,\end{equation} where $i:\X\to \sP$ is a closed regular immersion and $g:\sP\to \Y$ is a smooth morphism (not necessarily representable). Define
$$T_f:=[i^*T_{\sP/\Y}]-[N_{\X/\sP}]\in K^0(\X).$$ 
It is easy to show that $T_f$ is independent of the factorization $f=g\circ i$. There is a Grothendieck-Riemann-Roch formula for this kind of morphism, which is due to Toen \cite{T}. We begin with some definitions.

\begin{definition}[\cite{T}]
Define a map $\rho: K^0(I\X)\to K^0(I\X)$ as follows:
If a bundle $F$ on $I\X$ is decomposed into a direct sum $\oplus_\zeta F^{(\zeta)}$ of eigenbundles $F^{(\zeta)}$ with eigenvalue $\zeta$, then 
$$\rho(F):=\sum_\zeta \zeta F^{(\zeta)}\in K^0(I\X).$$
\end{definition}
\begin{definition}[\cite{T}]
Define $\tch: K^0(\X)\to H^*(I\X)$ to be the composite
$$K^0(\X)\overset{q_\X^*}{\longrightarrow} K^0(I\X)\overset{\rho}{\longrightarrow} K^0(I\X)\overset{ch}{\longrightarrow} H^*(I\X),$$
where $q_\X: I\X\to\X$ is the projection and $ch$ is the usual Chern character.
\end{definition}
\begin{definition}
Define an operation $\lambda_{-1}$ in K-theory as follows: for a vector bundle $V$, define $\lambda_{-1}(V):=\sum_{a\geq 0} (-1)^a\Lambda^aV$.
\end{definition}
\begin{definition}[Todd class]
Define $\ttd: K^0(\X)\to H^*(I\X)$ as follows: For a vector bundle $E$ on $\X$, $q_\X^*E$ is decomposed into a direct sum $(q_\X^*E)^{inv}\oplus (q_\X^*E)^{mov}$ where $(q_\X^*E)^{inv}$, the invariant part, is the eigenbundle with eigenvalue $1$, and $(q_\X^*E)^{mov}$, the moving part, is the direct sum of eigenbundles with eigenvalues not equal to $1$. Define 
$$\ttd(E):=\frac{Td((q_\X^*E)^{inv})}{ch(\rho\circ \lambda_{-1}(((q_\X^*E)^{mov})^\vee))}.$$
\end{definition}

The map $\ttd$ satisfies 
$$\ttd(V_1+V_2)=\ttd(V_1)\ttd(V_2),\quad \ttd(V_1-V_2)=\frac{\ttd(V_1)}{\ttd(V_2)}.$$

Recall that a stack has the resolution property if every coherent sheaf is a quotient of a vector bundle (see for instance \cite{To}).
\begin{theorem}[Grothendieck-Riemann-Roch formula \cite{T}]
Let $\X$ and $\Y$ be smooth Deligne-Mumford stacks with quasi-projective coarse moduli spaces and $f:\X\to \Y$ a proper morphism which factors as (\ref{lciq}). Assume that $\X$ and $\Y$ have the resolution property. Let $E\in K^0(\X)$, then 
$$\tch(f_*E)=If_*(\tch(E)\ttd(T_f)),$$
where $f_*$ is the K-theoretic pushforward and $If: I\X\to I\Y$ is the map induced by $f$.
\end{theorem}
\begin{remark}
The cohomological pushforward $If_*$ of  a non-representable morphism is defined by passing to a finite scheme cover of $I\X$, see \cite{Kr}.
\end{remark}


Restricting to the distinguished component $\Y\subset I\Y$, we obtain
\begin{corollary}\label{rr}
$$ch(f_*E)=If_*(\tch(E)\ttd(T_f)|_{If^{-1}(\Y)}).$$
\end{corollary}

\section{Properties of Virtual Bundles}\label{vbdle}
In this appendix we discuss some properties of the virtual bundle $F_{g,n,d}$. First note that the fact that $F_{g,n,d}$ is well-defined can be seen by factoring $f$ as in (\ref{lciq}) (which follows from the construction of the universal family in \cite{AGOT}). Note that $f$
is perfect, and resolution property implies that the $K$-theory of vector bundles coincides with the $K$-theory of perfect complexes.

We study how $F_{g,n,d}$ behaves under pulling back by the maps $f:\Mbar_{g,n+1}(\X,d)'\to\Mbar_{g,n}(\X,d)$, $s_j:\D_j\to \Mbar_{g,n+1}(\X,d)'$, and $i: \Z\to \Mbar_{g,n+1}(\X,d)'$. Let $\iota_{red}:\tilde{\Z}^{red}\to\Mbar_{g,n+1}(\X,d)'$ be the composition of double covering of $\Z^{red}$ and the inclusion into $\Mbar_{g,n+1}(\X,d)'$. Similarly we can define $\iota_{irr}:\tilde{\Z}^{irr}\to\Mbar_{g,n+1}(\X,d)$. By the definition of $\Z$ it is the universal gerbe at node over $f(\Z)\subset \Mbar_{g,n}(\X, d)$. According to \cite{AGV}, Proposition 5.2.1, we have  $$\tilde{\Z}^{red}=\coprod_{g_++g_-=g,n_++n_-=n,d_++d_-=d}\Mbar_{g_+,n_++1}(\X,d_+)\times_{I\X}\Mbar_{g_-,n_-+1}(\X,d_-), $$ and
$$\tilde{\Z}^{irr}=\Mbar_{g-1,n+2}(\X,d)\times_{I\X\times I\X}I\X.$$
Therefore we may view $\tilde{\Z}^{red}$ as the moduli stack which parametrizes pairs 
\begin{equation}\label{separating_node_locus_descrip}
(f_+: (\C_+, \{\Sigma_i\}_{1\leq i\leq n+}\cup\{\Sigma_+\})\to \X, f_-: (\C_-, \{\Sigma_i\}_{1\leq i\leq n_-}\cup\{\Sigma_-\})\to \X),
\end{equation}
where $[f_\pm]\in \Mbar_{g_\pm, n_\pm}(\X, d_\pm)$, such that $$[f_+|_{\Sigma_+}]= I([f_-|_{\Sigma_-}])\in I\X.$$ Here $I: I\X\to I\X$ is the involution defined in Section \ref{orbifold}

Similarly we may view $\tilde{\Z}^{irr}$ as the moduli stack which parametrizes maps 
\begin{equation}\label{non-separating_node_locus_descrip}
[f: (\C, \{\Sigma_i\}_{1\leq i\leq n}\cup\{\Sigma_+, \Sigma_-\})\to \X]\in \Mbar_{g-1, n+2}(\X, d),
\end{equation}
 such that $$[f|_{\Sigma_+}]=I([f|_{\Sigma_-}])\in I\X.$$

Let $ev_{node}: \tilde{\Z}\to I\X$ denote the evaluation map at the marked point of gluing in the description of $\Z$ above. More precisely, $ev_{node}$ is defined to map (\ref{separating_node_locus_descrip}) to $[f_+|_{\Sigma_+}]\in I\X$ and map (\ref{non-separating_node_locus_descrip}) to $[f|_{\Sigma_+}]\in I\X$. 

\begin{lemma}\label{split}
\hfill
\begin{enumerate}
\item \label{pull}$f^*F_{g,n,d}=F_{g,n+1,d}|_{\Mbar_{g,n+1}(\X,d)'}.$
\item \label{split1}$\iota_{red}^*F_{g,n+1,d}=p_+^*F_{g_+,n_++1,d_+}+p_-^*F_{g_-,n_-+1,d_-}-ev_{node}^*(q^*F)^{inv}.$
\item \label{split2}$\iota_{irr}^*F_{g,n+1,d}=F_{g-1,n+2,d}-ev_{node}^*(q^*F)^{inv}$. 
\end{enumerate}
\end{lemma}
\begin{proof}
The proofs are similar to those of the corresponding statements in \cite{CG}, \cite{C}. Let $\X=[M/G]$ be as in Assumption \ref{global_quotient_assumption}, where $M$ is a smooth quasi-projective variety and $G$ is a linear algebraic group. Choose a $G$-equivariant ample line bundle $L$ on $M$. The bundle $F$ corresponds to an equivariant vector bundle which we also denote by $F$. For $N$ sufficiently large we have the following exact sequence
$$
0\to Ker\to H^0(M,F\otimes L^N)\to F\otimes L^N\to 0.
$$
Tensoring with $L^{-N}$ yields an exact sequence
$$
0\to Ker\otimes L^{-N}\to H^0(M, F\otimes L^N)\otimes L^{-N}\to F\to 0.
$$
Let $A=H^0(M,F\otimes L^N)\otimes L^{-N}$ and $B=Ker\otimes L^{-N}$. These two bundles
induce two vector bundles on $\X$ which we denote by $\A$ and $\B$ respectively.
The above exact sequence implies that $F_{g,n,d}=\A_{g,n,d}-\B_{g,n,d}$.

If $d\neq 0$, then $R^0 f_*ev_{n+1}^*\A$ and $R^0 f_*ev_{n+1}^*\B$ both vanish for $N$ sufficiently large, and 
$-\A_{g,n,d},\,\,-\B_{g,n,d}$ are vector bundles.

We verify (\ref{split1}) for $\A_{g,n,d}$. Let $T$ be a scheme. Let  $$((f_+:\C_+\to\X),(f_-:\C_-\to\X))$$
 be a $T$-valued point of $$\overline{\M}_{g_+,n_++1}(\X,d_+)\times_{I\X}\overline{\M}_{g_-,n_-+1}(\X,d_-),$$
and $f:\C\to\X$ the stable map obtained by gluing. Denote by $\mathfrak{t}: \C\to T$, $\mathfrak{t}_\pm: \C_\pm\to T$ the structure maps, by $\nu: \C_+\cup \C_-\to \C$ the gluing morphism, and by $\Theta_{node}\subset \C$ the locus of the node formed by gluing. The restriction of $-\iota_{red}^*\A_{g,n+1,d}$
to the $T$-valued point $(f:\C\to\X)$ is $$R^1\mathfrak{t}_*f^*\A\simeq (R^0\mathfrak{t}_*(f^*\A^{\vee}\otimes\omega_{\C}))^{\vee}.$$
The restriction of $-p_\pm^*\A_{g_\pm,n_\pm+1,d_\pm}$ to the $T$-valued point $(f_\pm:\C_\pm\to\X)$ is $$R^1\mathfrak{t}_{\pm*}f_\pm^*\A\simeq (R^0\mathfrak{t}_{\pm*}(f_\pm^*\A^{\vee}\otimes\omega_{\C_\pm}))^{\vee}.$$ The relative dualizing sheaves of $\C,\C_+,\C_-$ are easily seen to fit into the following exact sequence,
$$0\to \omega_{\C/S}\to \nu_*(\omega_{\C_+/S}\oplus \omega_{\C_-/S})\to \sO_{\Theta_{node}}\to 0.$$
Tensoring by $f^*\A^\vee$ and applying $\mathfrak{t}_*$ give the following exact sequence:
\begin{equation*}
0\to R^0\mathfrak{t}_*(f^*\A^\vee\otimes \omega_{\C/T})\to R^0\mathfrak{t}_{+*}(f_+^*\A^\vee\otimes\omega_{\C_+/T})\oplus R^0\mathfrak{t}_{-*}(f_-^*\A^\vee\otimes \omega_{\C_-/T})\to R^0\mathfrak{t}_*(f^*\A^\vee\otimes \sO_{\Theta_{node}})\to 0.
\end{equation*}
 Note that $R^0\mathfrak{t}_*(f^*\A^\vee\otimes \sO_{\Theta_{node}})$ is the sheave of sections of $f^*\A^\vee$ which are {\em invariant} under the action of the stabilizer group of the node. Therefore $R^0\mathfrak{t}_*(f^*\A^\vee\otimes \sO_{\Theta_{node}})$ is the restriction to the $T$-valued point $((f_+:\C_+\to\X),(f_-:\C_-\to\X))$ of $ev_{node}^*((q^*\A^\vee)^{inv})$ (c.f. the proof of Lemma \ref{bundle}). Dualizing this sequence then proves (\ref{split1}) for $\A_{g,n,d}$. We can prove it for $\B_{g,n,d}$ in the same way. (\ref{split1}) thus hold for $F_{g,n,d}$ since $F_{g,n,d}=\A_{g,n,d}-\B_{g,n,d}$. 

If $d=0$, then $R^0f_*ev_{n+1}^*F$ is a trivial bundle and $R^1f_*ev_{n+1}^*F$ is a vector bundle. The same argument can be applied to this case.

(\ref{pull}) and (\ref{split2}) can be proved by a similar approach, we omit the details.
\end{proof}


\section{An Example of Quantized Operator}\label{quantized_operator}
In this Appendix we reproduce the calculation in \cite{C}, Example 1.3.3.1. 

Let $A=Bz^m$ be an infinitesimal symplectic transformation of $\sH$. Here $B: H^*(I\X)\to H^*(I\X)$ is a linear transformation. We write $B$ as a matrix $(B^\alpha_\beta)$ using the basis $\{\phi_\alpha\}$ of $H^*(I\X)$. Put $g_{\alpha\beta}=(\phi_\alpha,\phi_\beta)_{orb}$ and let $g^{\alpha\beta}$ denotes the matrix entry of the matrix inverse to $(g_{\alpha\beta})$. Then define $B_{\alpha\beta}=g_{\alpha\gamma}B^\gamma_\beta$ and $B^{\alpha\beta}=B^\alpha_\gamma g^{\gamma\beta}$.

A direct calculation shows that, 
$$\widehat{A}=\frac{1}{2\hbar}\sum_{0\leq k\leq -m-1} (-1)^{k+m}B_{\alpha\beta}q_k^\beta q_{-1-k-m}^\alpha-\sum_{k\geq -m}  B^\alpha_\beta q_k^\beta \frac{\partial}{\partial q_{k+m}^\alpha}, \quad \text{if } m< 0,$$
and 
$$\widehat{A}=-\sum_{k\geq 0} B^\alpha_\beta q_k^\beta \frac{\partial}{\partial q_{k+m}^\alpha}+\frac{\hbar}{2}\sum_{0\leq k\leq m-1} (-1)^k B^{\alpha\beta}\frac{\partial}{\partial q_k^\beta}\frac{\partial}{\partial q_{m-1-k}^\alpha},\quad \text{if } m> 0.$$

For $m=0$ we have $\widehat{A}=\sum_{k\geq 0} B^\alpha_\beta q_k^\beta \frac{\partial}{\partial q_{k}^\alpha}$.

We calculate 
\begin{equation*}
\begin{split}
&\left(\sum_k B^\alpha_\beta q_k^\beta \frac{\partial}{\partial q_{k+m}^\alpha}\right)\bq=\left(\sum_k B^\alpha_\beta q_k^\beta \frac{\partial}{\partial q_{k+m}^\alpha}\right)\left(\sum_l q_l^\gamma\phi_\gamma z^l\right)\\
&= \left[\sum_kB^\alpha_\beta q_k^\beta \phi_\alpha z^{k+m}\right]_+=[A\bq]_+.
\end{split}
\end{equation*}
This explains the appearance of (\ref{pqterm}).

Now suppose $m> 0$. We want to explain the double derivative terms in $\widehat{A}$ above, following \cite{C}, Example 1.3.3.1. Observe that the double derivative $$\frac{\partial}{\partial q_k^\beta}\frac{\partial}{\partial q_{m-1-k}^\alpha}$$ is the bivector field corresponding to $$\phi_\beta\bpsi_+^k\otimes \phi_\alpha\bpsi_-^{m-1-k} \simeq \sH_+\otimes\sH_+ \quad(\text{identifying }\bpsi_+, \bpsi_- \text{ with }z).$$
Note that for $m\geq 1$, we have $$\sum_{0\leq k\leq m-1}(-1)^k\bpsi_+^k\bpsi_-^{m-1-k}=\frac{\bpsi_+^m+(-1)^{m-1}\bpsi_-^m}{\bpsi_++\bpsi_-}.$$

Thus the term $\sum_{0\leq k\leq m-1} (-1)^k B^{\alpha\beta}\frac{\partial}{\partial q_k^\beta}\frac{\partial}{\partial q_{m-1-k}^\alpha}$ can be interpreted as the bivector field corresponding to $$\frac{B\bpsi_+^m+(-1)^{m-1}B\bpsi_-^m}{\bpsi_++\bpsi_-}.$$ 
This explains the appearance of (\ref{codim2}).

\section{Cocycle calculation}\label{cocycle_calculation}
In this appendix we calculate the cocycle (\ref{cocycle_GRR}). We begin with a lemma.

\begin{lemma}
Let $\Y$ be a smooth proper Deligne-Mumford stack. Denote by $q: I\Y\to \Y$ the natural projection. Let $A: H^*(\Y,\com)\to H^*(\Y,\com)$ be a linear operator defined by a class $a\in H^*(\Y,\com)$, i.e. $A(\gamma)=a\cdot \gamma$. Then 
$$\text{str}(A)=\int_{I\Y}q^*(a)\wedge e(T_{I\Y}).$$
\end{lemma}
\begin{proof}
Write the class $a$ as a sum of its degree zero part and positive degree part: $a=a_01+a'$ where $a'\in H^{>0}(\Y,\com)$. Since $H^*(\Y,\com)$ is a graded ring, the operator of multiplication by a positive degree element of $H^*(\Y,\com)$ has super-trace $0$. So $str(A)=str(a_01\cdot)=str(a_0 id)$. 

We find that
\begin{equation*}
\begin{split}
str(id)&=\chi(I\Y)\quad \text{by the Lefschetz trace formula (see e.g. \cite{Beh_coh})}\\
&=\int_{I\Y}e(T_{I\Y})\quad \text{by Gauss-Bonnet (see e.g. \cite{T2}, Corollaire 3.44)}.
\end{split}
\end{equation*}
Since $q^*a'\wedge e(T_{I\Y})=0$, we have 
$$\int_{I\Y}q^*(a)\wedge e(T_{I\Y})=\int_{I\Y}q^*(a_01)\wedge e(T_{I\Y})=str(a_0 id).$$
\end{proof}

(\ref{cocycle_GRR}) is obtained by applying this Lemma to each component $\X_i$ of $I\X$, and use the definition of {\em double inertia stack} $II\X:=I(I\X)=\cup_i I\X_i$. We denote the projection by $Iq: II\X\to I\X$.

\section{Proof of (TRR)}\label{pf_of_TRR}
In this Appendix we give a proof of the topological recursion relations (TRR) in genus $0$. In this proof we will use the moduli stack $\K_{g,n}(\X,d)$ instead of $\Mbar_{g,n}(\X,d)$. This is because the proof involves splitting nodal twisted curves along a node, and it is easier to express this using the stack $\K_{g,n}(\X,d)$. As pointed out in \cite{AGV2}, Section 6.1.3, orbifold Gromov-Witten invariants defined using $\K_{g,n}(\X,d)$ agree with those defined using $\Mbar_{g,n}(\X,d)$. We refer to \cite{AGV2} for properties of $\K_{g,n}(\X,d)$ used here.

Our proof is adopted from \cite{Ma}, Section VI.6.6. Let $\mathfrak{M}^{tw}_{0, 3+k}$ be the (Artin) stack of twisted curves of genus $0$ with $3+k$ marked gerbes (not trivialized) and denote by $p: \K_{0, 3+k}(\X,d)\to \mathfrak{M}^{tw}_{0,3+k}$ the forgetful morphism. For each partition $\{4,5,...,3+k\}=A\coprod B$ with $A,B$ nonempty we consider the stack $\mathfrak{D}^{tw}(0; \{1\}\cup A|0; \{2,3\}\cup B)$ defined in \cite{AGV2}, Section 5.1. Put $$\mathfrak{D}^{tw}:=\coprod_{A, B; A\coprod B=\{4,5,...,3+k\}}\mathfrak{D}^{tw}(0; \{1\}\cup A|0; \{2,3\}\cup B).$$ There is a natural gluing map $gl: \mathfrak{D}^{tw}\to \mathfrak{M}^{tw}_{0,3+k}$ as defined in \cite{AGV2}, Proposition 5.1.3. Form the following cartesian diagram 
$$\begin{CD}
\D(\X)@>\mu>> \K_{0,3+k}(\X,d)\\
@V{}VV @V{p}VV\\
\mathfrak{D}^{tw}@>gl>>\mathfrak{M}^{tw}_{0, 3+k}.
\end{CD}$$
Let $\D^{tw}\subset \mathfrak{M}^{tw}_{0, 3+k}$ denote the image of $\mathfrak{D}^{tw}$ under the map $gl$. Consider the forgetful maps $\mathfrak{M}^{tw}_{0, 3+k}\to \mathfrak{M}_{0, 3+k}\to \Mbar_{0,3}$, where the first map takes a twisted curve to its coarse curve, and the second map forgets all but the first three marked points and stabilizes the curves. Let $L_1$ be the line bundle over $\mathfrak{M}^{tw}_{0, 3+k}$ obtained by pulling back the first universal cotangent line bundle over $\mathfrak{M}_{0, 3+k}$, and $L_1'$ the line bundle over $\mathfrak{M}^{tw}_{0, 3+k}$ obtained by pulling back the first universal cotangent line bundle over $\Mbar_{0,3}$. (We slightly abuse notations here.) It is not hard to see that there is an exact sequence 
$$0\to L_1'\to L_1 \to \sO_{\D^{tw}}\to 0.$$
A standard intersection theory result (see e.g. \cite{Ma}, Chapter VI, equation (6.19)) shows that for any cycle class $\alpha$ on $\K_{0,3+k}(\X,d)$ we have 
$$ c_1(p^*L_1)\cap \alpha=c_1(p^*L_1')\cap \alpha +\mu_*gl^!\alpha.$$
Take $\alpha=[\K_{0,3+k}(\X,d)]^{vir}$ and use the fact that $c_1(L_1')=0$ (because $\Mbar_{0,3}$ is a point), we get 
\begin{equation}\label{trading_psi}
\bpsi_1 \cap [\K_{0,3+k}(\X,d)]^{vir}=\mu_*(gl^![\K_{0,3+k}(\X,d)]^{vir}).
\end{equation}
According to \cite{AGV2}, Proposition 5.2.2, we have $$\mathfrak{D}^{tw}(0; \{1\}\cup A|0; \{2,3\}\cup B)\times_{\mathfrak{M}^{tw}_{0, 3+k}}\K_{0,3+k}(\X,d)\simeq \coprod_{d_1+d_2=d}\K_{0, \{1\}\cup A\cup\bullet}(\X,d_1)\times_{\bar{I}\X}\K_{0, \{2,3\}\cup B\cup \star}(\X,d_2),$$
where $\bar{I}\X$ is the {\em rigidified} inertia stack of $\X$ (see \cite{AGV2}, Section 3.4). The diagonal map $\delta :\bar{I}\X\to \bar{I}\X\times \bar{I}\X$ fits into the following cartesian diagram
$$\begin{CD}
\K_{0, \{1\}\cup A\cup\bullet}(\X,d_1)\times_{\bar{I}\X}\K_{0, \{2,3\}\cup B\cup \star}(\X,d_2)@>>> \K_{0, \{1\}\cup A\cup\bullet}(\X,d_1)\times\K_{0, \{2,3\}\cup B\cup \star}(\X,d_2)\\
@V{ev_{node}}VV @V{ev_\bullet\times \ev_\star}VV\\
\bar{I}\X@>\delta>>\bar{I}\X\times \bar{I}\X.
\end{CD}$$
By the splitting result (\cite{AGV2}, Proposition 5.3.1), we get
\begin{equation}\label{splitting_[K]_vir} 
gl^![\K_{0,3+k}(\X,d)]^{vir}=\sum_{A\coprod B=\{4,5,...,3+k\}; d_1+d_2=d}\delta^!([\K_{0, \{1\}\cup A\cup\bullet}(\X,d_1)]^{vir}\times[\K_{0, \{2,3\}\cup B\cup \star}(\X,d_2)]^{vir}).
\end{equation}
We may apply (\ref{splitting_[K]_vir}) to (\ref{trading_psi}) and view the resulting equality in homology via cycle map. Integrate the resulting equality against $\gamma=\phi_{\alpha_1}\bpsi_1^{k_1}\prod_{i=2}^{3+k}\phi_{\alpha_i}\bpsi_i^{k_i}$ and  use an identification of $H^*(I\X)$ and $H^*(\bar{I}\X)$ (c.f. \cite{AGV2} section 6.1.3), we get 
$$
\<\phi_{\alpha_1}\bpsi_1^{k_1+1}\prod_{i=2}^{3+k}\phi_{\alpha_i}\bpsi_i^{k_i}\>_{0, 3+k,d} \\
$$
$$
=\sum_{A\coprod B=\{4, ...,3+k\}, d=d_1+d_2}\sum_\alpha \pm\<\phi_{\alpha_1}\bpsi_1^{k_1},\prod_{i\in A}\phi_{\alpha_i}\bpsi_i^{k_i} ,\phi_a\>_{0, |A|+2,d_1}\<\phi^a, \phi_{\alpha_2}\bpsi_2^{k_2}, \phi_{\alpha_3}\bpsi_3^{k_3},\prod_{i\in B}\phi_{\alpha_i}\bpsi_i^{k_i}\>_{0,|B|+3,d_2}.
$$
Here the sign come from the possibly different ordering of odd cohomology classes between the left and right sides. (TRR) follows as this is the equality of coefficients of the corresponding terms on the left and right sides of (TRR).

\end{document}